\documentclass[oneside,11pt,reqno]{amsart}
\usepackage[utf8]{inputenc}
\usepackage{stmaryrd}
\usepackage{bbm}
\usepackage[dvipsnames]{xcolor}
\usepackage{mathrsfs}
\usepackage{enumerate}
\usepackage{centernot}
\usepackage{latexsym,amsxtra}
\usepackage{dsfont}
\usepackage{slashed}
\usepackage[all]{xy}
\usepackage{amscd,graphics}
\usepackage{amsmath,amsfonts,amsthm,amssymb}
\usepackage{latexsym,amsmath}
\usepackage{graphicx,psfrag}
\usepackage{mathabx}
\usepackage{hyperref}
\usepackage{fancyhdr}

\usepackage{cleveref}
\usepackage{comment} 
\usepackage{braket}
\usepackage{caption}
\usepackage{enumitem}
\usepackage{kantlipsum}
\calclayout

\theoremstyle{plain}
\newtheorem{thm}{Theorem}[section]
\crefname{thm}{Theorem}{Theorems}
\Crefname{thm}{Theorem}{Theorems}
\newtheorem{pro}[thm]{Proposition}
\crefname{pro}{Proposition}{Propositions}
\Crefname{pro}{Proposition}{Propositions}
\newtheorem{lem}[thm]{Lemma}
\crefname{lem}{Lemma}{Lemmas}
\Crefname{lem}{Lemma}{Lemmas}
\newtheorem{cor}[thm]{Corollary}
\crefname{cor}{Corollary}{Corollaries}
\Crefname{cor}{Corollary}{Corollaries}
\newtheorem{conj}{Conjecture}
\crefname{conj}{Conjecture}{Conjectures}
\Crefname{conj}{Conjecture}{Conjectures}

\crefname{cons}{Construction}{Constructions}
\Crefname{cons}{Construction}{Constructions}

\crefname{claim}{Claim}{Claims}
\Crefname{claim}{Claim}{Claims}

\crefname{property}{Property}{Properties}
\Crefname{property}{Property}{Properties}

\crefname{problem}{Problem}{Problems}
\Crefname{problem}{Problem}{Problems}

\theoremstyle{definition}
\newtheorem{defi}[thm]{Definition}
\crefname{defi}{Definition}{Definitions}
\Crefname{defi}{Definition}{Definitions}
\newtheorem{nota}[thm]{Notation}
\crefname{nota}{Notation}{Notations}
\Crefname{nota}{Notation}{Notations}
\newtheorem {convention}[thm]{Convention}
\crefname{convention}{Convention}{Conventions}
\Crefname{convention}{Convention}{Conventions}

\crefname{cond}{Condition}{Conditions}
\Crefname{cond}{Condition}{Conditions}
\newtheorem{assum}[thm]{Assumption}
\crefname{assum}{Assumption}{Assumptions}
\Crefname{assum}{Assumption}{Assumptions}

\theoremstyle{remark}
\newtheorem{rmk}[thm]{Remark}
\crefname{rmk}{Remark}{Remarks}
\Crefname{rmk}{Remark}{Remarks}
\newtheorem{ex}[thm]{Example}
\crefname{ex}{Example}{Examples}
\Crefname{ex}{Example}{Examples}

\crefname{ques}{Question}{Questions}
\Crefname{ques}{Question}{Questions}

\crefname{section}{Section}{Sections}
\Crefname{section}{Section}{Sections}
\crefname{subsection}{Subsection}{Subsections}
\Crefname{subsection}{Subsection}{Subsections}
\crefname{figure}{Figure}{Figures}
\Crefname{figure}{Figure}{Figures}
\newcommand{\Emb}{\textnormal{Emb}^{\textnormal{lc}}}

\newcommand{\Diff}{\textnormal{Diff}}
\newcommand{\BDiff}{\textnormal{BDiff}}
\newcommand{\TDiff}{\textnormal{TDiff}}
\newcommand{\BTDiff}{\textnormal{BTDiff}}
\newcommand{\MCG}{\textnormal{MCG}^{+}}
\newcommand{\Homeo}{\textnormal{Homeo}}
\newcommand{\BHomeo}{\textnormal{BHomeo}}

\newcommand{\spinc}{\textnormal{spin}^{c}}
\newcommand{\ab}{\textnormal{ab}}

\newcommand{\SW}{\mathrm{SW}}
\newcommand{\SWbb}{\mathbb{SW}}
\newcommand{\SWcal}{\mathcal{SW}}
\newcommand{\SWbbtot}{\mathbb{SW}_{\mathrm{tot}}}
\newcommand{\SWcaltot}{\mathcal{SW}_{\mathrm{tot}}}
\newcommand{\SWbbhalftot}{\mathbb{SW}_{\mathrm{half\mathchar`- tot}}}
\newcommand{\SWcalhalftot}{\mathcal{SW}_{\mathrm{half\mathchar`- tot}}}
\newcommand{\EDiff}{\textnormal{EDiff}}
\newcommand{\id}{\textnormal{id}}
\newcommand{\Z}{\mathbb{Z}}

\newcommand{\R}{\mathbb{R}}

\newcommand{\Q}{\mathbb{Q}}

\newcommand{\CP}{\mathbb{CP}}
\newcommand{\fraks}{\mathfrak{s}}
\newcommand{\frakt}{\mathfrak{t}}
\newcommand{\frakM}{\mathfrak{M}}
\newcommand{\scrA}{\mathscr{A}}
\newcommand{\scrB}{\mathscr{B}}
\newcommand{\scrC}{\mathscr{C}}
\newcommand{\scrD}{\mathscr{D}}
\newcommand{\scrX}{\mathscr{X}}
\newcommand{\scrE}{\mathscr{E}}
\newcommand{\scrH}{\mathscr{H}}
\newcommand{\scrN}{\mathscr{N}}
\newcommand{\scrG}{\mathscr{G}}

\newcommand{\scrT}{\mathscr{T}}
\newcommand{\scrR}{\mathscr{R}}
\newcommand{\scrM}{\mathscr{M}}
\newcommand{\calM}{\mathcal{M}}
\newcommand{\calU}{\mathcal{U}}
\newcommand{\calL}{\mathcal{L}}

\newcommand{\cpt}{\mathrm{cpt}}
\newcommand{\ind}{\mathrm{ind}}
\newcommand{\Conj}{\mathrm{Conj}}
\DeclareMathOperator{\Aut}{Aut}

\DeclareMathOperator{\Map}{Map}
\DeclareMathOperator{\Int}{Int}
\DeclareMathOperator{\sign}{sign}
\newcommand{\Spinc}{\mathrm{Spin}^{c}}
\newcommand{\Spincns}{\mathrm{Spin}^{c}}
\newcommand{\circPi}{\mathring{\Pi}}
\newcommand{\del}{\partial}

\newcommand{\texts}{\textsf{s}}

\title[Homological instability for smooth 4-manifolds]{Homological instability for moduli spaces of smooth 4-manifolds}

\author{Hokuto Konno}
\address{Graduate School of Mathematical Sciences, the University of Tokyo, 3-8-1 Komaba, Meguro, Tokyo 153-8914, Japan \\and\\
RIKan iTHEMS, Wako, Saitama 351-0198, Japan}
\email{konno@ms.u-tokyo.ac.jp}

\author{Jianfeng Lin}
\address{Yau Mathematical Sciences Center\\ Tsinghua University\\ Beijing\\ China.}
\email{linjian5477@mail.tsinghua.edu.cn}

\begin{document}

\maketitle

\begin{abstract}
We prove that homological stability fails for the moduli space of any simply-connected closed smooth 4-manifold in any degree of homology, unlike what happens in all even dimensions $\neq 4$. We detect also the homological discrepancy between various moduli spaces, such as topological and smooth moduli spaces of 4-manifolds, and moduli spaces of 4-manifolds with positive scalar curvature metrics. Some of these discrepancies are special to dimension $4$. To prove these results, we use the Seiberg--Witten equations to construct a new characteristic class of families of 4-manifolds, which is unstable and detects the difference between the smooth and topological categories in dimension 4.
\end{abstract}


\section{Introduction}
\label{section Introduction}

\subsection{Background and setup}

Given a smooth manifold $M$, the classifying space $\BDiff(M)$ (also called the moduli space of $M$) classifies smooth fiber bundles with fiber $M$. Cohomology classes of $\BDiff(M)$ one-to-one correspond to characteristic classes for such bundles. For this reason, understanding the (co)homology group of $\BDiff(M)$ has been an important topic in manifold theory. A celebrated result of Harer~\cite{Harer85} proved homological stability for mapping class groups of oriented surfaces. This result can be interpreted in terms of
moduli spaces of Riemann surfaces. More recently, Galatius and Randal-Williams \cite{Galatius2018} established analogous homological stability for moduli spaces of manifolds of even dimension $\geq 6$, which is one of the highlights of the study of manifolds in the last decade. Similar homological stability results have also been established for moduli spaces of manifolds of odd dimension in \cite{Perlmutter16,Perlmutter18,Perlmutter16Linking,Krannich2019}  for $\dim \geq 9$ and in \cite{Lam15} for $\dim=3$, and also for mapping class groups of 3-manifolds in  \cite{Hatcher10}. As a natural question, one may ask what happens in the remaining dimensions $4, 5$ and $7$.
The main purpose of this paper is to answer this question in the only unresolved even dimension 4 by establishing a phenomenon in the opposite direction of homological stability, which we call {\it homological instability}, in broad generality. The second purpose is to establish comparison results between various moduli spaces of 4-manifolds and 4-manifolds with metrics. In particular, we will establish certain exotic phenomena that only appear in dimension 4.

To state our results, let us describe the standard setup to discuss homological stability for moduli spaces of manifolds.
Let $W$ be a compact smooth $2n$-dimensional manifold with non-empty boundary.
Let $\Diff_\del(W)$ denote the group of diffeomorphisms that are the identity near $\del W$.
We equip $\Diff_\del(W)$ with the $C^\infty$-topology.
One may define the {\it stabilization map}
\begin{align*}
s : \Diff_\del(W) \to \Diff_\del(W\#S^n\times S^n)
\end{align*}
as follows.
Form the connected sum of $W$ with $S^n \times S^n$ by letting $K = ([0,1] \times \del W)\#S^n\times S^n$ be the inner connected sum, and setting
$W \# S^n\times S^n = W \cup_{\del W=\{0\} \times \del W} K$.
Then the map $s$ is defined by extending diffeomorphisms lying in $\Diff_\del(W)$ by the identity of $S^n \times S^n$.
This induces the map, called the stabilization map, between homologies of moduli spaces
\[
s_\ast : H_k(\BDiff_\del(W);\Z) \to H_k(\BDiff_\del(W\#S^n\times S^n);\Z)
\]
for every $k \geq 0$. 
If $W$ is obvious from the context, for $N\geq 0$, we use 
\[
s_{N,\ast} : H_k(\BDiff_\del(W\#N S^n \times S^n);\Z) \to H_k(\BDiff_\del(W\#(N+1) S^n \times S^n);\Z)
\]
to denote the stabilization map for $W\#N S^n \times S^n$.

The aforementioned results of Harer~\cite{Harer85} for $2n=2$ and of Galatius and Randal-Williams~\cite{Galatius2018} for $2n \geq 6$ state that, for any simply-connected $W$ and any $k$, the stabilization map $s_{N,\ast}$ is an isomorphism for any $N>0$ large enough relative to $k$. 
As a consequence, the direct system \[\{(H_k(\BDiff_{\del}(W\#NS^n\times S^n);\Z), s_{N,\ast})\}_{N=0}^\infty\]
stabilizes after finitely many terms. This stability phenomenon is fundamental in the study of the homology group of moduli spaces. It provides a \emph{stable range} in which one can identify $H_{k}(\BDiff_{\del}(W\# N S^n\times S^n))$ with the stable homology group \[\lim_{N \to +\infty}H_k(\BDiff_{\del}(W\# N S^n\times S^n)),\] which can be computed purely in terms of homotopy theory \cite{MW07,Galatius17}. (Note that this method of computing stable homology group is valid in all even dimensions including 4.)

A closely related object of interest is the diffeomorphism group with the discrete topology, denoted by $\Diff_{\partial}(W)^{\delta}$. The Eilenberg--MacLane space $\BDiff_\del(W)^{\delta}$ classifies smooth fiber bundles with flat structures. 
Homology of $\BDiff_\del(W)^{\delta}$ is the \emph{group homology} of $\Diff_\del(W)$ and it reflects the purely algebraic (rather than topological) structure of $\Diff_{\partial}(W)$.

One can similarly define the stabilization map 
\[
s^{\delta}_{N,\ast}: H_k(\BDiff_\del(W\#N S^n \times S^n)^{\delta};\Z) \to H_k(\BDiff_\del(W\#(N+1) S^n \times S^n)^{\delta};\Z).
\]
When $W$ is a punctured $g( S^{n}\times S^{n})$ for $g>0$ with $2n\neq 4$, Nariman~\cite{Nariman17a,Nariman17} proved that $s^{\delta}_{N,\ast}$ is an isomorphism for any $N$ large enough relative to $k$.

\subsection{Homological instability}
Now we present the first main result of this paper.
Let us focus on when $2n=4$ and $W$ is a punctured 4-manifold.
For a closed smooth 4-manifold $X$, choose an embedded 4-disk $D^4$ in $X$ and set $\mathring{X} = X \setminus \Int(D^4)$.
We begin with a concise version of the theorem, which implies that, unlike dimensions $\neq 4$, homological stability fails for the diffeomorphism group of any simply-connected closed 4-manifolds, with either the $C^\infty$-topology or the discrete topology, in any positive degree of homology:
\begin{thm}[\bf{Homological instability: concise version}]
\label{thm: simpler main thm}
Let $X$ be a simply-connected closed smooth 4-manifold and let $k>0$.
Then there are infinitely many positive integers $N$ such that the stabilization maps
\[
s_{N,\ast} : H_k(\BDiff_\del(\mathring{X}\#N S^2 \times S^2);\Z) \to H_k(\BDiff_\del(\mathring{X}\#(N+1) S^2 \times S^2);\Z)
\]
are not isomorphic.
An analogous result holds also for $s_{N,\ast}^\delta$.
\end{thm}
In fact, we shall prove a much more precise result. To state it, we let $\Homeo_\del(\mathring{X})$ denote the group of homeomorphisms restricting to the identity near $\del \mathring{X}$ and let $i : \Diff_\del(\mathring{X}) \hookrightarrow \Homeo_\del(\mathring{X})$ be the inclusion maps.
We call an element of 
\[
\ker(i_\ast : H_\ast(\BDiff_\del(\mathring{X});\Z) \to H_\ast(\BHomeo_\del(\mathring{X});\Z))
\]
a {\it topologically trivial} homology class.
The following is a more precise version of the main result of this paper. It establishes a sequence of 2-torsion, topologically trivial unstable homology classes.
Also, the ``non-isomorphic" result for infinitely many stabilization maps in \cref{thm: simpler main thm} will be refined to non-injectivity and non-surjectivity results for infinitely many stabilization maps:
\begin{thm}[\bf{Homological instability: precise version}]
\label{thm: main cal}
Let $X$ be a simply-connected closed smooth 4-manifold and let $k>0$.
Then there is a sequence of positive integers $N_1<N_2<\cdots \to +\infty$ and a sequence of nonzero, 2-torsion elements \[\alpha_{i}\in H_k(\BDiff_{\del}(\mathring{X}\#N_iS^2 \times S^2);\Z)\] that satisfy the following properties.
\begin{description}
\item[(a) {[Topologically trivial]}] $\alpha_{i}$ is topologically trivial. 
\item[(b) {[Non-injective]}] $\alpha_{i}$ belongs to the kernel of the stabilization map $s_{N_{i},*}$. In particular, the map $s_{N_{i},*}$ is not injective.
\item[(c) {[Non-surjective]}] $\alpha_{i}$ does not belong to the image of the composition $$s_{N_{i}-1,*}\circ\cdots\circ s_{N_{i}-k-1,*}.$$ In particular, $s_{N_{i}-l,*}$ is not surjective for some $l\in \{1,\cdots,k+1\}$.
\item[(d) {[Discrete]}] $\alpha_{i}$ equals the image of a 2-torsion element 
$$
\alpha^{\delta}_{i}\in  H_k(\BDiff_{\del}(\mathring{X}\#N_i S^2 \times S^2)^{\delta};\Z)
$$
under the forgetful map
\[
H_k(\BDiff_{\del}(\mathring{X}\#N_iS^2 \times S^2)^{\delta};\Z)
\rightarrow H_k(\BDiff_{\del}(\mathring{X}\#N_iS^2 \times S^2);\Z).
\]
Furthermore, $\alpha^{\delta}_{i}$ belongs to the kernel of $s^{\delta}_{N_{i},*}$ and does not belong to the image of the composition $s^{\delta}_{N_{i}-1,*}\circ\cdots\circ s^{\delta}_{N_{i}-k-1,*}$.
\end{description}
\end{thm}

The above instability result shall be proven by using a new characteristic class of families of  4-manifolds, explained in \cref{subsection Intro SW char}. This characteristic class shall turn out to detect subtle differences between various moduli spaces of 4-manifolds, which is of independent interest.
From now on, we shall describe comparisons between moduli spaces derived from this characteristic class.




\subsection{Smooth and topological moduli spaces}

\cref{thm: main cal}~(a) reflects the homological discrepancy of the smooth and topological moduli spaces.
In fact, we prove a refinement which reflects an exotic phenomenon that is special to dimension $4$. The result holds both for a closed 4-manifold $X$ and the punctured manifold $\mathring{X}$. For simplicity, we focus on the former. We let $\Diff^+(X)$ and $\Homeo^+(X)$ denote the groups of orientation-preserving diffeomorphisms and homeomorphisms. As we mentioned, the spaces $\BHomeo^{+}(X)$ and $\BDiff^{+}(X)$  classify topological bundles and smooth bundles with fiber $X$ respectively. 
Hence we call $\BDiff^{+}(X)$ the moduli space of smooth manifold $X$ and call $\BHomeo^{+}(X)$ the moduli space of topological manifold $X$. There is actually a third moduli space $\BHomeo^{+}_{L}(X)$ as defined in \cref{defi: moduli space for formally smooth bundle}. To explain the geometric meaning of this moduli space, we recall that a microbundle is a generalization of a vector bundle, as defined by Milnor \cite{Milnormicrrobundles}. A lift of a microbundle to an actual vector bundle is called a linear structure. Following \cite{KK20}, a formally smooth bundle means a topological bundle $X\hookrightarrow E\to B$ equipped with a linear structure on its vertical tangent microbundle $\scrT(E/B)$ (see \cref{defi: vertical tangent microbundle}).  The moduli space $\BHomeo^{+}_{L}(X)$ classifies formally smooth bundles $X\hookrightarrow E\rightarrow B$ so is called the moduli space of formally smooth manifold $X$. We have canonical forgetful maps 
\[
\BDiff^{+}(X)\xrightarrow{\iota}\BHomeo^{+}_{L}(X)\xrightarrow{\tau}\BHomeo^{+}(X).
\]
The composition $\tau\circ \iota$ is homotopic to the map induced by the natural inclusion $i:\Diff^{+}(X)\hookrightarrow \Homeo^{+}(X)$.

For any manifold $X$ of $\dim\neq 4$, it follows from a celebrated theorem of Kirby--Siebenmann \cite{Kirby77} that the map $\iota : \BDiff^+(X) \to \BHomeo_L^+(X)$ is a weak homotopy equivalence to a path component of $\BHomeo^+_L(X)$ (see \cref{cor: iotawekequivalence}).
This implies, of course, that $\iota$ induces an isomorphism on homology of positive degree.
In contrast, an analogous statement in dimension 4 shall turn out to fail in any degree of homology. We shall prove:

\begin{thm}[\bf{BDiff vs. $\text{BHomeo}_{\text{L}}$}]
\label{thm: Diff Homeo sequence general}
Let $X$ be a simply-connected closed oriented smooth 4-manifold and let $k>0$.
Then there exists a sequence of positive integers $N_1<N_2<\cdots \to +\infty$ such that, for all $i$, none of the maps
\[
\iota_\ast : H_k(\BDiff^+(X\#N_i S^2 \times S^2);\Z)
\to H_k(\BHomeo_{L}^+(X\#N_iS^2 \times S^2);\Z)
\]
are injective and none of the maps 
\[
\iota^\ast : H^k(\BHomeo_{L}^+(X\#N_i S^2 \times S^2);\Z/2)
\to H^k(\BDiff^+(X\#N_iS^2 \times S^2);\Z/2)
\]
are surjective. 
\end{thm}

This \lcnamecref{thm: Diff Homeo sequence general} illustrates the special nature of moduli theory of manifolds in dimension 4, as well as the instability theorem (\cref{thm: main cal}).


\begin{rmk}[$\BDiff$ vs. $\BHomeo$]
\cref{thm: Diff Homeo sequence general} immediately implies that the induced map
\[
i_\ast : H_k(\BDiff^+(X\#N_i S^2\times S^2);\Z) \to H_k(\BHomeo^+(X\#N_i S^2\times S^2);\Z)
\]
is not injective.
This gives the first result that, in any degree $k>0$, there is a 4-manifold $M$ such that the map $i_\ast : H_k(\BDiff^+(M))
\to H_k(\BHomeo^+(M))$ is not an isomorphism. We also prove an analogous statement for the corresponding map $i_\ast^{\delta}$ with the discrete topologies on the diffeomorphism/homeomorphism groups (see \cref{thm: Diff Homeo closed 4-manifold}). 
\end{rmk}

\begin{rmk}[Two types of exotic phenomena] The failure of $\iota$ and $\tau$ being weak homotopy equivalences reflects two types of exotic phenomena in dimension 4. As suggested in \cite{KK20}, it may be fruitful to study them separately. Theorem \ref{thm: Diff Homeo sequence general} fits into this direction. It is expected that the formally smooth moduli space $\BHomeo^{+}_{L}(X)$ behaves much like the smooth moduli space of a higher dimensional manifold. While Theorem \ref{thm: main cal} shows that homological stability fails for the smooth moduli space, it may still be possible to prove homological stability for  $\BHomeo^{+}_{L}(X)$. Note that the unstable element $\alpha_{i}$ we established in  Theorem \ref{thm: main cal} becomes trivial when being mapped to the formally smooth moduli space. Also note that Watanabe \cite{Wa18} established exotic elements in $\pi_{i}(\BDiff^{+}(S^4))\otimes \mathbb{Q}$ using configuration space integrals. By \cite[Theorem 1.4]{LX23}, these elements survive in $\pi_{i}(\BHomeo_{L}^{+}(S^4))\otimes \mathbb{Q}$.
\end{rmk}

\begin{rmk}[Non-topological characteristic class]
The cohomological part of the statement of \cref{thm: Diff Homeo sequence general} implies that, in any positive degree $k$, there exists a non-topological characteristic class which is defined for \emph{all} orientable smooth fiber bundle whose fiber $X$ satisfies $b^{+}(X)>k$. It appears that characteristic classes with such properties were unknown before.
(See \cref{rem: relation with Rub and K,rmk: MMM class}.)
This also gives an explicit construction of a new unstable characteristic class for smooth 4-manifolds bundles (See \cref{subsection Intro SW char}).
\end{rmk}

\begin{rmk}
The results until here (\cref{thm: main cal,thm: Diff Homeo sequence general}) are closely related to an upcoming work by Auckly and Ruberman~\cite{AucklyRuberman}, which generalizes Ruberman's work \cite{Rub98,Rub99}.
A significant difference between \cite{AucklyRuberman} and this paper is which group of diffeomorphisms is considered.
While we consider $\Diff^+(X)$ in this paper, Auckly and Ruberman treat certain classes of proper subgroups $G(X)$ of $\Diff^+(X)$, including (subgroups of) the Torelli group 
\[
\TDiff(X) := \Set{f \in \Diff(X) | f_\ast=\id \text{ \ on\ } H_\ast(X;\Z)}.
\]
Also, for $G(X) \subset \TDiff(X)$, they detect torsion free elements in the kernels of the stabilization maps $H_*(BG(X);\Z) \to H_*(BG(X\#S^2\times S^2);\Z)$ for some $X$, while we detected 2-torsions for the whole diffeomorphism group $\Diff^+(X)$.
\end{rmk}

It is worth noting that the choice of group of diffeomorphisms can make essential differences in the homology groups of moduli spaces.
For instance, there is an example of a 4-manifold $X$ for which the kernel of the natural map
\[
H_\ast(\BTDiff(X);\Z) \to H_\ast(\BDiff^+(X);\Z)
\]
contains a subgroup isomorphic to $\Z^\infty=\bigoplus_{\Z}\Z$ (\cref{pro: Torelli vs the whole strong}).

\subsection{Abelianizations of diffeomorphism groups}

Next, we focus on homological instability in degree $k=1$, which is closely related to the abelianized mapping class group. Given $X$, we consider the smooth mapping class group $\MCG(X)=\pi_0(\Diff^+(X))$ and the topological mapping class group $\MCG_{\operatorname{Top}}(X)=\pi_0(\Homeo^+(X))$. We use $G_{\ab}$ to denote the abelianization of a group $G$. Then by definition, we have  
\[\begin{split}
H_{1}(\BDiff^+(X);\mathbb{Z})&=\MCG(X)_{\ab},\\
H_{1}(\BHomeo^+(X);\mathbb{Z})&=\MCG_{\operatorname{Top}}(X)_{\ab}.
\end{split} \]
By a series of deep results of Edwards--Kirby~\cite{Edwards71} and Thurston~\cite{Thurston74}, the unit components of $\Diff(X)$ and $\Homeo(X)$ are both simple. Therefore,  we also have 
\[\begin{split}
H_{1}(\BDiff^+(X);\mathbb{Z})&\cong \Diff^+(X)_{\ab},\\
H_{1}(\BHomeo^+(X);\mathbb{Z})&\cong \Homeo^{+}(X)_{\ab}.
\end{split}\]
Also note that when $X$ is simply-connected, the work of Freedman \cite{Fre82} and work of Quinn \cite{Q86} and Perron~\cite{P86} give the isomorphism
\begin{align}
\label{FreeQuinnPerron}
\MCG_{\operatorname{Top}}(X)\cong \Aut(Q_{X}).
\end{align}
Here $\Aut(Q_{X})$ denotes the group of automorphism on $H^{2}(X;\mathbb{Z})$ that preserves the intersection form $Q_{X}$.  Combining this isomorphism with a result of Wall \cite[Theorem 2]{Wall64D}, we see that the forgetful map
\begin{equation}\label{eq: MCG surjective}
\MCG(X\#N S^2 \times S^2)\rightarrow \MCG_{\operatorname{Top}}(X\#NS^2 \times S^2).    
\end{equation}
is always surjective when $N\geq 2$.  

Regarding the stabilization map for a closed manifold $X$, there is no well-defined stabilization map from $\BDiff^{+}(X)$ to $\BDiff^{+}(X\#S^2\times S^2)$ as there is no fixed ball to form connected sums. But as proved in \cite[Theorem 5.3]{DKPR}, there is indeed a well-defined stabilization map \[\MCG(X)\rightarrow \MCG(X\#S^2\times S^2).\]

With these results in mind, we can restate the $k=1$ case of \cref{thm: Diff Homeo sequence general,thm: main cal} as follows. 
\begin{thm}[\bf{Abelianizations of Diff and MCG}]
\label{thm: abelianisation noninjective}
Let $X$ be a simply-connected closed oriented smooth 4-manifold and let $k>0$.
Then there exists a sequence of positive integers $N_1<N_2<\cdots \to +\infty$ and a sequence of nonzero, 2-torsion elements $$\beta_{i}\in \Diff^+(X\#N_i S^2 \times S^2)_{\ab}\cong \MCG(X\#N_i S^2 \times S^2)_{\ab}$$
such that the following properties hold for all $i$. 
\begin{enumerate}[label=(\roman*)]
    \item $\beta_{i}$ can be represented by an exotic diffeomorphism (i.e. a diffeomorphism which is topologically isotopic to the identity but not smoothly so). In particular, $\beta_i$ belongs to the kernel of the map \[
\Diff^+(X\#N_i S^2 \times S^2)_{\ab}\rightarrow \Homeo^+(X\#N_iS^2 \times S^2)_{\ab}.
\]
\item $\beta_{i}$ belongs to the kernel of the map 
\[
\MCG(X\#N_i S^2 \times S^2)_{\ab}\rightarrow \MCG(X\#(N_i+1) S^2 \times S^2)_{\ab}.
\]
\item $\beta_{i}$ does not belong to the image of the map 
\[
\MCG(X\#(N_i-2) S^2 \times S^2)_{\ab}\rightarrow \MCG(X\#N_i S^2 \times S^2)_{\ab}.
\]
\end{enumerate}
\end{thm}
\cref{thm: abelianisation noninjective} also implies that for all simply-connected 4-manifolds, the degree-$1$ homological stability of smooth mapping class group fails.


\begin{ex} In dimension $\neq 4$, there are various tools to compute the mapping class groups and their abelianizations. In particular, let $W^{2n}_{g}=g(S^n\times S^n)$. Then $W^{2}_{g}$ is just a closed surface of genus $g$ and it is proved by Mumford \cite{Mumford67} and Powell \cite{Powell78} that
\[
\MCG(W^{2}_{g})_{\ab}\cong \begin{cases}
\mathbb{Z}/12 \quad &g=1,\\
\mathbb{Z}/10 \quad &g=2,\\
0 \quad &\text{otherwise}.
\end{cases}
\]
 For $n\geq 3$, Kreck  described the group $\MCG(W^{2n}_{g})$ up to extensional problems \cite{Kreck79}.
 These extension problems have been solved in many cases \cite{Galatius16,Krannich20}. As for the abelianization, Galatius and Randal-Williams~\cite{Galatius16} and Krannich~\cite{Krannich20} computed $\MCG(W^{2n}_{g})_{\ab}$ for all $n\geq 3$ and $g\geq 0$. As proved in \cite{Galatius16}, one has
$
\operatorname{Aut}(Q_{W^{4}_{g}})_{\ab}\cong \mathbb{Z}/2\oplus \mathbb{Z}/2
$
whenever $g\geq 2$. Theorem \ref{thm: abelianisation noninjective} applied to $X=S^{4}$ shows that $\MCG(W^{4}_{g})$ is strictly larger than $\mathbb{Z}/2\oplus \mathbb{Z}/2$ for infinitely many $g$. For such $g$, the kernel of the surjective map 
$$
f_{\ab}:\MCG(W^{4}_{g})_{\ab}\rightarrow \mathbb{Z}/2\oplus \mathbb{Z}/2
$$ contains a 2-torsion element that is represented by an exotic diffeomorphism on $W^{4}_{g}$ and is not stabilized from a mapping class on $W^{4}_{g-2}$.
\end{ex}

\begin{rmk}\label{rmk: twisted stabilization} All of our results also hold for twisted stabilizations (i.e. taking connected sum with $\mathbb{CP}^{2}\#\overline{\mathbb{CP}^2}$ instead of $S^2\times S^2$). 
If one replaces $S^2\times S^2$ with  $\mathbb{CP}^{2}\#\overline{\mathbb{CP}^2}$ in the statements of \cref{thm: simpler main thm}, \cref{thm: main cal}, \cref{thm: Diff Homeo sequence general} and \cref{thm: abelianisation noninjective}, the conclusions of these results still hold. See \cref{rmk: twisted stablization details} for a detailed explanation.
\end{rmk}

\subsection{Moduli spaces of manifolds with metrics}
\label{subsection:Moduli spaces of manifolds with metrics}

Until here we considered moduli spaces of manifolds.
Now we discuss moduli spaces of pairs of manifolds and metrics, in particular with positive scalar curvature.
Given an oriented smooth manifold $X$, let $\scrR(X)$ and $\scrR^+(X)$ denote the space of Riemannian metrics and the space of metrics with positive scalar curvature, respectively.
The group $\Diff^+(X)$ acts on $\scrR(X)$ by pull-back, and this action preserves  $\scrR^+(X)$.
Define $\scrM(X)$ and $\scrM^+(X)$ by 
\[
\scrM(X) = \EDiff^+(X) \times_{\Diff^+(X)} \scrR(X),\quad
\scrM^+(X) = \EDiff^+(X) \times_{\Diff^+(X)} \scrR^+(X).
\]

As we shall see in \cref{subsection Vanishing under positive scalar curvature},
$\scrM(X)$ classifies oriented fiber bundles with fiber $X$ equipped with fiberwise metrics, and similarly $\scrM^+(X)$ classifies fiber bundles with fiberwise positive scalar curvature metrics.
Let
$\iota : \scrM^+(X) \hookrightarrow \scrM(X)$
denote the injection induced from the inclusion $\scrR^+(X) \hookrightarrow \scrR(X)$.
The moduli space $\scrM(X)$ gives a model of $\BDiff^+(X)$, and $\iota$ is the same as the projection map $\scrM^+(X) \to \BDiff^+(X)$.
We shall prove the following comparison result:

\begin{thm}[\bf{Moduli spaces of manifolds with metrics vs. with psc metrics}]
\label{thm: psc main intro}
Let $X$ be a simply-connected closed oriented smooth 4-manifold and let $k>0$.
Then there exists a sequence of positive integers $N_1<N_2<\cdots \to +\infty$ such that the induced maps
\begin{align*}
\iota_\ast &: H_k(\scrM^+(X\# N_iS^2\times S^2);\Z) \to H_k(\scrM(X\# N_iS^2\times S^2);\Z)
\end{align*}
are not surjective for all $i$.
\end{thm}

We shall also prove analogous results for pointed 4-manifolds and punctured 4-manifolds in \cref{thm: psc main}, which are described in terms of observer moduli spaces $\scrM_{x_0}^+(X)$ and relative versions $\scrM_{\del}^+(\mathring{X})$.
See \cref{subsection Vanishing under positive scalar curvature} for the definition of these moduli spaces.

\cref{thm: psc main intro} gives the first result that in any degree $k$ there is a 4-manifold $X$ such that the map $\iota_\ast : H_k(\scrM^+(X))
\to H_k(\scrM(X))$ is not isomorphic.
Compared with higher dimensions, it is difficult to study the space of/moduli space of positive scalar curvature metrics on a 4-manifold because of the failure of surgery techniques.
As one of few known results,
Ruberman~\cite{Rub01} proved that $\scrR^+(X)$ are not connected for some 4-manifolds $X$.
In fact, the above non-surjectivity result (\cref{thm: psc main intro}) for $k=1$ implies that $\scrR^+(X\#N_i S^2 \times S^2)$ is disconnected (\cref{lem psc disconnected}).
Thus \cref{thm: psc main intro} can be regarded as a higher degree generalization of  Ruberman's result.

It is also worth noting that a recent work by Botvinnik and Watanabe~\cite{BW22} combined with Watanabe's work on the Kontsevich characteristic class~\cite{Wa18} implies a result for $D^4$ in a complementary direction to the above, detection of the image of (a relative version of) $\iota_\ast$.

\begin{rmk}
If $X$ is spin and $\sign(X)\neq0$, then $X\# NS^2\times S^2$ do not admit positive scalar curvature metrics for all $N\geq0$. In this case, the claim of \cref{thm: psc main} is equivalent to the nonvanishing result $H_k(\scrM(X\# N_iS^2\times S^2);\Z)\neq 0$. 
On the other hand, if $X$ is either non-spin or $\sign(X)=0$, it follows from Wall's theorem \cite{Wall64} that $X\#N S^2 \times S^2$ admit positive scalar curvature metrics for all $N \gg 0$.
\end{rmk}

\subsection{Tool: characteristic classes from Seiberg--Witten theory}\label{subsection Intro SW char}

To derive the results explained until the last subsection,
we shall construct and calculate a new characteristic class, which is based on Seiberg--Witten theory. 
The setup is as follows.
For $k>0$, let $X$ be a closed oriented smooth 4-manifold with $b^{+}(X)>k+1$.
We shall construct characteristic classes
\begin{align*}
\SWbbtot^k(X) &\in H^k(\BDiff^+(X); \Z_{\EDiff^+(X)}),\\
\SWbbhalftot^k(X) &\in H^k(\BDiff^+(X); \Z/2),
\end{align*}
where $\Z_{\EDiff^+(X)}$ is a certain local coefficient system with fiber $\Z$, determined by the monodromy action on what is called the homology orientation.
We call $\SWbbtot^k(X)$ and $\SWbbhalftot^k(X)$ the {\it $k$-th total Seiberg--Witten characteristic class} and the {\it $k$-th half-total Seiberg--Witten characteristic class} respectively.

For some computational reason explained later, $\SWbbhalftot^k(X)$ will be used rather than $\SWbbtot^k(X)$ in our applications.
Here are notable properties of $\SWbbhalftot^k(X)$.
(Similar results hold also for $\SWbbtot^k(X)$, but we omit them here.)
\begin{itemize}
\item $\SWbbhalftot^k(X)$ is a cohomology class of the whole moduli space $\BDiff^{+}(X)$, rather than $\textnormal{B}G$ for a proper subgroup $G\subset\Diff^{+}(X)$. So this characteristic class can be defined for all oriented smooth bundles with fiber $X$, with no constraint on the monodromy.
This is a major difference from known gauge-theoretic characteristic classes \cite{K21}.
\item $\SWbbhalftot^k(X)$ is {\it unstable} under stabilizations by $S^2\times S^2$ (\cref{cor: vanishing}).
The non-surjectivity result in \cref{thm: main cal} is a consequence of this property.
\item Just as the Seiberg--Witten invariant detects exotic structures of 4-manifolds,  $\SWbbhalftot^k(X)$ detects subtle differences between smooth families of 4-manifolds.
This yields the non-injectivity result in \cref{thm: main cal} and the comparison results on the smooth and topological categories (\cref{thm: main cal} (a), \cref{thm: Diff Homeo sequence general}).
\item $\SWbbhalftot^k(X)$ can be non-trivial on $\BDiff^+(X)^\delta$. 
Precisely, the pull-back of $\SWbbhalftot^k(X)$ to $H^k(\BDiff^+(X)^\delta;\Z/2)$ under the map induced from the identity  $\Diff^+(X)^\delta \to \Diff^+(X)$ is non-zero for some $X$.
This makes it useful in the study of the algebraic structure of the discrete group $\Diff^{+}(X)^{\delta}$ and yields homological instability for $\BDiff^+(X)^{\delta}$ (\cref{thm: main cal} (d)). 
\item $\SWbbhalftot^k(X)$ vanishes for families that admit fiberwise positive scalar curvature metrics (\cref{thm: vanishing by psc,cor: vanishing for fiberwise psc}). 
This yields a comparison result on moduli spaces of Riemannian metrics and of positive scalar curvature metrics (\cref{thm: psc main intro,thm: psc main}).
\end{itemize}


From a technical point of view, the upshot of the characteristic classes $\SWbbtot^k(X)$ and $\SWbbhalftot^k(X)$ is that they are free from a choice of $\spinc$ structure, despite the dependence of the Seiberg--Witten equations on a $\spinc$ structure. 
The construction is modeled on a characteristic class defined by the first author~\cite{K21} based on the Seiberg--Witten equations, but the characteristic class in \cite{K21} does depend on the choice of a $\spinc$ structure.
To get rid of such constraint, we use an idea due to Ruberman~\cite{Rub01} to define a numerical Seiberg--Witten-type invariant of a diffeomorphism that does not necessarily preserve a given $\spinc$ structure.
See \cref{rem: relation with Rub and K} for a detailed comparison of this work with \cite{K21} and \cite{Rub01}.

Here is a rough idea of the definition of $\SWbbtot^k(X)$.
Let $d(\fraks)$ denote the formal dimension of the Seiberg--Witten moduli space for $\fraks$, and let $\Spinc(X,k)$ denote the set of $\spinc$ structures on $X$ with $d(\fraks)=-k$.
The diffeomorphism group $\Diff^+(X)$ naturally acts on $\Spinc(X,k)$.
Once we fix a pair $\sigma$ of a fiberwise metric and perturbation for the universal bundle $\EDiff^+(X) \to \BDiff^+(X)$, we have a family over $\BDiff^+(X)$ of ``collection (indexed by elements of $\Spinc(X,k)$) of the moduli spaces" of solutions to the Seiberg--Witten equations.
We define a $k$-cochain $\SWcaltot^k(X,\sigma)$ on $\BDiff^+(X)$ by counting this family of collections of moduli spaces over each $k$-cell of $\BDiff^+(X)$.
This cochain shall be shown to be a cocycle and the cohomology class $\SWbbtot^k(X) = [\SWcaltot^k(X,\sigma)]$ will turn out to be independent of the choice of $\sigma$.

The other characteristic class $\SWbbhalftot^k(X)$ is a variant of $\SWbbtot^k(X)$,  corresponding to the quotient of all Seiberg--Witten moduli spaces divided by the charge conjugation on $\spinc$ structures.
For our purpose to prove the homological instability, $\SWbbhalftot^k(X)$ shall be effectively used because the conjugation symmetry annihilates $\SWbbtot^k(X)$ in a lot of situations. See \cref{subsection Motivation} for more detailed motivation for introducing $\SWbbhalftot^k(X)$.


\subsection{Outline}
The remaining sections of this paper are as follows.
In \cref{sec: The Seiberg--Witten characteristic classes}, we construct characteristic classes $\SWbbtot^\bullet(X), \SWbbhalftot^\bullet(X)$, and study some of the basic properties.
The definition is given in \cref{subsec cochain SWtot}, and subsections until there are devoted to giving technical preliminaries.
We prove the vanishing of these characteristic classes under stabilizations and positive scalar curvatures in \cref{subsectionVanishing under stabilizations,subsection Vanishing under positive scalar curvature}, respectively.
In \cref{sec Calculation} we calculate the characteristic class $\SWbbhalftot$.
A key computational result is \cref{thm key computation source}.
In \cref{section BDiffL}, we define the moduli space $\BHomeo^{+}_{L}(X)$ using simplical sets. The main results of this paper, \cref{thm: main cal,thm: Diff Homeo sequence general,thm: abelianisation noninjective,thm: psc main intro}, are proven in \cref{subsectionProof of the main instability theorem}. 
Main ingredients of the proofs of them are the key computation (\cref{thm key computation source}) and a result on geography involving the usual (i.e. unparameterized) Seiberg--Witten invariant, which is given as \cref{thm: 4-mfds that dissolves}.
\cref{construction of 4-manifolds} is devoted to proving \cref{thm: 4-mfds that dissolves}. 

\subsection{Acknowledgement} The authors would like to thank Dave Auckly, David Baraglia, Yi Gu, Sander Kupers, Ciprian Manolescu, Mark Powell, Oscar Randal-Williams, Danny Ruberman, Shun Wakatsuki, and Tadayuki Watanabe for many enlightening discussions and comments on an earlier version of the paper.
H.~Konno was partially supported by JSPS KAKENHI Grant Numbers 19K23412 and 21K13785, and Overseas Research Fellowships.
This work was completed while the first author was in residence at the Simons Laufer Mathematical Sciences Institute in Berkeley, California, during the Fall 2022 semester.

\section{The Seiberg--Witten characteristic classes}
\label{sec: The Seiberg--Witten characteristic classes}

\subsection{Motivation}
\label{subsection Motivation}
In this \lcnamecref{sec: The Seiberg--Witten characteristic classes}, given $k\geq0$ and a 4-manifold $X$ with $b^+(X)>k+1$, we define characteristic classes $\SWbbtot^k(X), \SWbbhalftot^k(X)$
mentioned in \cref{subsection Intro SW char}.
As the construction of $\SWbbhalftot^k(X)$ is rather complicated than that of $\SWbbtot^k(X)$,
here we clarify why we need the characteristic class $\SWbbhalftot^k(X)$ to motivate readers.
First we describe why we adopt (twisted) $\Z$-coefficient for  $\SWbbtot$. 
If one tries to define a characteristic class using the Seiberg--Witten equations in full generality, it would be natural to define a refinement $\widetilde{\mathbb{SW}}_{\mathrm{tot}}^k(X) \in H^k(\BDiff^+(X;\tilde{\Z}_E^\infty)$ of $\SWbbtot^k(X)$,
where $\tilde{\Z}_E^\infty$ is a local system with fiber $\bigoplus_{\fraks \in \Spinc(X,k)}\Z$ determined by the monodromy action for $\EDiff^+(X)$ on $\Spinc(X,k)$ and on the homology orientation.
The class $\SWbbtot$ can be recovered from $\widetilde{\mathbb{SW}}_{\mathrm{tot}}$ via a homomorphism 
$\bigoplus_{\fraks \in \Spinc(X,k)}\Z  \to \Z$ defined by summing up entries.
However, the action of $\Diff^+(X)$
on $\Spinc(X,k)$ may be quite complicated, and we do not have control of the local system $\tilde{\Z}_E^\infty$.
Thus we eventually pass to $\SWbbtot$ for a computational reason.

Moreover, even after passing to $\SWbbtot$, one still has to work with the local system $\Z_E=\Z_{\EDiff^+(X)}$ because the action of $\Diff^+(X)$ on the homology orientation is nontrivial in general. Since our main interest is the homology of moduli spaces for untwisted coefficient, we may hope to study $\SWbbtot$ after the mod 2 reduction.

However, $\SWbbtot$ corresponds to summing up the counts of the moduli spaces over $\Spinc(X,k)$, so almost every time the charge conjugation symmetry of the Seiberg--Witten equations kills the $\Z/2$-coefficient $\SWbbtot$.
Now we arrive at the reason why we need to take the ``half" of $\Spinc(X,k)$, namely $\Spinc(X,k)/\Conj$, to get interesting results for a broad class of 4-manifolds.

Nevertheless, in this paper we shall also describe the construction of  $\SWbbtot$ not only $\SWbbhalftot$, because the construction of $\SWbbtot$ gives a good guide in the construction of the more complicated object $\SWbbhalftot$.

\subsection{Virtual neighborhoods for families}
\label{subsection: Virtual neighborhoods for families}

To define the characteristic classes $\SWbbtot^\bullet(X)$ and $\SWbbhalftot^\bullet(X)$ for a 4-manifold $X$,
we need to treat the issue of transversality in counting parameterized moduli spaces over $\BDiff^+(X)$, which is not a finite-dimensional smooth manifold. This may cause a complication, in particular, in proving that cochains we construct to define $\SWbbtot^\bullet(X)$ and  $\SWbbhalftot^\bullet(X)$ are cocycles (corresponding to \cref{SWtot cocycle}). Since we are studying the (co)homology group of $\BDiff^+(X)$, we may assume that $\BDiff^+(X)$ is a CW complex by passing to its CW approximation. 
Then we shall adopt the virtual neighborhood technique along Ruan~\cite{Ru98} to bypass this transversality issue.
This technique is easily generalized to a parameterized setup \cite[Subsection~5.2]{K21}.
Further, in \cref{subsectionFamily of collections of Hilbert bundles}, we need to introduce a slight generalization of virtual neighborhoods for families.
For the reader's convenience,
in this \lcnamecref{subsection: Virtual neighborhoods for families}, we summarize necessary facts on and constructions of virtual neighborhoods for families along \cite[Subsection~5.2]{K21}.

\begin{rmk}
Another potential way to handle the transversality issue could be taking a particular model of $\BDiff^+(X)$ that is an infinite-dimensional manifold such as the space of embeddings of $X$ into $\R^\infty$ divided by $\Diff^+(X)$ as in,  e.g., \cite{Galatius2018}, but the authors never checked details along this line.
\end{rmk}

\begin{convention}
In this paper, unless otherwise stated, a {\it section} of a continuous fiber bundle over a topological space means a continuous section.
\end{convention}

We first fix our notation.
Recall that the Seiberg--Witten equations give rise to a Fredholm section of a certain Hilbert bundle over a certain Hilbert manifold.
We shall consider a parameterized version of this setup, and abstract it as follows (see \cite[Subsection~5.2]{K21} for the precise setup).
The objects we are going to introduce are summarized in the diagram
\begin{equation*}
\vcenter{
\xymatrix{
    \scrE \ar[dd]\ar[rd]&{}\\
    {}& \scrX \ar[ld]^{\pi} \ar@/_1pc/[lu]_{s}\\
    B. &
}
}
\end{equation*}
First, $B$ is a normal space, which is used as the base space or parameter space in this \lcnamecref{subsection: Virtual neighborhoods for families}.
Next, $\scrX = \bigcup_{b \in B}\scrX_b$ is a continuous (locally trivial) family of Hilbert manifolds over $B$ with the projection $\pi : \scrX \to B$, and $\scrE$ is a parameterized Hilbert bundle over $B$.
Namely, $\scrE = \bigcup_{b \in B}\scrE_b$ is a continuous (locally trivial) family of vector bundles $\scrE_b \to \scrX_b$ with Hilbert space fiber $\scrH_b$.
Finally, $s : \scrX \to \scrE$ is a parameterized Fredholm section.
Namely, $s$ is a continuous map commuting with the projections onto $B$ such that, at each $b \in B$, the restriction $s_b : \scrX_b \to \scrE_b$ of $s$ is a smooth Fredholm section.

Set $\calM_b = s_b^{-1}(0)$ for each $b \in B$, called the {\it unparameterized moduli space} at $b$.
Let $\calM$ denote the union of the unparameterized moduli space:
\[
\calM
= \bigcup_{b \in B}\calM_b.
\]
We call $\calM$ the {\it parameterized moduli space} over $B$.

\begin{assum}
Throughout this \lcnamecref{subsection: Virtual neighborhoods for families}, suppose that $\calM$ is compact.
\end{assum}

Suppose also that the index of $d(s_{\pi(x)})_x$ 
is independent of $x \in \calM$. We denote by $\ind{s}$ the index.
Let $\det{ds} \to \scrX$ denote the parameterized determinant line bundle of the parameterized linearized Fredholm sections $ds$.
Let $\calL$ be the local system on $\scrX$ with fiber $\Z$ given as the orientation local system of $\det{ds}$.
Pick a point $b \in B$, and suppose that the monodromy action of $\pi_1(\scrX_b)$ on a fiber of $\det{ds_b}$  is trivial.
In this case, $\calL$ is isomorphic to the pull-back of a local system over $B$, denoted also by $\mathcal{L}$, along $\scrX \to B$.

An output from virtual neighborhood techniques we need can be summarized as follows.
Let $H^\ast_{\mathrm{cpt}}(-)$ denote the compactly supported cohomology group.

\begin{pro}
\label{pro: summary of vn}
Under the above assumptions,
we obtain a cohomology class 
\[
\frakM(s) \in H^{- \ind{s}}_\cpt(B; \calL)
\]
with the following properties:
\begin{itemize}
\item If $\calM=\emptyset$, then $\mathfrak{M}(s)=0$.
\item Suppose that $\ind{s} \leq 0$.
If $B$ is a smooth oriented closed manifold of dimension $-\ind{s}$ and $\calL$ is a trivial local system, and further if $s$ is generic so that $\calM$ is a (0-dimensional) smooth manifold, then we have
\[
\#\calM = \left< \frakM(s), [B] \right>.
\]
\end{itemize}

Moreover, if the section $s$ is nowhere vanishing over a subset $B' \subset B$, we may obtain a relative cohomology class
\[
\frakM(s; B') \in H^{- \ind{s}}_\cpt(B, B'; \calL)
\]
with properties analogous to the above.
\end{pro}

We shall use the following naturality result,
which was given in \cite[Corollary~5.19]{K21} (there is an obvious typographical error in \cite[Corollary~5.19]{K21}, corrected below):

\begin{lem}
\label{lem: first naturality}
Let $\scrX \to B, \scrE \to \scrX \to B, s : \scrX \to \scrE$ be as above, described at the beginning of this \lcnamecref{subsection: Virtual neighborhoods for families}.
Let $A$ be a normal space and $f : A \to B$ be a continuous map.
Suppose that, for the pull-back section $f^\ast s : f^\ast \scrX \to f^\ast \scrE$,
the parameterized zero set $(f^\ast s)^{-1}(0)$ is also compact as well as $s^{-1}(0)$.
Then we have
\[
f^\ast \frakM(s)
= \frakM(f^\ast s)
\]
in $H^{-\ind{s}}_\cpt(A;f^\ast\calL)$. 
Similarly, for subsets $A' \subset A$, $B'\subset B$ with $f(A') \subset B'$, if $s$ is nowhere vanishing over $B'$, we have
\[
f^\ast \frakM(s; B')
= \frakM(f^\ast s; A')
\]
in $H^{-\ind{s}}_\cpt(A, A';f^\ast\calL)$.
\end{lem}

In the rest of this \lcnamecref{subsection: Virtual neighborhoods for families}, we sketch the construction of the cohomology classes $\frakM(s)$ and $\frakM(s;B')$.
First, under the above assumptions, we shall construct a (typically non-locally trivial) family of smooth manifolds
\[
\calU = \bigcup_{b \in B}\calU_b
\]
over $B$ such that $\calM$ is embedded in the interior of $\calU$ for a natural topology on $\calU$.
This space $\calU$ is called a {\it families virtual neighborhood} for $\calM$ (or for $s$) \cite[Definition~5.10]{K21}.

The first step to construct $\mathcal{U}$ is to find a natural number $N$ and a fiberwise smooth map
\[
\varphi : \scrX \times \R^N \to \scrE
\]
that satisfies that $\scrX \times \{0\} \subset \varphi^{-1}(0)$ and the following property: set $\tilde{s} = s+\varphi$.
Then, for each $x \in \calM$, the differential
\[
d\left(\tilde{s}_{\pi(x)}\right)_{(x,0)} : T_x\scrX_{\pi(x)} \oplus \R^N \to \scrH_{\pi(x)}
\]
is surjective.

Such $N$ and $\varphi$ can be found by using the compactness of $\calM$: at each point $x \in \calM$, one can easily find a natural number $N_x$ and a linear map $f_x : \R^{N_x} \to \scrH_{\pi(x)}$ such that $ds_{\pi(x)}+f_x : T\scrX_{\pi(x)} \oplus \R^{N_x}\to \scrH_{\pi(x)}$ is surjective.
Since the surjectivity is an open condition, we may take an open neighborhood $U_x$ of $x$ in $\scrX$ so that, for all $y \in U_x \cap \calM$, $ds_{\pi(y)}+f_x$ are surjective.
By the compactness of $\calM$, we may cover $\calM$ with finitely many $U_x$.
Let $N$ be the sum of $N_x$ for such $x$, and then we may obtain $\varphi$ above by summing the pull-back of $f_x$ to $\scrX \times \R^N$ with multiplying cut-off functions on $U_x$.

Again since the surjectivity is an open condition, we may find an open neighborhood $\scrN$ of $\calM \times \{0\}$ in $\scrX \times \R^N$ such that the differential of $\tilde{s}_{(\pi(x),0)}$ along the fiber-direction is surjective at any point $x \in \tilde{s}^{-1}(0) \cap \scrN$.

\begin{defi}
Set
\[
\calU = \calU(s, N, \varphi, \scrN) := \tilde{s}^{-1}(0) \cap \scrN.
\]
We call $\calU$ a {\it families virtual neighborhood} for $\mathcal{M}$.
\end{defi}
Denote by $\calU_b$ the fiber of $\calU$ over $b \in B$ under the obvious map $\calU \to B$.
By the implicit function theorem, $\calU_b$ is a smooth manifold of dimension $\ind(ds_b) + N$ for each $b$.
Clearly, $\calU$ is a neighborhood of $\calM \times \{0\}$ in $\scrX \times \R^N$, and $\calU_b$ is a neighborhood of $\calM_b$ in $\scrX_b \times \R^N$ for each $b \in B$.

\begin{rmk}
Note that $\calU$ is not uniquely determined by $s$: there are various choices of $N, \varphi$, and $\scrN$.
Note also that $\calU$ is not necessarily a locally trivial family over $B$.
For example, $\calU$ is supported only over $\pi(\calM) \subset B$: if we have $\calM_b=\emptyset$ for $b \in B$, then (a small neighborhood of $\calM_b$ in) $\calU_b$ is also empty.
\end{rmk}


Next, we consider a certain relative Euler class.
Let $h_\calU : \calU \to \R^N$ denote the `height' function, defined as the restriction of the second projection $\scrX \times \R^N \to \R^N$ to $\calU \subset \scrX \times \R^N$.
The parameterized moduli space $\calM$ is the level set for $h_{\calU}$ of height zero: $\calM \times \{0\} = h_\calU^{-1}(0)$.
The map $h_\calU: \calU \to \R^N$ defines a section of the trivial bundle $\calU \times \R^N \to \calU$ by $u \mapsto (u,h_\calU(u))$, which we denote by $h_\calU' : \calU \to \calU \times \R^N$.
Consider the relative Euler class for $h_\calU'$, namely
\[
e_\calU := (h_\calU')^\ast \tau(\calU \times \R^N \to \calU) \in H^N_\cpt(\calU; \Z),
\]
where $\tau(\calU \times \R^N \to \calU)$ denotes the Thom class of the trivial bundle $\calU \times \R^N \to \calU$.
Note that $e_\calU$ is compactly supported since $\calM$ was supposed to be compact.

As the last input, we have  the integration along fibers for the family $\calU \to B$, while $\calU \to B$ is not necessarily locally trivial (see \cite[explanations below Equation (21)]{K21}):
\[
\pi_! : H^\ast_\cpt(\calU; \Z) \to H^{\ast - (\ind{s}+N)}_\cpt(B; \calL).
\]
Note that $\ind{s}+N$ coincides with the dimension of $\calU_b$ for any $b \in B$.

\begin{defi}
\label{defi: vn general coh}
We define a cohomology class 
$\frakM(s) \in H^{- \ind{s}}_\cpt(B; \calL)$
by $\frakM(s) = \pi_!(e_\calU)$.
\end{defi}


Now, suppose that the section $s$ is nowhere vanishing over a subset $B' \subset B$.
Then $\frakM(s)$ is supported outside $B'$, and thus the following definition makes sense:

\begin{defi}
\label{defi: vn general coh rel}
We define a cohomology class
$\frakM(s; B') \in H^{- \ind{s}}_\cpt(B, B'; \calL)$
by $\frakM(s; B') = \pi_{!}(e_{\mathcal{U}})$.
\end{defi}

As proven in \cite[Lemma~5.15]{K21}, $\frakM(s)$ and $\frakM(s;B')$ are independent of the choice of $\calU$.

\subsection{Family of collections of Hilbert bundles}
\label{subsectionFamily of collections of Hilbert bundles}

For the purpose of this paper, it is important to note that all constructions in \cref{subsection: Virtual neighborhoods for families} can be generalized to a Fredholm section of a `family of collections' of Hilbert bundles, defined below.
To avoid confusion, throughout we shall use the word `family' for the case that the parameter space is not discrete, and the word `collection' for the case that the parameter space is discrete.

Let $\Lambda$ be a set equipped with the discrete topology. 
(In our application to Seiberg--Witten theory, we shall finally take $\Lambda$ to be the set of $\spinc$ structures on a 4-manifold with a fixed formal dimension.)
Suppose that we are given a collection of (unparameterized) Hilbert bundles over Hilbert manifolds:
\[
(\scrE_\lambda^0 \to \scrX_\lambda^0)_{\lambda \in \Lambda}.
\]
(The superscript ``$0$" indicates that the object is unparameterized.)
An {\it automorphism} of the collection $(\scrE_\lambda^0 \to \scrX_\lambda^0)_{\lambda}$ consists of a bijection $\alpha : \Lambda \to \Lambda$ and a collection of isomorphisms of Hilbert bundles 
\[
\xymatrix{
    \scrE_\lambda^0 \ar[r]\ar[d] & \scrE_{\alpha(\lambda)}^0\ar[d]\\
    \scrX_\lambda^0 \ar[r] & \scrX_{\alpha(\lambda)}^0.
}
\]
Let $\Aut\left((\scrE_\lambda^0)_\lambda\right)$ denote the set of all automorphisms of $(\scrE_\lambda^0 \to \scrX_\lambda^0)_{\lambda}$, which forms a group under composition.
We call an $\Aut\left((\scrE_\lambda^0)_\lambda\right)$-bundle a {\it family of collections of Hilbert bundles} modeled on $(\scrE_\lambda^0 \to \scrX_\lambda^0)_\lambda$.

Let $\scrE \to \scrX \to B$ be a family over a normal space $B$ of collections of Hilbert bundles modeled on $(\scrE_\lambda^0)_\lambda$.
For each $b \in B$, let $\scrX_b$ denote the fiber of $\scrX \to B$ over $b$, and similarly define $\scrE_b$.
Then, on each $b \in B$, we have a pair of uncanonical identifications
\begin{align}
\label{eq: indentifications}
\scrX_b \cong \bigsqcup_{\lambda \in \Lambda}\scrX_\lambda^0,\quad
\scrE_b \cong \bigsqcup_{\lambda \in \Lambda}\scrE_\lambda^0,
\end{align}
and the set of choices of identification is given by $\Aut\left((\scrE_\lambda^0)_\lambda\right)$.
Note that we have a natural map $\scrE \to \scrX$ induced from the collection of the projections $\scrE_\lambda^0 \to \scrX_\lambda^0$.
Thus we have the notions of a (parameterized) section $s : \scrX \to \scrE$, which is a family of sections $s_b : \scrX_b \to \scrE_b$.
Similarly, we have the notion that $s$ is (fiberwise) Fredholm.
This is equivalent to say that, for each $b \in B$,
once we fix a pair of identifications \eqref{eq: indentifications}, 
the induced sections $(s_\lambda^0)_b : \scrX_\lambda^0 \to \scrE_\lambda^0$ are Fredholm for all $\lambda$.

Let $s :\scrX \to \scrE$ be a parameterized Fredholm section.
For each $b \in B$, set $\calM_{b} = s_b^{-1}(0)$.
Again, if we fix a pair of identifications \eqref{eq: indentifications} and use the induced sections, there is a homeomorphism
\[
\calM_{b} 
\cong \bigsqcup_{\lambda \in \Lambda}(s_\lambda^0)_b^{-1}(0).
\]
Define
\[
\calM := 
\bigcup_{b \in B}\calM_{b} \subset \scrX.
\]

Suppose that $\calM$ is compact, and the Fredholm indices of the linearizations of $s_\lambda$ to the fiber-direction are independent of $\lambda$.
Then we may repeat the constructions of a virtual neighborhood $\calU$ for $s$ and the cohomology class $\frakM(s)$ in \cref{subsection: Virtual neighborhoods for families} word for word, and claims analogous to \cref{pro: summary of vn,lem: first naturality} hold without any essential changes.

\subsection{Collection of Seiberg--Witten equations}
\label{subsection Collection of Seiberg--Witten equations}

We shall apply the abstract machinery explained until the last \lcnamecref{subsectionFamily of collections of Hilbert bundles} to the Seiberg--Witten equations.
First we fix our notation for the Seiberg--Witten equations.
Let $X$ be a closed oriented smooth 4-manifold.
As in \cite[Section~2]{K21}, we define a $\spinc$ structure on $X$ using the double covering of $\mathrm{SL}(4,\R)$ instead of $\mathrm{Spin}(4)$ to avoid using a Riemannian metric on $X$.
Once  we take a metric on $X$, then a $\spinc$ structure in our sense gives rise to a $\spinc$ structure in the usual sense, and the isomorphism classes of $\spinc$ structures in our sense bijectively correspond to those of $\spinc$ structures in the usual sense.
Let $\Spinc(X)$ denote the set of isomorphism classes of $\spinc$ structures on $X$.

Given a Riemannian metric $g$ and a $\spinc$ structure $\texts$ on $X$, let $L$ denote the determinant line bundle of $\texts$ and let $S^+(X,\texts,g)$ and $S^-(X,\texts,g)$ denote the positive and negative spinor bundles respectively.
Let $\Lambda^+_g(X)$ denote the self-dual part of $\Lambda^2(X)$ for the metric $g$.
Fix $l > 1$, and let $\scrA(X,\texts,g)_{L^2_l}$ and $L^2_l(S^\pm(X,\texts,g))$ denote the space of $\mathrm{U}(1)$-connections of $L$ and sections of $S^\pm(X,\texts,g)$ completed by the $L^2_l$-norms.
Set
\begin{align*}
 &\scrC(X,\texts,g) := \scrA(X,\texts,g)_{L^2_l} \times L^2_l(S^+(X,\texts,g)),\\
 &\scrC^\ast(X,\texts,g) := \scrA(X,\texts,g)_{L^2_l} \times (L^2_l(S^+(X,\texts,g)) \setminus \{0\}),\\
 &\scrD(X,\texts,g) := L^2_{l-1}(i\Lambda^+_g(X)) \times L^2_{l-1}(S^-(X,\texts,g)).
\end{align*}
Let $\scrR(X)$ denote the space of Riemannian metrics on $X$ and define
\[
\Pi_\ast(X) = \bigcup_{g \in \scrR(X)} L^2_{l-1}(i\Lambda^+_g(X)),
\]
which is a Hilbert bundle over $\scrR(X)$.
(We shall consider a subspace of $\Pi_\ast(X)$ with a boundedness condition, which will be denoted by $\Pi(X)$ dropping $\ast$.)
Let $\pi : \Pi_\ast(X) \to \scrR(X)$ denote the natural projection.
Set
\[
\circPi_\ast(X) = \bigcup_{g \in \scrR(X)} L^2_{l-1}(i\Lambda^+_g(X)) \setminus \mathrm{Im}(d^+_g : L^2_{l}(i\Lambda^1(X)) \to L^2_{l-1}(i\Lambda^+_g(X))).
\]
This space $\circPi_\ast(X)$ is a subspace of $\Pi_\ast(X)$ of codimension-$b^+(X)$.
Since $\Pi_\ast(X)$ is contractible, $\circPi_\ast(X)$ is $(b^+(X)-2)$-connected.

For each $\mu \in \Pi_\ast(X)$,
let 
\[
\tilde{s}_\mu : \scrC(X,\texts,\pi(\mu)) \to \scrD(X,\texts,\pi(\mu))
\]
be the map corresponding to the $\mu$-perturbed Seiberg--Witten equations:
\[
\tilde{s}_\mu(A,\Phi)
= (F^{+_{\pi(\mu)}}_A - \mu - \sigma(\Phi, \Phi), D_A\Phi),
\]
where $\sigma(\Phi, \Phi)$ denotes the quadratic term for spinors in the Seiberg--Witten equations and $D_A$ is the $\spinc$ Dirac operator.
Set $\scrG = L^2_{l+1}(X, \mathrm{U}(1))$ and 
\begin{align*}
&\scrB(X,\texts,g) = \scrC(X,\texts,g)/\scrG,\\
&\scrB^\ast(X,\texts,g) = \scrC^\ast(X,\texts,g)/\scrG,\\
&\scrE(X,\texts,g) = (\scrC^\ast(X,\texts,g) \times \scrD(X,\texts,g))/\scrG.
\end{align*}
Then $\scrE(X,\texts,g) \to \scrB^\ast(X,\texts,g)$ is a Hilbert bundle with fiber $\scrD(X,\texts,g)$.
Since $\tilde{s}_\mu$ is $\scrG$-equivariant, it gives rise to a section
\[
s_\mu : \scrB^\ast(X,\texts,\pi(\mu)) \to \scrE(X,\texts,\pi(\mu)),
\]
which is a Fredholm section of index
\[
d(\texts) = \frac{1}{4}(c_1(\texts)^2 -2\chi(X) -3\sign(X)),
\]
which is the formal dimension of the Seiberg--Witten moduli space for $\texts$.

Now we discuss isomorphisms of $\spinc$ structures.
Let $\Aut(X,\texts)$ denote the automorphism group of a $\spinc$ 4-manifold $(X,\texts)$, and let $\fraks$ denote the isomorphism class of $\texts$.
Then we have an exact sequence
\[
1 \to \Map(X,\mathrm{U}(1)) \to \Aut(X,\texts) \to \Diff^+(X).
\]
The image of the last map is given as $\Diff(X, \fraks)$, the group of orientation-preserving diffeomorphisms that preserve the isomorphism class $\fraks$.
Let $g$ be a metric on $X$ and let $\texts'$ be a $\spinc$ structure on $X$ that is isomorphic to $\texts$.
Once we choose an isomorphism from $\texts$ to $\texts'$, we have an isomorphism
\begin{align}
\label{eq: iso scrC}
\scrC(X,\texts,g) \to \scrC(X,\texts',g)
\end{align}
of Hilbert manifolds, and it gives rise to an isomorphism
\begin{align}
\label{eq: iso scrB}
\scrB^\ast(X,\texts,g) \to \scrB^\ast(X,\texts',g)
\end{align}
of Hilbert manifolds.
While the isomorphism \eqref{eq: iso scrC} depends on the choice of isomorphism from $\texts$ to $\texts'$, the isomorphism \eqref{eq: iso scrB} between quotients by $\scrG$ is canonically determined by $\texts$ and $\texts'$, since the choice of isomorphism from $\texts$ to $\texts'$ is in $\scrG$.
The same remark also applies to $\scrE(X,\texts,g)$.
For this reason, we write $\scrB^\ast(X,\fraks,g)$ and $\scrE(X,\fraks,g)$ for $\scrB^\ast(X,\texts,g)$ and $\scrE(X,\texts,g)$ respectively.

Now we consider the collection of the above objects indexed by $\fraks$.
Set
\begin{align*}
&\scrB^\ast(X) = \bigsqcup_{\fraks \in \Spinc(X)}
\bigcup_{g \in \scrR(X)}
\scrB^\ast(X,\fraks,g),\\
&\scrE(X) = \bigsqcup_{\fraks \in \Spinc(X)}
\bigcup_{g \in \scrR(X)}
\scrE(X,\fraks,g).
\end{align*}
Here $\bigcup_{g \in \scrR(X)}
\scrB^\ast(X,\fraks,g)$ is regarded as a locally trivial family of Hilbert manifolds over $\scrR(X)$, and $\scrB^\ast(X)$ is a disjoint union of such families of Hilbert manifolds.
We note that $\Diff^+(X)$ acts on $\scrB^\ast(X)$ and on $\scrE(X)$ in a natural way.
To see this, let $f \in \Diff^+(X)$, $\fraks \in \Spinc(X)$ and $g \in \scrR(X)$.
If we take a representative $\texts$ of $\fraks$ and an isomorphism $\tilde{f} : f^\ast \texts \to \texts$ that covers $f : X \to X$,
we obtain an isomorphism
\begin{align}
\label{eq: pullback}
f^\ast : \scrB^\ast(X,\fraks,g) \to \scrB^\ast(X,f^\ast\fraks,f^\ast g)
\end{align}
of Hilbert manifolds.
However, for the same reason that the isomorphism \eqref{eq: iso scrB} is canonically determined by $\texts$ and $\texts'$, the isomorphism \eqref{eq: pullback} is independent of the choice of $\texts$ and $\tilde{f}$.
Thus we obtain a canonical isomorphism \eqref{eq: pullback} determined only by $f$, and hence $f$ yields a canonical isomorphism of $\scrB^\ast(X)$, a disjoint union of families of Hilbert manifolds.
Thus $\Diff^+(X)$ acts on $\scrB^\ast(X)$, and similarly also on $\scrE(X)$ so that the projection $\scrE(X) \to \scrB^\ast(X)$ is $\Diff^+(X)$-equivariant.

Also, we note that the above action of $\Diff^+(X)$ is compatible with the Fredholm section corresponding to the Seiberg--Witten equations.
Notice that we have a natural map $\scrE(X) \to \scrB^\ast(X)$, and the sections $s_\mu : \scrB^\ast(X,\texts,\pi(\mu)) \to \scrE(X,\texts,\pi(\mu))$ yield a section
\[
s : \scrB^\ast(X) \to \scrE(X)
\]
of this natural map.
Since $s_\mu$ was induced from a $\scrG$-equivariant map, again for the same reason as the canonicality of \eqref{eq: pullback}, $s$ is $\Diff^+(X)$-equivariant under the above actions on $\scrB^\ast(X), \scrE(X)$.

For the isomorphism class $\fraks$ of a $\spinc$ structure $\texts$ on $X$, we set $d(\fraks) = d(\texts)$.
Given an integer $k$, let $\Spinc(X,k)$ denote the set of isomorphism classes $\fraks$ of $\spinc$ structures on $X$ with $d(\fraks)=-k$.
We define subsets $\scrB^\ast(X, k), \scrE(X, k)$ of $\scrB^\ast(X), \scrE(X)$ to be
\begin{align*}
&\scrB^\ast(X,k) = \bigsqcup_{\fraks \in \Spinc(X,k)}
\bigcup_{g \in \scrR(X)}
\scrB^\ast(X,\fraks,g),\\
&\scrE(X,k) = \bigsqcup_{\fraks \in \Spinc(X,k)}
\bigcup_{g \in \scrR(X)}
\scrE(X,\fraks,g).
\end{align*}
Then $\Diff^+(X)$ acts on $\scrB^\ast(X,k), \scrE(X,k)$ as well, and the section $s^k : \scrB^\ast(X,k) \to \scrE(X,k)$ obtained as the restriction of the above $s$ is $\Diff^+(X)$-equivariant.

From now on, we shall consider a parameterized setup. Let $X \to E \to B$ be a fiber bundle with fiber $X$ with structure group $\Diff^+(X)$.

\begin{nota}
\label{notation data for bundle}
Let $A(X,a)$ be an object determined by $X$ and a datum $a$, and suppose that $A(X,a)$ is acted by $\Diff^+(X)$.
We denote by the associated bundle with $E$ with fiber $A(X,a)$ by 
\[
A(X,a) \to A(E,a) \to B.
\]
\end{nota}

For example, the all objects above, which are acted by $\Diff^+(X)$, give rise to fiber bundles associated with $E$:
\begin{align*}
\scrB^\ast(X,k) \to \scrB^\ast(E,k) \to B,&\\
\scrE(X,k) \to \scrE(E,k) \to B,&\\
\Pi_\ast(X) \to \Pi_\ast(E) \to B,&\\
\circPi_\ast(X) \to \circPi_\ast(E) \to B,&\\
\scrR(X) \to \scrR(E) \to B.&
\end{align*}

\begin{rmk}
The notation $\scrR(E)$ does not mean ``the space of metrics on the total space $E$". 
(Note that even the total space of $E$ does not have a natural manifold structure if $B$ is not a manifold.)
Instead, a section of $\scrR(E) \to B$ corresponds to a fiberwise metric on $E$.
Similar remarks apply to many families over $B$ that shall appear later.
\end{rmk}

We denote the fiber of each bundle, say $\scrB^\ast(E,k)$, over a point $b \in B$ by $\scrB^\ast(E_b,k)$.
Note that, since $s^k : \scrB^\ast(X,k) \to \scrE(X,k)$ is $\Diff^+(X)$-equivariant, this induces a section
\[
s^k : \scrB^\ast(E,k) \to \scrE(E,k)
\]
that commutes with the projections $\scrB^\ast(E,k) \to B$, $\scrE(E,k) \to B$.

To get some compactness of moduli spaces later, we consider bounded perturbations.
Note that the fiber of $\Pi_\ast(E) \to B$ over a point $b \in B$ is of the form
\[
(\Pi_\ast(E))_b 
= \bigcup_{g \in \scrR(E_b)} L^2_{l-1}(i\Lambda^+_g(E_b)).
\]
Let $\mathring{D}(i\Lambda^+_g(E_b))$ denote the open ball centered at the origin of radius 1 with respect to the $L^2_{l-1}$-norm.
Define a subset $\Pi(E)$ of $\Pi_\ast(E)$ by
\[
\Pi(E) = \bigcup_{b \in B}
\bigcup_{g \in \scrR(E_b)} \mathring{D}(i\Lambda^+_g(E_b)).
\]
We also set 
\[
\circPi(E) = \circPi_\ast(E) \cap \Pi(E).
\]


Now we choose a fiberwise metric and self-dual 2-forms on $E$.
Recall that we denoted by $\pi$ the projection $\Pi_\ast(X) \to \scrR(X)$.
This gives rise to a map (denoted by the same symbol) $\pi : \Pi_\ast(E) \to \scrR(E)$.
Restricting this, we have a map $\pi : \circPi(E) \to \scrR(E)$.
Taking a section $\sigma : B \to \circPi(E)$ corresponds to taking a fiberwise metric and a fiberwise self-dual 2-form (with fiberwise boundedness) on $E$ that continuously vary over $B$.
Let us consider the projections $\scrB^\ast(X,k) \to \scrR(X)$ and $\scrE(X,k) \to \scrR(X)$.
For example, the fiber of the projection $\scrB^\ast(X,k) \to \scrR(X)$ over a metric $g \in \scrR(X)$ is the collections (over $\Spinc(X,k)$) of all configurations for the common metric $g$.
These projections give rise to natural maps
\[
p_{\scrB} : \scrB^\ast(E,k) \to \scrR(E),\quad
p_{\scrE} : \scrE(E,k) \to \scrR(E).
\]
We set
\[
\scrB^\ast(E,k,\sigma)
:= p_{\scrB}^{-1}(\pi\circ\sigma(B)),\quad
\scrE^\ast(E,k,\sigma)
:= p_{\scrE}^{-1}(\pi\circ\sigma(B)).
\]

The above section $s^k : \scrB^\ast(E,k) \to \scrE(E,k)$ corresponds to a union of families over $B$ of Seiberg--Witten equations over all fiberwise metrics and fiberwise self-dual 2-forms.
Restricting this to the image of 
$\sigma$, we obtain a section
\[
s^k_{\sigma} : \scrB^\ast(E,k,\sigma) \to \scrE(E,k,\sigma).
\]
The resulting section $s^k_{\sigma}$ is a parameterizd Fredholm section of a family of collections of Hilbert bundles, in the sense of \cref{subsectionFamily of collections of Hilbert bundles}.

\begin{defi}
Set
\begin{align}
\label{eq: SWtot moduli}
\calM_{\sigma} := \bigcup_{b \in B}(s^k_{\sigma})_{b}^{-1}(0).
\end{align}
We call $\calM_{\sigma}$ the {\it family of total collections of moduli spaces} associated with $\sigma$.
\end{defi}

\begin{lem}
\label{lem: finiteness1}
If $B$ is compact, $\calM_{\sigma}$ is also compact.
\end{lem}

\begin{proof}
Recalling that we use a common fiberwise metric $\pi \circ \sigma$ for all $\spinc$ structures to define $\calM_\sigma$, and $\circPi(E)$ consists of bounded perturbations, the claim follows from the standard compactness of the Seiberg--Witten moduli space: the compactness of the moduli space for a fixed $\spinc$ structure and the finiteness of the numbers of $\spinc$ structures for which the Seiberg--Witten moduli spaces are non-empty.
\end{proof}

Thanks to \cref{lem: finiteness1}, we may apply the virtual neighborhood technique in \cref{subsection: Virtual neighborhoods for families,subsectionFamily of collections of Hilbert bundles} to $s_{\sigma}$.
In particular, as described in \cref{defi: vn general coh}, we obtain a cohomology class
\begin{align}
\label{eq: coh compact base tot}
\frakM(s^k_{\sigma}) \in H^{k}(B; \Z_E)
\end{align}
for compact $B$.
Here $\Z_E$ is a local system over $B$ with fiber $\Z$ that is determined by the monodromy action for $E$ on the homology orientation, the orientation of $H^1(X;\R) \oplus H^+(X)$, where $H^+(X)$ denotes a maximal-dimensional positive-definite subspace of $H^2(X;\R)$ with respect to the intersection form.
Similarly, as in \cref{defi: vn general coh rel}, if the section $s^k_{\sigma}$ is nowhere vanishing over a subset $B' \subset B$, we have a relative cohomology class
\begin{align}
\label{eq: coh compact base tot rel}
\frakM(s^k_{\sigma}; B') \in H^{k}(B, B'; \Z_E)
\end{align}
for compact $B$.


\subsection{Charge conjugation}
\label{subsection Complex conjugation}

In this \lcnamecref{subsection Complex conjugation}, we shall consider the charge conjugation on the Seiberg--Witten equations in our context, which is a preliminary to the half-total characteristic class $\SWbbhalftot^\bullet(X)$.
Henceforth, when there is no risk of confusion, we call an isomorphism class $\fraks \in \Spinc(X)$ just a $\spinc$ structure.

Let $X$ be a closed oriented smooth 4-manifold.
Recall that, for each $\fraks \in \Spinc(X)$, there is a unique $\spinc$ structure $\bar{\fraks}$ on $X$ called the charge conjugation of $\fraks$, which satisfies $c_1(\bar{\fraks})=-c_1(\fraks)$.
The group $\Z/2$ acts on $\Spinc(X)$ by the charge conjugation.
Let $\Spinc(X)/\Conj$ denote the quotient.

Let $\fraks \in \Spinc(X)$, and
take a metric $g$ on $X$ and a representative $\texts$ of the isomorphism class $\fraks$.
Let $\bar{\texts}$ denote the charge conjugation of $\texts$.
The charge conjugation induces isomorphisms between Hilbert (affine) spaces
\[
c : \scrA(X,\texts,g)_{L^2_l} \to \scrA(X,\bar{\texts},g)_{L^2_l},\quad
c' : L^2_l(S^\pm(X,\texts,g)) \to L^2_l(S^\pm(X,\bar{\texts},g)).
\]
Let $(-1) : \Omega^+_g(X) \to \Omega^+_g(X)$ be the $(-1)$-muptiplication.
Then the Seiberg--Witten equations are compatible with these maps in the sense that the following diagram commutes for each $\mu \in L^2_{l-1}(i\Lambda^+_g(X))$:
\[
\xymatrix{
    \scrC(X,\texts,g) \ar[r]^{\tilde{s}_\mu} \ar[d]_{c \times c'} & \scrD(X,\texts,g) \ar[d]^{(-1) \times c'} \\
\scrC(X,\bar{\texts},g) \ar[r]^{\tilde{s}_{-\mu}} & \scrD(X,\bar{\texts},g). 
}
\]
This implies that the diagram
\[
\xymatrix{
    \scrB^\ast(X,\texts,g) \ar[r]^{s_\mu} \ar[d] & \scrE(X,\texts,g) \ar[d] \\
\scrB^\ast(X,\bar{\texts},g) \ar[r]^{s_{-\mu}} & \scrE(X,\bar{\texts},g)
}
\]
commutes, where the vertical maps are the isomorphisms induced from the charge conjugation.

Now we consider collections over $\spinc$ structures.
Fix $k \in \Z$ henceforth.
To avoid unnecessary complications, let us make the following assumption for a while:

\begin{assum}
\label{assum: no spin structure}
Assume that $\Spinc(X,k)$ does not contain a $\spinc$ structure that comes from a spin structure.
\end{assum}

Under this assumption, $\Z/2$ acts freely on $\Spinc(X)$.
In \cref{rmk: if spin exists} we describe how to remove this assumption.

Then the charge conjugation (together with the $(-1)$-multiplication) makes $\scrE(X) \to \scrB^\ast(X)$ a $\Z/2$-equivariant bundle, and the above observation implies that
$s : \scrB^\ast(X) \to \scrE(X)$ is a $\Z/2$-equivariant section.
Similarly, for $k>0$, $s^k : \scrB^\ast(X,k) \to \scrE(X,k)$ is also a $\Z/2$-equivariant section.

Set 
\begin{align*}
\scrB^\ast(X,k)' 
&:= (\scrB^\ast(X,k))/(\Z/2),\\
\scrE(X,k)' 
&:= (\scrE(X,k))/(\Z/2),\\
\circPi(X,k)' &:= (\Spincns(X,k) \times \circPi(X))/(\Z/2).
\end{align*}
The diffeomorphism group $\Diff^+(X)$ acts on these three spaces via the pull-back, since the pull-back commutes with both of the charge conjugation and the $(-1)$-multiplication.
Note that all of these spaces are equipped with the natural forgetful map onto $\Spincns(X,k)/\Conj$.
Also, $\Diff^+(X)$ acts on $\Spincns(X,k)/\Conj$ via the pull-back.
Thus we obtain fiber bundles associated with $E$ whose fibers are these objects, denoted by
\begin{align*}
\scrB^\ast(X,k)' \to
\scrB^\ast(E,k)' \to B,&\\
\scrE(X,k)' \to
\scrE(E,k)' \to B,&\\
\circPi(X,k)' \to \circPi(E,k)' \to B,&\\
\Spinc(X,k)/\Conj \to \Spinc(E,k)/\Conj \to B.&
\end{align*}

Now, the forgetful map $\scrB^\ast(X,k)' \to \Spincns(X,k)/\Conj$
gives rise to a surjection
\begin{align*}
\scrB^\ast(E,k)'
\to \Spincns(E,k)/\Conj,
\end{align*}
and similarly, we have surjections
\begin{align}
\scrE(E,k)'
\to \Spincns(E,k)/\Conj,\nonumber\\
\label{eq: map with connected fiber}
\circPi(E,k)'
\to \Spincns(E,k)/\Conj.
\end{align}
Evidently, all of these maps commute with the projections onto $B$.

Since the natural map $\scrE(X,k) \to \scrB^\ast(X,k)$ is $\Diff^+(X)$-equivariant and $\Z/2$-equivariant, this induces a map
\[
\scrE(E,k)' \to \scrB^\ast(E,k)'
\]
that commutes both with the projections onto $B$ and projections onto $\Spincns(E,k)/\Conj$.
Moreover, since $s^k$ is  $\Diff^+(X)$-equivariant too, and also $\Z/2$-equivariant as seen above, it induces a section
\[
s^k_{\mathrm{half}} : \scrB^\ast(E,k)' \to \scrE(E,k)'
\]
that commutes both with the projections onto $B$ and projections onto $\Spincns(E,k)/\Conj$.


In \cref{subsection Collection of Seiberg--Witten equations}, we considered a families perturbation $\sigma$.
Here we consider a counterpart $\sigma'$ which takes into account the charge conjugation.
Just like $\pi, p_{\scrB}, p_{\scrE}$, we have natural maps
\[
\pi': \circPi(E,k)' \to \scrR(E), \quad
p_{\scrB}' : \scrB^\ast(E,k)' \to \scrR(E),\quad
p_{\scrE}' : \scrE(E,k)' \to \scrR(E).
\]
Let $\widetilde{g}: B\rightarrow \scrR(E)$ be a section and let 
\[
\sigma' : \Spincns(E,k)/\Conj
\to \circPi(E,k)'
\]
be a section that fits into the following diagram:
\begin{align}
\label{eq: diagram half-tot section}
\begin{split}
\xymatrix{
\Spincns(E,k)/\Conj\ar[r]^-{\sigma'}\ar[d]& \circPi(E,k)'\ar[d] \\
B\ar[r]^{\widetilde{g}} & \scrR(E).}
\end{split}
\end{align}
We set
\[
\scrB^\ast(E,k,\sigma')'
:= (p_{\scrB}')^{-1}(\widetilde{g}),\quad
\scrE^\ast(E,k,\sigma)
:= (p_{\scrE}')^{-1}(\widetilde{g}).
\]
Then the restriction of $s^k_{\mathrm{half}}$ gives rise to a section
\[
(s^k_{\mathrm{half}})_{\sigma'} : \scrB^\ast(E,k,\sigma')' \to \scrE(E,k,\sigma')'.
\]

The section $(s^k_{\mathrm{half}})_{\sigma'}$ constructed now corresponds to the collection of the families Seiberg--Witten equations for the `half' of $\Spinc(X,k)$, namely $\Spinc(X,k)/\Conj$.

\begin{defi}
Set 
\[
\calM_{\sigma', \mathrm{half}} := \bigcup_{b \in B}((s^k_{\mathrm{half}})_{\sigma'})_{b}^{-1}(0).
\]
\end{defi}

Just as in \cref{lem: finiteness1}, we have:

\begin{lem}
If $B$ is compact, $\calM_{\sigma', \mathrm{half}}$ is also compact.
\end{lem}

We call $\calM_{\sigma', \mathrm{half}}$ the {\it family of half collections of moduli spaces} associated with $\sigma'$.
Thus, as described in \cref{defi: vn general coh}, we obtain a cohomology class
\begin{align}
\label{eq: half coh 1}
\frakM((s^k_{\mathrm{half}})_{\sigma'}) \in H^{k}(B; \Z/2)
\end{align}
through the virtual neighborhood technique.
If the section $(s^k_{\mathrm{half}})_{\sigma'}$ is nowhere vanishing over a subset $B' \subset B$, we have also a relative cohomology class
\begin{align}
\label{eq: half coh 2}
\frakM((s^k_{\mathrm{half}})_{\sigma'}; B') \in H^{k}(B, B'; \Z/2)
\end{align}
as in \cref{defi: vn general coh rel}.
Note that we exclude $\widetilde{g}$ in our notation since it is determined by $\sigma'$.

\begin{rmk}
\label{rem: why mod 2}
The pull-back of $\Z_E$ to $\Spincns(E,k)/\Conj$ along a natural map $\Spincns(E,k)/\Conj \to B$ does not necessarily give an appropriate local system over $\Spincns(E,k)/\Conj$ that coherently takes into account the collection of the orientations of the moduli spaces,
essentially by the same reason that the Seiberg--Witten invariant for $\texts$ equals that for $\bar{\texts}$ {\it only up to sign} in general.
Thus here we consider $\Z/2$-coefficient cohomology rather than $\Z_E$-coefficient.
\end{rmk}

Finally, we discuss how to remove \cref{assum: no spin structure}:

\begin{rmk}
\label{rmk: if spin exists}
Let us suppose that there are spin structures on $X$ with formal dimension $-k$.
In this case, let us decompose $\Spinc(X,k)$ into two parts:
\[
\Spinc(X,k) = \Spinc_0(X,k) \sqcup \Spinc_1(X,k).
\]
Here $\Spinc_0(X,k)$ denotes the subset of $\Spinc(X,k)$ consisting of $\spinc$ structures coming from spin structures, and $\Spinc_1(X,k)$ is the complement of $\Spinc_0(X,k)$ in $\Spinc(X,k)$.
Note that the $\Z/2$-action on $\Spinc(X,k)$ preserves this decomposition, and $\Z/2$ acts trivially on $\Spinc_0(X,k)$.
Thus we can carry out the construction of \cref{subsection Collection of Seiberg--Witten equations} for collections over $\Spinc_0(X,k)$, and the construction of this \lcnamecref{subsection Complex conjugation} for collections over $\Spinc_1(X,k)$.
Thus we may define the moduli space $\calM_{\sigma', \mathrm{half}}$ and the cohomology classes \eqref{eq: half coh 1} and \eqref{eq: half coh 2} as well.
\end{rmk}

\subsection{Inductive families perturbations}
\label{subsection Inductive sections}

In the constructions of the Seiberg--Witten characteristic classes,
we shall use a specific way to construct families perturbations.
We summarize it in this \lcnamecref{subsection Inductive sections} with basic examples.

Throughout this \lcnamecref{subsection Inductive sections}, we fix $k>0$ and let $X$ be a closed oriented smooth 4-manifold with $b^+(X) \geq k+2$.
Let $B$ be a CW complex, and let $X \to E \to B$ be a fiber bundle with structure group $\Diff^+(X)$.
For an $i$-cell $e \subset B$, let $\varphi_{e} : D_{e}^{i} \to \bar{e} \subset B$ denote the characteristic map of $e$, where $D^{i}_{e}$ is the standard $i$-dimensional disk indexed by $e$.

First, as in \cite[Subsection~6.1]{K21},
we inductively construct a section
\[
\sigma^{(i)} : B^{(i)} \to \circPi(E)|_{B^{(i)}}
\]
for $i \leq k+1$ as follows.
Choosing a generic point in $\circPi(E_{b})$ for each $b \in B^{(0)}$, we have $\sigma^{(0)} : B^{(0)} \to \circPi(E)|_{B^{(0)}}$.
Assume that we have constructed $\sigma^{(i-1)} : B^{(i-1)} \to \circPi(E)|_{B^{(i-1)}}$ for $i \leq k+1$ such that $\calM_{\sigma^{(i-1)}} = \emptyset$, where $\calM_{\sigma^{(i-1)}}$ is the family of total collections of moduli spaces associated with $\sigma^{(i-1)}$, introduced in \eqref{eq: SWtot moduli}.
Let $e$ be an $i$-cell of $B$.
Since the pull-back bundle $\varphi_{e}^{\ast} \circPi(E) \to D_{e}^{i}$ under the characteristic map is trivial, there is a trivialization 
\[
\psi_{e} : \varphi_{e}^{\ast} \circPi(E) \to D_{e}^{k} \times \circPi(X).
\]
Let $p : \del D_{e}^{i} \times \circPi(X) \to \circPi(X)$ denote the projection.
Since $b^{+}(X) \geq k+2$, we can extend the continuous map 
\[
p \circ \psi_{e} \circ (\varphi_{e}|_{\del D_{e}^{i}})^{\ast}\sigma^{(i-1)}
: \del D_{e}^{i} \to (\varphi_{e}|_{\del D_{e}^{i}})^{\ast} \circPi(E) \to \del D_{e}^{i} \times \circPi(X) \to \circPi(X)
\]
to a map from $D_{e}^{i}$ into $\circPi(X)$, rather than $\Pi(X)$, which corresponds to choosing a family of self-dual 2-forms avoiding the wall.
This extended map gives rise to a section of $\circPi(E)|_{\bar{e}} \to \bar{e}$.
We may choose the extension generically so that $\calM_{\sigma^{(i)}|_{e}} = \emptyset$ for $i<k$ (see \cref{rmk: generic ext} below).
Thus we obtain a section $\sigma^{(i)} : B^{(i)} \to \Pi(E)|_{B^{(i)}}$ for $i \leq k+1$, which satisfies that $\calM_{\sigma^{(i)}} = \emptyset$ if $i<k$.
We call such $\sigma^{(\bullet)}$ an {\it inductive families perturbation} (called an inductive section in \cite{K21}).

\begin{rmk}
\label{rmk: generic ext}
Let us explain why one can extend a {\it continuous} map
\[
\Phi := p \circ \psi_{e} \circ (\varphi_{e}|_{\del D_{e}^{i}})^{\ast}\sigma^{(i-1)}
: \del D_{e}^{i} \to \circPi(X)
\]
to a generic map 
\[
D_{e}^{i} \to \circPi(X)
\]
so that the parameterized moduli space over $D_{e}^{i}$ is empty for $i < k$.
First, recall that the parameterized moduli space over $\del D_{e}^{i}$ was supposed to be empty.
Take a smooth map $\Phi' : \del D_{e}^{i} \to \circPi(X)$ that is close enough to $\Phi$.
Since $\Phi'$ is smooth, by a standard argument one can find a smooth and generic extension of $\Phi'$ to a map $\tilde{\Phi}' : D_{e}^{i} \to \circPi(X)$ for which the parameterized moduli space is empty.

On the other hand, since the emptiness of the moduli space is an open condition,
we may suppose that there is a continuous homotopy $\Phi_t$ between $\Phi$ and $\Phi'$ such that the moduli space is empty for all $\Phi_t$.
Pick a collar neighborhood $N(\del D_{e}^{i})$ of $\del D_{e}^{i}$ in $D_{e}^{i}$ and define a continuous extension of $\Phi$ to a map $\hat{\Phi} : N(\del D_{e}^{i}) \to \circPi(X)$ using $\Phi_t$ by assigning $t$ to the radial coordinate.
Gluing $\hat{\Phi}$ with (rescaled) $\tilde{\Phi}'$, we may obtain an extension of $\Phi$ to $D_{e}^i$ for which the parameterized moduli space over $D_{e}^{i}$ is empty for $i < k$.
\end{rmk}

The above way to construct an inductive section shall be used to define the total characteristic class $\SWbbtot^\bullet$.
We also describe a similar construction to define the half-total characteristic class $\SWbbhalftot^\bullet$.
In this case, we need to construct a section
\[
\sigma'^{(i)} : \left(\Spincns(E,k)/\Conj\right)^{(i)}
\to \circPi(E,k)'
\]
for $i \leq k+1$ inductively. 
First, note that $\Spincns(E,k)/\Conj$ is a covering space of $B$ with fiber $\Spincns(X,k)/\Conj$ with the natural projection $\Spincns(E,k)/\Conj \to B$.
We equip $\Spincns(E,k)/\Conj$ with the CW complex structure induced from that of $B$: the lifted CW complex structure along the map $\Spincns(E,k)/\Conj \to B$.
For an $i$-cell $e \subset B$ with the characteristic map $\varphi_e : D^i_e \to B$, let
\[
(e_{[\fraks]})_{[\fraks] \in \Spincns(X,k)/\Conj}
\]
be the corresponding lifted $i$-cells of $\Spincns(E,k)/\Conj$ with the characteristic maps
\[
\hat{\varphi}_{e_{[\fraks]}} : D^i_e \to \Spincns(E,k)/\Conj,
\]
which are lifts of $\varphi_e$.
Note that the fiber of the natural map $\circPi(E,k)' \to \Spincns(E,k)/\Conj$ is given by $\circPi(X)$.
Also, a trivialization $\psi_{e} : \varphi_{e}^{\ast} \circPi(E) \to D_{e}^{k} \times \circPi(X)$ induces trivializations
\[
\hat{\psi}_{e_{[\fraks]}} : \hat{\varphi}_{e_{[\fraks]}}^{\ast} \circPi(E,k)' \to D_{e}^{k} \times \circPi(X).
\]
Thus we may repeat the above construction of an inductive families perturbation using $(e_{[\fraks]}, \hat{\varphi}_{e_{[\fraks]}}, \hat{\psi}_{e_{[\fraks]}})$ in place of $(e, \varphi_e, \psi_e)$, and obtain a section $\sigma'^{(i)} : \Spincns(E,k)/\Conj
\to \circPi(E,k)'$ with the property that 
$\calM_{\sigma'^{(i)}, \mathrm{half}} = \emptyset$ if $i<k$.
We call this $\sigma'^{(i)}$ an {\it inductive half-total families perturbation}.

We give a few examples of (inductive) half-total families perturbations.

\begin{ex}
\label{ex: trivial monodromy}
Let $X$ and $k$ be as above, and let $X \to E \to B$ be a $\Diff^+(X)$-bundle with trivial monodromy on $\Spincns(X,k)$.
Then the covering space $\Spincns(E,k)/\Conj \to B$ is just the product with an obvious projection: 
$B \times \Spincns(X,k)/\Conj \to B$.
Therefore a section $\sigma' : \Spincns(E,k)/\Conj \to \circPi(E,k)'$ corresponds to a collection of sections
\[
\left(\sigma_{[\fraks]} : B \to \circPi(E)\right)_{[\fraks] \in \Spincns(X)/\Conj}.
\]
There is no constraint on this collection of sections;
for example, we may take all $\sigma_{[\fraks]}$ to be a common section $\sigma : B \to \circPi(E)$.
\end{ex}

\begin{ex}
\label{ex: perturbation single mapping torus 1}
Let $X$ and $k$ be as above, and let $f \in \Diff^+(X)$.
Let $X \to X_f \to S^1$ denote the mapping torus of $f$ with fiber $X$.
Then a continuous section 
\[
\sigma' : \Spincns(X_f,k)/\Conj \to \circPi(X_f,k)'
\]
corresponds to a family of continuous sections
\[
\left(\sigma_{t} : (\Spincns(X,k)/\Conj) \to \circPi(X,k)' \right)_{t \in [0,1]}
\]
that satisfies the condition $f^\ast \circ \sigma_{0}=\sigma_1 \circ f^\ast$.

We can describe such $\sigma'$ in terms of maps to $\circPi(X)$, once we fix a section 
\[
\tau : \Spincns(X,k)/\Conj \to \Spincns(X,k)
\]
of the surjection $\Spincns(X,k) \to \Spincns(X,k)/\Conj$.
Let
\begin{align*}
\left(\sigma_{[\fraks]} : [0,1] \to \circPi(X)
\right)_{[\fraks] \in \Spincns(X,k)/\Conj}
\end{align*} 
be a collection of continuous maps which satisfies that, 
 for $[\fraks] \in \Spincns(X,k)/\Conj$,
\[
f^\ast\sigma_{[\fraks]}(0) =
\left\{
\begin{array}{ll}
    \sigma_{f^\ast[\fraks]}(1) \quad &\text{if} \quad f^\ast\tau([\fraks]) \in \tau(\Spincns(X,k)/\Conj),\\
    -\sigma_{f^\ast[\fraks]}(1) \quad &\text{if} \quad f^\ast\tau([\fraks]) \notin \tau(\Spincns(X,k)/\Conj).
\end{array}
\right.
\]
Here $ -\sigma_{f^\ast[\fraks]}(1)$ means keeping the metric unchanged and multiplying the perturbation by $-1$.
Then it is easy to check that the collection $(\sigma_{[\fraks]})_{[\fraks]}$ gives rise to a half-total families perturbation \[\sigma' : \Spincns(X_f,k)/\Conj \to \circPi(X_f,k)'\] for the mapping torus $X_f$.
Not as in the above \cref{ex: trivial monodromy}, we may not take all $\sigma_{[\fraks]}$ to be a common map $[0,1] \to \circPi(X)$ in general.
\end{ex}

\begin{ex}
\label{ex: perturbation multi map torus}
Let $X$ and $k>0$ be as above, and let $f_1, \ldots, f_k \in \Diff^+(X)$ and suppose that they mutually commute.
Let $X \to E \to T^k$ denote the multiple mapping torus of $f_1, \ldots, f_k$ with fiber $X$.
The above \cref{ex: perturbation single mapping torus 1} can be easily generalized to the multiple mapping torus as follows.

For $i_1, \ldots, i_l \in \{1, \ldots, k\}$ with $i_1<\cdots<i_l$, let $f_{i_1, \ldots, i_l}^\ast$ denote the composition of pull-backs $f_{i_1}^\ast, \ldots,  f_{i_l}^\ast$.
For $\mathbf{t} = (t_1, \ldots, t_n) \in [0,1]^k$  with $t_{i_1}=\dots=t_{i_l}=0$, define
$\overline{\mathbf{t}}^{i_1, \ldots, i_l} \in [0,1]^k$
by 
$\overline{\mathbf{t}}^{i_1, \ldots, i_l}= (t_1',\dots, t_k')
$,
where 
\begin{align*}
t_j'
=
\left\{
\begin{array}{ll}
t_j \quad &\text{if} \quad j \notin \{i_1, \ldots, i_l\},\\
1 \quad &\text{if} \quad j \in \{i_1, \ldots, i_l\}.
\end{array}
\right.
\end{align*}
Fix a section $\tau : \Spincns(X,k)/\Conj \to \Spincns(X,k)$.
Let
\begin{align*}
\left(\sigma_{[\fraks]} : [0,1]^k \to \circPi(X)
\right)_{[\fraks] \in \Spincns(X,k)/\Conj}
\end{align*}
be a collection of continuous maps.
Suppose that, for all
\begin{itemize}
\item $[\fraks] \in \Spincns(X,k)/\Conj$,
\item $i_1, \ldots, i_l \in \{1, \ldots, k\}$ with $i_1<\cdots<i_l$ and,
\item $\mathbf{t} = (t_1, \ldots, t_n) \in [0,1]^k$ with $t_{i_1}=\dots=t_{i_l}=0$,
\end{itemize}
we have
\begin{equation}
\label{eq: peruturb general mapping tori}
f_{i_1, \ldots, i_l}^\ast\sigma_{[\fraks]}(\mathbf{t}) =
\left\{
\begin{array}{ll}
    \sigma_{f_{i_1, \ldots, i_l}^\ast[\fraks]}(\overline{\mathbf{t}}^{i_1, \ldots, i_l}) \quad &\text{if} \quad f_{i_1, \ldots, i_l}^\ast\tau([\fraks]) \in \tau(\Spincns(X,k)/\Conj),\\
    -\sigma_{f_{i_1, \ldots, i_l}^\ast[\fraks]}(\overline{\mathbf{t}}^{i_1, \ldots, i_l}) \quad &\text{if} \quad f_{i_1, \ldots, i_l}^\ast\tau([\fraks]) \notin \tau(\Spincns(X,k)/\Conj).
\end{array}
\right.
\end{equation}
Then the collection $(\sigma_{[\fraks]})_{[\fraks]}$ gives rise to a half-total families perturbation
\[
\sigma' : \Spincns(E,k)/\Conj \to \circPi(E,k)'
\]
for the multiple mapping torus $E$.

If we equip $T^k$ with the standard cell structure, we may construct such a collection $(\sigma_{[\fraks]})_{[\fraks]}$ inductively from $0$-cells, and it gives rise to an inductive half-total families perturbation $\sigma'$.
\end{ex}

\subsection{The characteristic classes $\SWbbtot$ and $\SWbbhalftot$}
\label{subsec cochain SWtot}

In this \lcnamecref{subsec cochain SWtot}, we shall define the characteristic classes $\SWbbtot$ and $\SWbbhalftot$.

To help our intuition,
first we define the 0-th degree part $\SWbbtot^0(X)$ and $\SWbbhalftot^0(X)$, which are numerical invariants.

\begin{defi}
Let $X$ be a closed oriented smooth 4-manifold with $b^+(X) \geq 2$.
Fix a homology orientation of $X$.
Let $\tau : \Spincns(X,k)/\Conj \to \Spincns(X,k)$ be a section.
We define
\begin{align*}
\SWbbtot^0(X) &= \sum_{\fraks \in \Spinc(X,0)}\SW(X,\fraks) \in \Z,\\
\SWbbhalftot^0(X) &= \sum_{\fraks \in \tau(\Spincns(X,0)/\Conj)}\SW(X,\fraks) \in \Z/2.    
\end{align*}
Here $\SW(X,\fraks)$ denotes the Seiberg--Witten invariant or its mod 2 reduction.
\end{defi}

Note that the $\Z/2$-valued invariant $\SWbbhalftot^0(X)$ is independent of $\tau$ because of
the formula $\SW(X,\fraks) = \pm\SW(X,\bar{\fraks})$ in $\Z$.

\begin{ex}
\label{ex: minimal algebraic surfaces of general type}
A typical example of a 4-manifold $X$ with $\SWbbhalftot^0(X) \neq 0$ is a minimal algebraic surface of general type, since it has a unique Seiberg--Witten basic class up to conjugation.
One may produce more examples by the fiber sum operation, which shall be discussed in \cref{construction of 4-manifolds}.
(See, for example, \cref{pro: fiber sum has nontrivial SW-tot}.)
\end{ex}

Now we define cohomological generalizations of the above numerical invariants, which are the characteristic classes $\SWbbtot^\bullet$ and $\SWbbhalftot^\bullet$ with $\bullet>0$.
Henceforth, throughout this \lcnamecref{subsec cochain SWtot}, let $k>0$ and let $X$ be a closed oriented smooth 4-manifold with $b^+(X) \geq k+2$ and with a fixed homology orientation.
Let $B$ be a CW complex, and let $X \to E \to B$ be a fiber bundle with structure group $\Diff^+(X)$.
We shall define cohomology classes $\SWbbtot^k(E) \in H^k(B, \Z_E)$ and $\SWbbhalftot^k(E) \in H^k(B;\Z/2)$, where $\Z_E$ is the local system with fiber $\Z$ defined in \cref{subsection Collection of Seiberg--Witten equations}, which is determined by the monodromy action for $E$ on the homology orientation.

Let
\[
C_\ast(B;\Z_E),\quad C^\ast(B;\Z_E),\quad
C_\ast(B),\quad
C^\ast(B)
\]
denote the cellular chain and cochain complexes with local coefficients $\Z_E$, and the cellular chain and cochain complexes with coefficient $\Z/2$, respectively.
The boundary operator and coboundary operator are denoted by $\del$ and $\delta$ for any (co)chain complex.
Here we adopt the following classical model due to Steenrod~\cite{Steenrod43} as the definition of (co)chain complex with local coefficients:
For each cell $e$, we choose a reference point $x(e)$ in $e$.
Let $\Z_E(e)$ denote the fiber of $\Z_E$ over $x(e)$. 
The $i$-th chain group $C_i(B;\Z_E)$ is defined be the set of formal finite sums $\Sigma_{e}a_{e}e$, where $e$ runs over $i$-cells of $B$ and $a_{e} \in \Z_E(e)$.
The cochain group $C^i(B;\Z_E)$ is defined to be the set of functions that send each $i$-cell $e$ to an element of $\Z_E(e)$.
The (co)boundary operator is defined by twisting the ordinal definition by the automorphism of $\Z$ corresponding to the paths connecting the respective reference points.

The first step of the construction of $\SWbbtot^k(E)$ is to define a $k$-cochain of $B$, which will be shown to be a cocycle and whose cohomology class is $\SWbbtot^k(E)$.
Let $e$ be a $k$-cell of $B$, and let $\varphi_e : D^k_e \to B^{(k)}$ be the characterisric map of $e$.
Take an inductive families perturbation
\[
\sigma = \sigma^{(k)} : B^{(k)} \to \circPi(E)|_{B^{(k)}}
\]
over the $k$-skeleton of $B$,
explained in \cref{subsection Inductive sections}.
Consider the pull-back section
\[
\varphi_e^\ast \sigma : D^k_e \to \varphi_e^\ast\circPi(E).
\]
Note that $\varphi_e^\ast\circPi(E) \subset \circPi(\varphi_e^\ast E)$, and hence $\varphi_e^\ast \sigma$ gives a section of $\circPi(\varphi_e^\ast E) \to D^k_e$.
Applying the construction in \cref{subsection Collection of Seiberg--Witten equations} to  this pull-back section $\varphi_e^\ast \sigma$, we obtain a collection of Fredholm sections
\[
s^k_{\varphi_e^\ast \sigma} : 
\scrB^\ast(\varphi_e^\ast E,k,\varphi_e^\ast\sigma) \to \scrE(\varphi_e^\ast E,k,\varphi_e^\ast\sigma)
\]
and the parameterized moduli space $\calM_{\varphi_e^\ast \sigma}$, parameterized over $D^k_e$.
Since any inductive families perturbation was constructed to satisfy $\calM_{\sigma^{(i)}} = \emptyset$ for $i<k$, we have $\calM_{\varphi_e^\ast \sigma|_{\del D^k_e}} = \emptyset$.
Thus, as explained in \cref{subsection Collection of Seiberg--Witten equations}, we obtain a relative cohomology class 
\[
\frakM(s^k_{\varphi_e^\ast \sigma}; \del D^k_e) \in H^{k}(D^k_e, \del D^k_e; \varphi_e^\ast \Z_E),
\]
which is an application of Equation~\eqref{eq: coh compact base tot rel} to $(B,B',\sigma) = (D^k_e, \del D^k_e,\varphi_e^\ast \sigma)$.

Note that there is a canonical isomorphism
$H^{k}(D^k_e, \del D^k_e; \varphi_e^\ast \Z_E) \cong \Z_E(e)$
via $\varphi_e$, without ambiguity of $\Aut(\Z) \cong \Z/2$.
This isomorphism gives us a pairing 
\[
\left<-, - \right> : H^{k}(D^k_e, \del D^k_e; \varphi_e^\ast \Z_E)
\otimes
H_{k}(D^k_e, \del D^k_e; \Z)
\to \Z_E(e).
\]

Instead of an inductive families perturbation $\sigma$,
if we take an inductive half-total families perturbation
\[
\sigma'=\sigma'^{(k)} : \left(\Spincns(E,k)/\Conj\right)^{(k)}
\to \circPi(E,k)',
\]
the `half-total version' of the above argument works:
For each $k$-cell $e$,
the pull-back under $\varphi_e$ gives rise to a section
\[
\varphi_e^\ast \sigma' : \varphi_e^\ast (\Spincns(E,k)/\Conj) \to \varphi_e^\ast (\circPi(E)').
\]
Here we identified $\varphi_e^\ast(\Spincns(E,k)/\Conj)$ with $\Spincns(\varphi_e^\ast E,k)/\Conj$.
Applying the construction of \cref{subsection Complex conjugation} to this pull-back section,
we may obtain a relative cohomology class
\[
\frakM((s_{\mathrm{half}}^k)_{\varphi_e^\ast \sigma'}; \del D^k_e) 
\in H^{k}(D^k_e, \del D^k_e; \Z/2),
\]
which is an application of Equation~\eqref{eq: half coh 2} to $(B,B',\sigma') = (D^k_e, \del D^k_e,\varphi_e^\ast \sigma')$.

\begin{defi}
\label{defi: cochain}
For an inductive families perturbation $\sigma$ and an inductive half-total families perturbation $\sigma'$, we define cochains
\begin{align*}
\SWcaltot^k(E, \sigma) &\in C^k(B;\Z_E),\\  \SWcalhalftot^k(E, \sigma') &\in C^k(B)
\end{align*}
by
\begin{align*}
\SWcaltot^k(E, \sigma)(e)
&=\langle\frakM(s^k_{\varphi_e^\ast \sigma}; \del D^k_e), [(D^k_e, \del D^k_e)]\rangle \in \Z_E(e),\\
\SWcalhalftot^k(E, \sigma')(e)
&=\langle\frakM((s_{\mathrm{half}}^k)_{\varphi_e^\ast \sigma'}; \del D^k_e), [(D^k_e, \del D^k_e)]\rangle \in \Z/2
\end{align*}
for each $k$-cell $e$ of $B$, respectively.
\end{defi}

In subsequent \cref{subsection Well-definedness and naturality}, we shall prove:

\begin{pro}
\label{SWtot cocycle}
The cochains $\SWcaltot^k(E, \sigma), \SWcalhalftot^k(E, \sigma')$ are cocycles.
\end{pro}

\begin{pro}
\label{SWtot indep}
The cohomology classes $[\SWcaltot^k(E, \sigma)] \in H^k(B;\Z_E)$, $[\SWcalhalftot^k(E, \sigma')] \in H^k(B;\Z/2)$ are independent of the choice of $\sigma$ and $\sigma'$ respectively.
\end{pro}

Thus we arrive at the following definition:

\begin{defi}
Define 
\begin{align*}
\SWbbtot^k(E) &\in H^k(B;\Z_E),\\  \SWbbhalftot^k(E) &\in H^k(B;\Z/2)
\end{align*}
by
\begin{align*}
\SWbbtot^k(E) &= [\SWcaltot^k(E, \sigma)],\\
\SWbbhalftot^k(E) &= [\SWcalhalftot^k(E, \sigma')],
\end{align*}
where $\sigma : B^{(k)} \to \circPi(E)|_{B^{(k)}}$ is an inductive families perturbation, and $\sigma' : \left(\Spincns(E,k)/\Conj\right)^{(k)}
\to \circPi(E,k)'$ is an inductive half-total families perturbation, respectively.
We define also 
\begin{align*}
\SWbbtot^k(X) &\in H^k(\BDiff^+(X);\Z_{\EDiff^+(X)}),\\  \SWbbhalftot^k(X) 
&\in H^k(\BDiff^+(X);\Z/2)
\end{align*}
by
\begin{align*}
\SWbbtot^k(X) &:= \SWbbtot^k(\EDiff^+(X)),\\ 
\SWbbhalftot^k(X) &:= \SWbbhalftot^k(\EDiff^+(X)).
\end{align*}
\end{defi}


\begin{rmk}
\label{rem: relation with Rub and K}
There are two preceding constructions deeply related to $\SWbbtot^k(X)$ and $\SWbbhalftot^k(X)$.
Given (the isomorphism class of) a $\spinc$ structure  $\fraks$ with $d(\fraks)=-k$, let $\Diff^+(X,\fraks)$ denote the group of diffeomorphisms that preserve orientation and $\fraks$.
A characteristic class $\SWbb(X,\fraks) \in H^k(\BDiff^+(X, \fraks))$ was constructed by the first author \cite{K21}, which is defined by counting moduli spaces for $\fraks$ over $k$-cells.
A significant difference between this paper from \cite{K21} is that $\SWbbtot^k(X)$ and $\SWbbhalftot^k(X)$ can be defined over the whole $\BDiff^+(X)$ in other words, defined for all oriented fiber bundles with fiber $X$. On the other hand,  $\SWbb(X,\fraks)$ in \cite{K21} can be defined only for fiber bundles whose monodromy preserve $\fraks$.

The idea of the ``total" construction is inspired by Ruberman's construction of his total Seiberg--Witten invariant of diffeomorphisms $\mathrm{SW}_{\rm{tot}}(f,\fraks)$ \cite{Rub01}.
The invariant $\mathrm{SW}_{\rm{tot}}(f,\fraks)$ is a numerical invariant, and our invariant $\SWbbtot^k(X)$ is, roughly, a cohomological refinement of  $\mathrm{SW}_{\rm{tot}}(f,\fraks)$.
More precisely, although $\mathrm{SW}_{\rm{tot}}(f,\fraks)$ was defined by summing up the moduli spaces over the orbit $\{(f^n)^\ast \fraks\}_{n \in \Z}$, to define an invariant $\SWbbtot^k(X)$ independent of $\fraks$, we shall use the sum over all $\spinc$ structures $\fraks$ with $d(\fraks)=-k$. This is why we need to consider also $\SWbbhalftot^k(X)$, not as in Ruberman's situation \cite{Rub01}: in most cases, the contribution from the orbit $\{(f^n)^\ast \fraks\}_{n \in \Z}$ cancels out that from the conjugate orbit $\{(f^n)^\ast \bar{\fraks}\}_{n \in \Z}$ over $\Z/2$.
\end{rmk}


\subsection{Well-definedness}
\label{subsection Well-definedness and naturality}

In this \lcnamecref{subsection Well-definedness and naturality}, we prove \cref{SWtot cocycle} (cocycle) and  \cref{SWtot indep} (well-definedness).
Some parts of the proofs are simple adaptations of arguments in \cite[Subsection~6.2]{K21}, and in such cases we just indicate corresponding propositions of \cite{K21}.

First, we shall note naturality results.
As in \cref{subsec cochain SWtot}, let $k>0$ and let $X$ be a closed oriented smooth 4-manifold with $b^+(X) \geq k+2$. 
Let $B$ be a CW complex, and let $X \to E \to B$ be a fiber bundle with structure group $\Diff^+(X)$.
The following is a straightforward consequence of 
\cref{lem: first naturality}:

\begin{lem}
\label{lem: frakM calSW}
In the setup of \cref{defi: cochain}, suppose that $B$ is compact.
Then we have 
\begin{align*}
\frakM(s^k_{\sigma}; B^{(k-1)})
&= \SWcaltot^k(E,\sigma),\\
\frakM((s^k_{\mathrm{half}})_{\sigma'}; B^{(k-1)})
&= \SWcalhalftot^k(E,\sigma')
\end{align*}
in $H^k(B^{(k)}, B^{(k-1)};\Z_E) = C^k(B;\Z_E)$ and in $H^k(B^{(k)}, B^{(k-1)};\Z/2) = C^k(B)$, respectively.
\end{lem}

\begin{proof}
Using \cref{lem: first naturality}, we can adapt the proof of \cite[Lemma~6.13]{K21}:
we may repeat the proof using $\SWcaltot^k(E,\sigma)$ and $\SWcalhalftot^k(E,\sigma')$, in place of $\mathcal{A}(E,\sigma)$ in the proof of \cite[Lemma~6.13]{K21}.
\end{proof}

This \lcnamecref{lem: frakM calSW} implies the naturality for cochains:

\begin{pro}
\label{prop: cochain natural}
Let $B'$ be a CW complex and $f : B' \to B$ be a cellular map.
Then we have 
\begin{align*}
f^\ast\SWcaltot^k(E,\sigma)
&= \SWcaltot^k(f^\ast E,f^\ast \sigma),\\
f^\ast\SWcalhalftot^k(E,\sigma)
&= \SWcalhalftot^k(f^\ast E,f^\ast \sigma)
\end{align*}
in $C^k(B;\Z_{f^\ast E})$ and in $C^k(B)$, respectively.
\end{pro}

\begin{proof}
Using \cref{lem: first naturality} and \cref{lem: frakM calSW},
we can just adapt the proof of \cite[Proposition~6.14]{K21}.
\end{proof}

It directly follows  from \cref{prop: cochain natural} and the cellular approximation theorem that:

\begin{cor}
\label{SWtot functoriality}
Let $B'$ be a CW complex and $f : B' \to B$ be a continuous map.
Then
\[
f^\ast\SWbbtot^k(E) = \SWbbtot^k(f^\ast E),\quad
f^\ast\SWbbhalftot^k(E) = \SWbbhalftot^k(f^\ast E)
\]
hold in $H^k(B'; f^\ast \Z_E) (\cong H^k(B'; \Z_{f^\ast E}))$ and in $H^k(B';\Z/2)$, respectively.
\end{cor}

Note that,
applying \cref{SWtot functoriality} to the identity map, we have that the cohomology class $\SWbbtot^k(E)$ is independent of the choice of CW structure of $B$.
It also follows from a formal argument that $\SWbbtot^k(E)$ can be defined for any oriented smooth fiber bundle with fiber $X$ over an arbitrary topological space, not only over a CW complex (see \cite[Lemma  6.8]{K21}).

\begin{proof}[Proof of \cref{SWtot cocycle}]
This corresponds to \cite[Proof of Proposition~6.2]{K21}.
As the statement for $\SWcaltot$ is proven by a rather direct adaptation of \cite{K21}, here we give a proof of the statement for $\SWcalhalftot$ for the reader's convenience.
As noted in \cref{subsection Inductive sections}, we can extend the construction of an inductive half-total families perturbation to the $(k+1)$-skeleton, thanks to the condition that $b^+(X) \geq k+2$.
Let 
\[
\sigma' : \left(\Spincns(E,k)/\Conj\right)^{(k+1)}
\to \circPi(E,k)'
\]
be an extension of $\sigma'$ to the $(k+1)$-skeleton.

Let $e$ be a $(k+1)$-cell of $B$  with the characteristic map $\varphi_e : D^{k+1}_e \to B$.
Fix a homeomorphism between $D^{k+1}_e$ and $\Delta^{k+1}$, the standard $(k+1)$-simplex.
Let $D^{k+1}_e$ be equipped with a CW structure using the standard CW structure on $\Delta^{k+1}$ via this homeomorphism.
By cellular approximation, we can suppose that $\varphi_e$ is cellular.

We rewrite $\delta \SWcalhalftot(E, \sigma')(e)$ as follows.
First, it follows from a direct computation with \cref{prop: cochain natural} that
\begin{align}
\label{eq: proof cocycle 1}
\delta \SWcalhalftot(E, \sigma')(e)
= \SWcalhalftot(\varphi_e^\ast E, \varphi_e^\ast\sigma')(\del \Delta^{k+1}).
\end{align}
Next, it follows from \cref{lem: frakM calSW} that 
\begin{align}
\label{eq: proof cocycle 2}
\frakM((s^k_{\mathrm{half}})_{\varphi_e^\ast\sigma'|_{(\Delta^{k+1})^{(k)}}}; (\Delta^{k+1})^{(k-1)})
= \SWcalhalftot(\varphi_e^\ast E, \varphi_e^\ast \sigma'),
\end{align}
where $\varphi_e^\ast\sigma'|_{(\Delta^{k+1})^{(k)}}$ denotes the restricted section
\[
\varphi_e^\ast\sigma'|_{(\Delta^{k+1})^{(k)}}:
\left(\Spincns(\varphi_e^\ast E,k)/\Conj\right)^{(k+1)}|_{(\Delta^{k+1})^{(k)}}
\to \circPi(\varphi_e^\ast E,k)'|_{(\Delta^{k+1})^{(k)}}.
\]
As a last piece of the proof, because of the vanishing of the moduli space over the $(k-1)$-skeleton,
we may obtain a cohomology class 
\[
\frakM((s^k_{\mathrm{half}})_{\varphi_e^\ast\sigma'}; (\Delta^{k+1})^{(k-1)}) \in H^{k}(\Delta^{k+1}, (\Delta^{k+1})^{(k-1)}; \Z/2).
\]
Let
\[
i : ((\Delta^{k+1})^{(k)}, (\Delta^{k+1})^{(k-1)})
\hookrightarrow (\Delta^{k+1}, (\Delta^{k+1})^{(k-1)})
\]
denote the inclusion.
Then we have 
\begin{align}
\label{eq: proof cocycle 3}
i^\ast \frakM((s^k_{\mathrm{half}})_{\varphi_e^\ast\sigma'}; (\Delta^{k+1})^{(k-1)})
= \frakM((s^k_{\mathrm{half}})_{\varphi_e^\ast\sigma'|_{(\Delta^{k+1})^{(k)}}}; (\Delta^{k+1})^{(k-1)}).
\end{align}

Now, a straightforward calculation using \eqref{eq: proof cocycle 1}, \eqref{eq: proof cocycle 2}, \eqref{eq: proof cocycle 3} with $i_\ast \del \Delta^{k+1}=0$ verifies that
$\delta \SWcalhalftot(E, \sigma')(e)
= 0$, which proves that $\SWcalhalftot(E, \sigma')$ is a cocycle.
\end{proof}

\begin{proof}[Proof of \cref{SWtot indep}]
This is an adaptation of \cite[Proof of Theorem~6.4]{K21}.
Again we just give a proof of the statement for $\SWcalhalftot$.
For $j=0,1$, let
\[
\sigma_j' = \sigma_j'^{(k)}: \left(\Spincns(E,k)/\Conj\right)^{(k)}
\to \circPi(E,k)'
\]
be inductive half-total families perturbations.
Set 
\[
\SWcal_j = \SWcalhalftot^k(E,\sigma_j').
\]
Let $p : B \times [0,1] \to B$ be the projection.
For each cell $e$ of $B$ equipped with a characteristic map $\varphi_e : D^k_e \to B$, we fix a trivialization
\[
\tilde{\psi}_e : (\varphi_e \times \id)^\ast p^\ast \circPi(E) \to D^{\dim e} \times [0,1] \times \circPi(X).
\]
As in \cref{subsection Inductive sections}, we regard $\Spincns(E,k)/\Conj$ as a covering space over $B$ with fiber $\Spincns(X,k)/\Conj$.
Let $\hat{p} : (\Spincns(E,k)/\Conj) \times [0,1] \to \Spincns(E,k)/\Conj$ be the projection.
Recall that $\varphi_e$ lifts to characteristic maps $(\hat{\varphi}_{e_{[\fraks]}})_{[\fraks] \in \Spincns(E,k)/\Conj}$ of cells $(e_{[\fraks]})_{[\fraks] \in \Spincns(E,k)/\Conj}$ of $\Spincns(E,k)/\Conj$ which are lifts of $e$.
Also, $\tilde{\psi}_e$ induces trivializations
\[
\hat{\tilde{\psi}}_{e_{[\fraks]}} : (\hat{\varphi}_{e_{[\fraks]}} \times \id)^\ast \hat{p}^\ast \circPi(E,k)' \to D^{\dim e} \times [0,1] \times \circPi(X).
\]

We inductively construct a section
\[
\tilde{\sigma}' = \tilde{\sigma}'^{(k)}: (\Spincns(E,k)/\Conj)^{(k)} \times [0,1] \to \hat{p}^\ast \circPi(E,k)'
\]
as follows.
Assume that we have constructed a section
$\tilde{\sigma}'^{(i-1)}$ for $i \leq k$ such that 
\[
\tilde{\sigma}'^{(i)}|_{(\Spincns(E,k)/\Conj)^{(i)} \times {\{j\}}}
= \sigma_j'^{(i)}
\]
for $j=0,1$ and that $\calM_{\tilde{\sigma}'^{(i-1)}, \mathrm{half}} = \emptyset$.
(This can be easily done for $i=0$.)
For an $i$-cell $e$ of $B$, we take a section 
\[
\tilde{\sigma}'^{(i)} : (\Spincns(E,k)/\Conj)|_{\bar{e}} \times [0,1] \to \hat{p}^\ast \circPi(E,k)'
\]
so that:
\begin{itemize}
\item $\tilde{\sigma}'^{(i)}|_{(\Spincns(E,k)/\Conj)|_{\bar{e}} \times {\{j\}}}
= \sigma_j'^{(i)}|_{(\Spincns(E,k)/\Conj)|_{\bar{e}}}$
for $j=0,1$,
\item 
$\tilde{\sigma}'^{(i)}|_{(\Spincns(E,k)/\Conj)|_{\bar{e} \setminus \{e\}} \times [0,1]}
= \sigma'^{(i-1)}|_{(\Spincns(E,k)/\Conj)|_{\bar{e} \setminus \{e\}} \times [0,1]}$, and
\item the composition
\begin{align*}
&D^i_{e_{[\fraks]}} \times (0,1) \xrightarrow{(\hat{\varphi}_{e_{[\fraks]}} \times \id)^\ast \tilde{\sigma}'^{(i)}} (\hat{\varphi}_{e_{[\fraks]}} \times \id)^\ast \hat{p}^\ast \circPi(E,k)'\\ &\xrightarrow{\hat{\tilde{\psi}}_{e_{[\fraks]}}} D^{i} \times [0,1] \times \circPi(X) \xrightarrow{p_2} \circPi(X)
\end{align*}
is generic for every $[\fraks] \in \Spincns(X,k)/\Conj$, where $p_2$ is the second projection.
\end{itemize}
We may take the last composition map to be a map to $\circPi(X)$ rather than to $\Pi(X)$ since we have $b^+(X) \geq k+2$.
Thus we obtain a section
\[
\tilde{\sigma}'^{(i)} : (\Spincns(E,k)/\Conj)^{(i)} \times [0,1] \to \hat{p}^\ast \circPi(E,k)'.
\]

Equip $[0,1]$ with the standard CW structure and let $I$ denote the unique $1$-cell of $[0,1]$.
We define a cochain $\widetilde{\SWcal} \in C^{k-1}(B)$ by 
\[
\widetilde{\SWcal} = \SWcalhalftot^k(p^\ast E, \tilde{\sigma}')(e \times I)
\]
for a $(k-1)$-cell $e$.

We claim that $\delta \widetilde{\SWcal} = \SWcal_1 + \SWcal_0$, which completes the proof.
First, we note the isomorphisms
\[
C^\ast(B \times [0,1], B \times \{0,1\})
\cong C^\ast(B) \otimes C^1([0,1])
\cong C^\ast(B).
\]
Here the first isomorphism is the K{\"u}nneth isomorphism and the second one is the map $C^\ast(B) \to C^\ast(B) \otimes C^1([0,1])$ given by $u \mapsto u \otimes I$.
Let $e$ be a $k$-cell of $B$.
Using that $\SWcalhalftot^k$ is a cocycle (\cref{SWtot cocycle}) and the naturality (\cref{prop: cochain natural}), we have
\begin{align*}
0=&\delta \SWcalhalftot^k(p^\ast E, \tilde{\sigma}')(e \otimes I)\\
=&\SWcalhalftot^k(p^\ast E, \tilde{\sigma}')(\del (e \otimes I))\\
=& \widetilde{\SWcal}(\del e)
+ \SWcalhalftot^k(p^\ast E, \tilde{\sigma}')(e \otimes 1) + \SWcalhalftot^k(p^\ast E, \tilde{\sigma}')(e \otimes 1)\\
=& \delta\widetilde{\SWcal}(e)
+ \SWcal_1(e) + \SWcal_0(e),
\end{align*}
which verifies the claim.
\end{proof}

\subsection{Vanishing under stabilizations}
\label{subsectionVanishing under stabilizations}

In this \lcnamecref{subsectionVanishing under stabilizations}, 
we give vanishing theorems for $\SWbbtot$ and $\SWbbhalftot$ under stabilizations.
Let $M, N$ closed oriented smooth 4-manifolds.
Choose smoothly embedded open 4-disks in $M, N$, and let $\mathring{M}, \mathring{N}$ denote the manifolds obtained from $M, N$ by removing the open disks.
We denote by $\Diff_\del(\mathring{M})$ the group of diffeomorphisms that are the identities on collar neighborhoods of $\del\mathring{M}$. 
We consider the stabilization map  
\[
s_{M,M\# N} : \Diff_\del(\mathring{M}) \to \Diff^+(M\# N)
\]
defined by extending elements in $\Diff_\del(\mathring{M})$ with the identity of $N$.

\begin{thm}
\label{thm: vanishing}
Let $k\geq0$ and let $M,N$ be closed oriented smooth 4-manifolds with  $b^+(M) \geq k+1$ and $b^{+}(N)\geq 1$.
Then we have
\[
s_{M,M\# N}^\ast\SWbbtot^k(M\#N)=0.
\]
\end{thm}

\begin{proof}
Set $X=M\#N$, $B=\BDiff_\del(M)$ and $E(\mathring{M}) = \EDiff_\del(\mathring{M})$.
By gluing $E(\mathring{M})$ with the trivial bundle $B \times \mathring{N}$ along $B \times S^3$, we obtain a bundle with fiber $X$, which we denote by $E(X) \to B$.
Let $X_b$ denote the fiber of $E(X)$ over $b \in B$.
We use similar notations also for other bundles.

By the construction of $\SWbbtot$, the claim of the \lcnamecref{thm: vanishing} follows once we construct an inductive families perturbation $\sigma : B^{(k)} \to \circPi(E(X))|_{B^{(k)}}$ with the following property:
for every $b \in B^{(k)}$ and
every $\spinc$ structure $\fraks_{{X}_b} \in \Spinc(X_b)$ with $d(\fraks_{{X}_b}) = -k$, the Seiberg--Witten equations for $(X_b,\fraks_{X_b})$ perturbed by $\sigma(b)$ do not have a solution.


To construct such $\sigma$, we shall consider a fiberwise neck stretching argument.
Let us describe a setup to discuss it.
Let $E(\mathring{M}^\ast) \to B$ denote the bundle with fiber $\mathring{M}^\ast = \mathring{M} \cup ([0,\infty)\times S^3)$ obtained by gluing $E(\mathring{M})$ with $B \times ([0,\infty)\times S^3)$.
We fix the standard positive scalar metric on the cylinder $[0,\infty)\times S^3$.
By requiring the restriction to the cylinder is given by this standard metric and the zero perturbation 2-form, one has straightforward variants, such as $\circPi(\mathring{M}^\ast)$, of spaces of perturbations introduced in \cref{subsection Collection of Seiberg--Witten equations} for the cylindrical-end family $\mathring{M}^\ast$.
Whenever we consider the Seiberg--Witten equations on a cylindrical-end 4-manifold, say $\mathring{M}^\ast$, we impose the standard boundary condition: asymptotic to the unique reducible solution on $S^3$.

Similarly, we may also consider the space of perturbations, say $\circPi(\mathring{M})$, for a punctured 4-manifold by imposing that the restriction to neighborhoods of $\del \mathring{M}$ is given by the standard positive scalar curvature metric on $[0,\epsilon) \times S^3$ with the zero perturbation 2-form for some $\epsilon>0$.
Whenever we consider the Seiberg--Witten equations on a punctured 4-manifold, we impose the Atiyah--Patodi--Singer boundary condition.

By the assumption that $b^+(M)\geq k+1$, as in \cref{subsection Inductive sections}, we can construct an inductive families perturbation $\sigma_{\mathring{M}} : B^{(k)} \to \circPi(\mathring{M})$, with an underlying families metric, 
that satisfies the following property:
for every $b \in B^{(k)}$ and
every $\spinc$ structure $\fraks_{{M}_b} \in \Spinc(M_b)$ with $d(\fraks_{{M}_b}) < -k$, the Seiberg--Witten equations for $(\mathring{M}_b,\fraks_{M_b})$ perturbed by $\sigma_{\mathring{M}}(b)$ have neither reducible nor irreducible solution.

On the other hand, by $b^+(N)\geq 1$, we can find a perturbation  $\sigma_{\mathring{N}}^0 \in \circPi(\mathring{N})$, with an underlying metric,
that satisfies the following property:
the Seiberg--Witten equations on $\mathring{N}$ perturbed by $\sigma_{\mathring{N}}^0$ do not have reducible solution for every $\spinc$ structure $\fraks_{N} \in \Spinc(N)$, and have neither reducible nor irreducible solution for every $\spinc$ structure $\fraks_{N} \in \Spinc(N)$ with $d(\fraks_N)<0$.
Define $\sigma_{\mathring{N}} : B^{(k)} \to \circPi(B \times \mathring{N})$ to be the constant families metric/perturbation for the trivial family $B \times \mathring{N}$ given by $\sigma_{\mathring{N}}^0$. 

We equip $\mathring{M}_b$ for $b \in B$ and $\mathring{N}$ with the underlying metrics of $\sigma_{\mathring{M}}(b)$ and $\sigma_{\mathring{N}}^0$.
For $L>0$, define a Riemannian manifold $M_b\#_{L}N$ by inserting a cylinder of length $2L$ between $\mathring{M}_b$ and $\mathring{N}$,
\[
M_b\#_{L}N = 
\mathring{M}_b \cup ([-L,L] \times S^3) \cup \mathring{N},
\]
where the cylinder $[-L,L] \times S^3$ is equipped with the standard positive scalar curvature metric.
By zero extension, we obtain from $\sigma_{\mathring{M}}(b)$ and $\sigma_{\mathring{N}}^0$ a new perturbation $\sigma_{\mathring{M}}(b)\#_L\sigma_{\mathring{N}}^0 \in \circPi(E(X)_b)$ whose underlying metric is given by the metric on $M_b\#_L N$ described above.

Now we define a function $f : B^{(k)} \to [0,\infty) \cup \{\infty\}$ by
\[
f(b) = \sup\Set{L \in [0,\infty) | \text{Property }(*)_{b,L}},
\]
where Property~$(*)_{b,L}$ is that $M_b\#_L N$ admits a solution to the Seiberg--Witten equations with respect to $\sigma_{\mathring{M}}(b) \#_L \sigma_{\mathring{N}}^0$ for some $\spinc$ structure $\fraks$ on $M_b\#N$ with $d(\fraks) = -k$.
(Here we adopt the convention that $\sup\emptyset=0$.)

Now we claim that $f(b)<\infty$ for any $b \in B^{(k)}$, so $f$ is regarded to be a function $f : B^{(k)} \to [0,\infty)$, and further $f$ is upper semi-continuous.
Once we verify this, it is straightforward to inductively construct a continuous function $\mathcal{L} : B^{(k)} \to [0,\infty)$ with $\mathcal{L} > f$ over $B^{(k)}$.
Then the perturbations/metrics $\sigma_{\mathring{M}}(b)\#_{\mathcal{L}(b)}\sigma_{\mathring{N}}^0$ on  $M_b\#_{\mathcal{L}(b)}N$ yield the desired families perturbation $\sigma : B^{(k)} \to \circPi(E(X))$, which completes the proof of the \lcnamecref{thm: vanishing}.

To verify that $f$ is finite-valued and is upper semicontinuous, it is enough to get a contradiction from the following assumption: there are sequences $\{b_i\}_{i=1}^\infty \subset B^{(k)}$ and $\{L_i\}_{i=1}^\infty \subset \R$ such that:
\begin{itemize}
\item $\{b_i\}$ converges to some point $b_{\infty} \in B^{(k)}$, and $L_i \to +\infty$ as $i \to \infty$,
\item for each $i$, there exists a solution $\gamma_i$ to the Seiberg--Witten equations for 
\[
(M_{b_i}\#_{L_i}N, \sigma_{\mathring{M}}(b_i)\#_{L_i}\sigma_{\mathring{N}}^0)
\]
and for some $\spinc$ structure $\fraks_i$ on $M_{b_i}\#N$ with $d(\fraks_i) = -k$.
\end{itemize}

To get a contradiction, first note that we have $\#\Set{\fraks_i|i \in \mathbb{N}}<\infty$, as in the proof of the finiteness of the basic classes.
Thus, after passing to a subsequence, we may suppose that all $\fraks_i$ coincide with a common $\spinc$ structure, denoted by $\fraks$.
Then the solutions $\{\gamma_i\}$ converge to a pair of solutions $\gamma_\infty = (\gamma_\infty^M, \gamma_\infty^N)$ on $(\mathring{M}_{b_\infty}^\ast,  \mathring{N}^\ast)$ for the pair of perturbations $(\sigma_{\mathring{M}}(b_\infty), \sigma_{\mathring{N}}^0)$, extended to $\mathring{M}^\ast$ and $\mathring{N}^\ast$ by zero, and for the pair of $\spinc$ structures $(\fraks|_{M}, \fraks|_{N})$.
The existence of the solution $\gamma_{\infty}^N$ and the choice of $\sigma_{\mathring{N}}^0$ implies that $d(\fraks|_N) \geq 0$.
Now $d(\fraks|_{M})<-k$ follows from this combined with $d(\fraks) = -k$ and a general formula $d(\fraks|_{M})+d(\fraks|_{N})+1=d(\fraks)$.
However, then the existence of the solution $\gamma_\infty^M$ contradicts the choice of $\sigma_{\mathring{M}}$.
Thus we completed to prove the existence of the desired families perturbation $\sigma$ and this completes the proof of the \lcnamecref{thm: vanishing}.
\end{proof}

For $\SWbbhalftot$, 
a statement analogous to \cref{thm: vanishing} does not hold in general.
The upshot is that, for $\SWbbhalftot$, one may not take a constant families perturbation for $N$ as in the proof of \cref{thm: vanishing}.

\begin{ex}
To get a counterexample to a  statement for $\SWbbhalftot$ analogous to \cref{thm: vanishing}, let us take $N$ to be $S^2\times S^2$ and $M$ to be a complete intersection of general type defined by real coefficient polynomials so that $M$ admits an involution $f : M \to M$ given by the complex conjugation.
Then we have
\begin{align}
\label{eq: counter ex}
s_{M,M\#N}^\ast(\SWbbhalftot^1(M\#N)) \neq 0.
\end{align}
As (\ref{eq: counter ex}) is not needed in the proofs of our results, we only sketch the argument and summerize the upshot.
By isotopy, we may assume that $f$ fixes a 4-disk in $M$.
Let $M_f \to S^1$ denote the mapping torus with fiber $M$ defined by $f$.
Let $BM_f : S^1 \to \BDiff_\del(M)$ be the classifying map of $M_f$ and define $\alpha := (BM_f)_\ast([S^1]) \in H_1(\BDiff^+(M);\Z/2)$.
The non-vanishing \eqref{eq: counter ex} follows from the fact that 
\begin{align}
\label{eq: counter ex pairing}
\left<\SWbbhalftot^1(M\#N),(s_{M,M\#N})_\ast(\alpha)\right> \neq 0,
\end{align}
which is derived by an argument similar \cref{thm key computation source} which we shall prove.
The point is, since $f$ sends the canonical class of $M$ to its conjugate, one needs to take a families perturbation for $N$ that causes the wall-crossing in the calculation of the pairing \eqref{eq: counter ex pairing}.
\end{ex}

However, we still have the following vanishing for $\SWbbhalftot$:

\begin{thm}
\label{thm: vanishing2}
Let $k\geq0$ and let $M,N$ be closed oriented smooth 4-manifolds with  $b^+(M),b^{+}(N)\geq k+1$.
Then we have
\[
s_{M,M\# N}^\ast\SWbbhalftot^k(M\#N)=0.
\]
\end{thm}

\begin{proof}
The proof is quite similar to that of \cref{thm: vanishing}, so we describe only the upshot.
We use notations in the proof of \cref{thm: vanishing}.
Also, set $E(N) = B \times N$ and use notations $E(\mathring{N})$, $E(\mathring{N}^\ast)$ in a similar manner to $M$.

First, note that the $\Z/2$-actions induced by the charge conjugation make
a natural bijection 
\begin{align}
\label{eq: identify}
\Spinc(X,k) \cong \bigtimes_{k=p+q+1} (\Spinc(M,p) \times \Spinc(N,q))
\end{align} 
$\Z/2$-equivariant, where the  $\Z/2$-action on the right-hand side is the diagonal action.
With this in mind,
we may modify the definition of (spaces of)  half-total families perturbations in \cref{subsection Complex conjugation,subsection Inductive sections} to a pair of (punctured) 4-manifolds.
We define
\[
\circPi(\mathring{M}, \mathring{N},k)' = \left(\bigtimes_{k=p+q+1}\left((\Spinc(M,p) \times \circPi(\mathring{M})) \times (\Spinc(N,q) \times \circPi(\mathring{N}))\right)\right)/(\Z/2),
\]
where $\Z/2$ acts as the diagonal action.
Then, using the $\Z/2$-equivariant bijection \eqref{eq: identify}, $E(M)$ and $E(N)$ induce a fiber bundle
\[
\circPi(E(\mathring{M}), E(\mathring{N}),k)' \to \Spinc(E(X),k)/\Conj
\]
with fiber $\circPi(\mathring{M}, \mathring{N},k)'$.
A half-total families perturbation is formulated as a section 
\[
\sigma'_{M,N} : \Spinc(E(X),k)/\Conj \to \circPi(E(\mathring{M}), E(\mathring{N}),k)'
\]
that covers, as in the diagram \eqref{eq: diagram half-tot section}, a pair of families metrics $(\tilde{g}_M, \tilde{g}_N)$ for $E(\mathring{M})$, $E(\mathring{N})$ with cylindrical metrics near boundary.
We require that the self-dual 2-forms that appear in a half-total families perturbation is zero near the boundaries of $\mathring{M}_b, \mathring{N}_b$.
Then we may glue $\sigma_{M,N}'$ near $\del \mathring{M}_b, \del \mathring{N}_b$ to get a half-total families perturbation for $E(X)$ in the sense of \cref{subsection Inductive sections}. 
We denote by 
\[
\sigma_{\mathring{M}}'\#_L\sigma_{\mathring{N}}' : \Spinc(E(X),k)/\Conj \to \circPi(E(X),k)'
\]
the glued half-total families perturbation obtained by inserting the cylinder $[-L,L] \times S^3$.
Note also that $\sigma_{M,N}'$ induces half-total families perturbations
\begin{align*}
&\sigma_{\mathring{M}}' : \Spinc(E(M),p)/\Conj \to \circPi(E(\mathring{M}),p)',\\
&\sigma_{\mathring{N}}' : \Spinc(E(N),q)/\Conj \to \circPi(E(\mathring{N}),q)'
\end{align*}
for $E(\mathring{M})$ and $E(\mathring{N})$ in the sense of \cref{subsection Inductive sections}, only with the modification that we allow punctured 4-manifolds not only closed 4-manifolds.

Let $\calM_{\sigma'_{\mathring{M}}, \mathrm{half},k}$ and $\calM_{\sigma'_{\mathring{N}}, \mathrm{half},k}$ be the paramaterized moduli spaces  over $B^{(k)}$ for $\sigma'_{\mathring{M}}$ and $\sigma'_{\mathring{N}}$  defined by considering collections over $\bigsqcup_{l > k}\Spinc(M,l)/\Conj$ and $\bigsqcup_{l > k}\Spinc(N,l)/\Conj$, respectively.
By $b^+(M), b^+(N) \geq k+1$,
one can take an inductive half-total families perturbation
$\sigma'_{M,N}$ so that
$\calM_{\sigma'_{\mathring{M}}, \mathrm{half},k}$ and $\calM_{\sigma'_{\mathring{N}}, \mathrm{half},k}$ are empty.

Let $\pi : \Spinc(E(X),k)/\Conj \to B$ be the projection.
Now we define a function $f' : B^{(k)} \to [0,\infty) \cup \{\infty\}$ by
\[
f'(b) = \sup\Set{L \in [0,\infty) | \text{Property }(*)'_{b,L}},
\]
where Property~$(*)'_{b,L}$ is that $M_b\#_L N_b$ admits a solution to the Seiberg--Witten equations with respect to some perturbation that lies in $(\sigma_{\mathring{M}}' \#_L \sigma_{\mathring{N}}')(\pi^{-1}(b))$ and for some $\spinc$ structure $\fraks$ on $M_b\#N_b$ with $d(\fraks) = -k$.

Just as in the proof of \cref{thm: vanishing}, we have that $f'<+\infty$ over $B^{(k)}$ and $f'$ is upper semicontinuous, and there is a continuous function $\mathcal{L}' : B^{(k)} \to [0,\infty)$ with $\mathcal{L}'>f'$ over $B^{(k)}$, and $\sigma_{\mathring{M}}'\#_{\mathcal{L}'(b)} \sigma_{\mathring{N}}' : \Spinc(E(X),k)/\Conj \to \circPi(E(X),k)'$ is an inductive half-total families perturbation for which all moduli spaces relevant to $s_{M,M\# N}^\ast\SWbbhalftot^k(M\#N)$ are empty, which completes the proof.
\end{proof}

Now, for an oriented closed smooth 4-manifold $X$, let 
\[
s_X : \Diff_\del(\mathring{X})
\to \Diff^+(X\# S^2\times S^2)
\]
denote the stabilization map for $X$ by $S^2 \times S^2$.
\cref{thm: vanishing,thm: vanishing2} immediately imply the following, which states that $\SWbbtot$ and $\SWbbhalftot$ are {\it unstable} characteristic classes with respect to stabilizations by $S^2\times S^2$:

\begin{cor}
\label{cor: vanishing}
Let $k\geq0$ and let $X$ be a closed oriented smooth 4-manifolds with $b^+(X) \geq k+1$. 
Then we have
\begin{align*}
&s_{X}^\ast\SWbbtot^k(X\#S^2\times S^2)=0,\\
&s_{X}^\ast \circ \dots \circ s_{X\#kS^2\times S^2}^\ast \SWbbhalftot^k(X\#(k+1)S^2\times S^2)=0.
\end{align*}
\end{cor}

Note that, in the above proof of \cref{thm: vanishing2}, the triviality of the bundle $E(N)=B\times N$ is not used, and a similar (actually rather simpler) argument applies also to $\SWbbtot$.
We record this below, while we do not use it in the rest of this paper.
Let $M, N$ be closed oriented smooth 4-manifolds, and $B$ be a CW complex.
Let $E(M) \to B$ and $E(N) \to B$ be oriented smooth fiber bundles with fiber $M$ and $N$ with structure group $\Diff_\del(\mathring{M})$ and $\Diff_\del(\mathring{N})$, respectively.
Then we can form an oriented fiber bundle $E(M)\#_fE(N) \to B$ with fiber $M\#N$ obtained by taking the fiberwise connected sum of $E(M)$ and $E(N)$ along the trivial families of 4-disks in $E(M), E(N)$ corresponding to the punctures of $M, N$. 

\begin{thm}
In the above setup, let $k\geq0$ and suppose that $b^+(M), b^+(N)\geq k+1$. Then we have
\[
\SWbbtot^k(E(M)\#_fE(N))=0,\quad
\SWbbhalftot^k(E(M)\#_fE(N))=0.
\]
\end{thm}

\subsection{Positive scalar curvature}
\label{subsection Vanishing under positive scalar curvature}

The vanishing of the Seiberg--Witten invariant for a positive scalar curvature metric extends to invariants for families, and this generalization is useful to study the topology of spaces of positive scalar curvature metrics  \cite{Rub01,K19}.
In this \lcnamecref{subsection Vanishing under positive scalar curvature}, we formulate a corresponding vanishing theorem in our framework.

Given an oriented smooth manifold $X$,  recall that we defined notations
\[
\scrR(X),\ \scrR^+(X),\ \scrM(X),\  \scrM^+(X)
\]
in \cref{subsection:Moduli spaces of manifolds with metrics}.
We start by noting that the space $\scrM(X)$ classifies pairs $(E,g_E)$ of oriented fiber bundles $E$ with fiber $X$ and fiberwise metrics $g_E$ on $E$ up to deformation equivalence, and $\scrM^+(X)$ classifies pairs $(E,g_E)$ of oriented fiber bundles $E$ with fiber $X$ and fiberwise positive scalar curvature metrics $g_E$ up to deformation equivalence. 
Indeed, let $B$ be a CW complex and let $X \to E \to B$ be an oriented fiber bundle with fiber $X$ over $B$.
Let $\scrR(X) \to \scrR(E) \to B$ denote the fiber bundle associated with $E$ with fiber $\scrR(X)$.
Suppose that $E$ is equipped with a fiberwise metric $g_E$.
This corresponds to a section of $\scrR(E) \to B$, which gives a lift of the classifying map $B \to \BDiff^+(X)$ of $E$ to $\EDiff^+(X) \times_{\Diff^+(X)} \scrR(X)$.
Similarly, a fiberwise positive scalar curvature metric $g_E$ on $E$ corresponds to a lift of the classifying map $B \to \BDiff^+(X)$ of $E$ to $\EDiff^+(X) \times_{\Diff^+(X)} \scrR^+(X)$:
\begin{align}
\label{eq: moduli of psc}
\begin{split}
\xymatrix{
     & \EDiff^+(X)\times_{\Diff^+(X)} \scrR^+(X) = \scrM^+(X) \ar[d] \\
    B \ar[ru]^-{g_E} \ar[r]_-{E} & \BDiff^+(X).
   }   
\end{split}  
\end{align}

There are also other kinds of natural moduli spaces.
Fix a point $x_0 \in X$, and let $\Diff_{x_0}(X)$ denote the group of diffeomorphisms that fix $x_0$ and act on $T_{x_0}X$ trivially.
If we suppose that $X$ is a closed manifold, then $\Diff_{x_0}(X)$ acts freely on $\scrR(X)$ (see, e.g. \cite[Lemma~1.2]{BHSW10}).
The quotients 
\[
\scrM_{x_0}(X) = \scrR(X)/\Diff_{x_0}(X),\quad
\scrM_{x_0}^+(X) = \scrR^+(X)/\Diff_{x_0}(X)
\]
are called the {\it observer moduli spaces} of Riemannian metrics and metrics of positive scalar curvature, respectively.
These moduli spaces were introduced in \cite{AB02}, and have been extensively studied by several authors for higher-dimensional manifolds.
See such as \cite{BHSW10,HSS14}.

Relative versions are also natural moduli spaces.
Choose an embedded 4-disk $D^4$ in $X$, and set $\mathring{X} = X \setminus \mathrm{Int}(D^4)$.
Let $\scrR_{\del}(\mathring{X})$ denote the space of Riemannian metrics on $X$ for which the restrictions to collar neighborhoods of $\del\mathring{X}$ coincide with the standard cylindrical metric of positive scalar curvature on $[0,1] \times S^3$.
Set $\scrR_{\del}^+(\mathring{X}) = \scrR_{\del}(\mathring{X}) \cap \scrR^+(\mathring{X})$.
It is easy to see that $\Diff_{\del}(\mathring{X})$ acts freely on $\scrR_{\del}^+(\mathring{X})$.
Define
\[
\scrM_{\del}(\mathring{X}) = \scrR_\del(\mathring{X})/\Diff_{\del}(\mathring{X}),\quad
\scrM_{\del}^+(\mathring{X}) = \scrR^+_{\del}(\mathring{X})/\Diff_{\del}(\mathring{X}).
\]


Let us compare the various moduli spaces above.
Let
\[
\iota : \scrM^+(X) \hookrightarrow \scrM(X),\quad
\iota_{x_0} : \scrM_{x_0}^+(X) \hookrightarrow \scrM_{x_0}(X),\quad
\iota_{\del} : \scrM_{\del}^+(\mathring{X}) \hookrightarrow \scrM_{\del}(\mathring{X})
\]
be injections induced from the inclusion $\scrR^+(X) \hookrightarrow \scrR(X)$.
Since the actions of $\Diff_{x_0}(X)$ and $\Diff_{\del}(\mathring{X})$ on spaces of metrics are free, the honest quotients by these groups are homotopy equivalent to corresponding homotopy quotients.
Also, let $D^4$ be equipped with the restriction of the standard round metric of $S^4$.
Gluing relative metrics with this round metric on $D^4$,
we have natural extension maps $\scrR_{\del}(\mathring{X}) \hookrightarrow \scrR(X)$ and $\scrR_{\del}^+(\mathring{X}) \hookrightarrow \scrR^+(X)$.
With homotopy inverses from honest quotients to homotopy quotients, the injections $\Diff_{\del}(\mathring{X}) \hookrightarrow \Diff_{x_0}(X) \hookrightarrow \Diff^+(X)$ and extensions $\scrR_{\del}(\mathring{X}) \hookrightarrow \scrR(X)$,  $\scrR_{\del}^+(\mathring{X}) \hookrightarrow \scrR^+(X)$ give rise to a homotopy commutative diagram
\begin{align}
\label{eq: diagram three psc moduli comparison}
\begin{split}
\xymatrix{
\scrM_\del(\mathring{X})
\ar[r]& 
\scrM_{x_0}(X)
\ar[r] &
\scrM(X)\\
\scrM_\del^+(\mathring{X})
\ar[r] \ar[u]^{\iota_{\del}}& 
\scrM_{x_0}^+(X)
\ar[r] \ar[u]^{\iota_{x_0}}&
\scrM^+(X). \ar[u]^{\iota}
}
\end{split}
\end{align}

Now we discuss the vanishing of the Seiberg--Witten characteristic classes.
Let $k>0$ and suppose that $X$ is an oriented closed smooth 4-manifold with $b^+(X) \geq k+2$ so that the Seiberg--Witten characteristic classes are well-defined.
Since $\scrR(X)$ is contractible, $\EDiff^+(X) \times \scrR(X)$ is also contractible, and $\EDiff^+(X) \times \scrR(X)$ is acted by $\Diff^+(X)$ freely.
Thus $\scrM(X)$ gives a model of the classifying space of $\Diff^+(X)$, and hence we may regard the Seiberg--Witten characteristic classes as cohomology classes on $\scrM(X)$:
\[
\SWbbtot^k(X) \in H^k(\scrM(X);\Z_{\EDiff^+(X)}),\quad
\SWbbhalftot^k(X) \in H^k(\scrM(X);\Z/2).
\]

\begin{thm}
\label{thm: vanishing by psc}
We have
\[
\iota^\ast\SWbbtot^k(X)=0 \text{\quad and\quad }
\iota^\ast\SWbbhalftot^k(X)=0
\]
in $H^k(\scrM^+(X);\iota^\ast\Z_{\EDiff^+(X)})$ and in $H^k(\scrM^+(X);\Z/2)$ respectively.
\end{thm}

\begin{proof}
Fix a model of the universal principal $\Diff^+(X)$-bundle $\EDiff^+(X) \to \BDiff^+(X)$.
Then the bundle $\EDiff^+(X) \times \scrR(X) \to \scrM(X)$  gives another model of the universal pricipal $\Diff^+(X)$-bundle.
Let $E \to \scrM^+(X)$ denote the pull-back of 
this bundle $\EDiff^+(X) \times \scrR(X) \to \scrM(X)$ under $\iota$.
The claim of the \lcnamecref{thm: vanishing by psc} is equivalent to that $\SWbbtot^k(E)=0$ and $\SWbbhalftot^k(E)=0$.

We first prove $\SWbbtot^k(E)=0$.
Let $E_b$ denote the fiber of $E$ over a point $b \in \scrM^+(X)$.
Recall the notation
$\Pi_\ast(X) = \bigcup_{g \in \scrR(X)} L^2_{l-1}(i\Lambda^+_g(X))$ and
regard $\scrR(X)$ as a subset of $\Pi_\ast(X)$.
From the construction of $E$, we have $E = \EDiff^+(X) \times \scrR^+(X)$.
It follows from this that a $\Diff^+(X)$-equivariant map $E \to \Pi_\ast(X)$ gives rise to a section of the fiber bundle
\[
\Pi_\ast(X) \to \Pi_\ast(E) \xrightarrow{\pi} \scrM^+(X),
\]
which is a fiber bundle associated to $E$ with fiber $\Pi_\ast(X)$ defined in  \cref{subsection Collection of Seiberg--Witten equations}.
Similarly, we have an associated fiber bundle $\scrR^+(E) \to \scrM^+(X)$ with fiber $\scrR^+(X)$.
Note that we have a fiberwise inclusion $\scrR^+(E) \hookrightarrow \Pi_\ast(E)$ induced from the inclusion $\scrR^+(X) \hookrightarrow \Pi_\ast(X)$.

Let $\phi : E \to \Pi_\ast(X)$ be a $\Diff^+(X)$-equivariant map defined as the composition of the projection $\EDiff^+(X) \times \scrR^+(X) \to \scrR^+(X)$ and the inclusion $\scrR^+(X) \hookrightarrow \Pi_\ast(X)$.
Let 
\[
\sigma: \scrM^+(X) \to \Pi_\ast(E)
\]
denote the section induced from $\phi$.
Then $\sigma$ factors through $\scrR^+(E)$, and hence, for each $b \in \scrM^+(X)$, the Seiberg--Witten equations for $(E_b, \sigma(b))$ have no irreducible solutions for all $\spinc$ structures because of the standard vanishing for positive scalar curvature metrics.

Now we slightly perturb $\sigma$ to avoid reducibles.
First, the non-existence of irreducible solutions is an open condition, hence there is an open neighborhood $\mathcal{N} \subset \Pi_\ast(E)$ of the image of $\sigma$ for which, for every $n \in \mathcal{N}$, the Seiberg--Witten equations for $n$ have no irreducible solutions for all $\spinc$ structures.
We may suppose that $\mathcal{N}$ lies in $\Pi(E)$, the space of perturbations with uniformly bounded norms.
Since $\mathcal{N} \cap \circPi(E_b)$ is a codimension-$b^+(X)$ subspace of $\mathcal{N}_b = \Set{n \in \mathcal{N} | \pi(n)=b}$,
$\mathcal{N} \cap \circPi(E_b)$ is $(b^+(X)-2)$-connected, as well as $\circPi(X)$.
Thus, as in \cref{subsection Inductive sections},
we may construct an inductive families perturbations $\tilde{\sigma} : \scrM^+(X)^{(k)} \to \circPi(E)$ whose image lies in $\mathcal{N}$.
This implies that $\SWbbtot^k(E)=0$.

The proof of that $\SWbbhalftot^k(E)=0$ is similar to the above.
Let us take the subset $\mathcal{N} \subset \Pi_\ast(E)$ as above.
Define $\Pi_\ast(X,k)' = (\Spinc(X,k) \times \Pi_\ast(X))/(\Z/2)$.
Recall that we let $\Z/2$ act on $\Pi_\ast(E)$ as the fiberwise $(-1)$-multiplication, where we regard $\Pi_\ast(E)$ as a fiber bundle with fiber $L^2_{l-1}(i\Lambda^+_g(X))$ for a metric $g$.
By taking $\mathcal{N}$ smaller, we may assume that $\mathcal{N}$ is preserved under the $\Z/2$-action.
This $\mathcal{N}$ induces an open neighborhood $\mathcal{N}' \subset \Pi_\ast(E,k)'$ of $(\Spinc(X,k) \times \scrR^+(X))/(\Z/2)$ for which, for every $n \in \mathcal{N}'$, the Seiberg--Witten equations for $n$ have no irreducible solutions.
Then we can construct an inductive half-total families perturbation
$\sigma' : \left(\Spincns(E,k)/\Conj\right)^{(i)}
\to \circPi(E,k)'$
whose image lies in $\mathcal{N}'$.
This proves that $\SWbbhalftot^k(E)=0$.
\end{proof}

\begin{cor}
\label{cor: vanishing for fiberwise psc}
Let $X \to E \to B$ be an oriented smooth fiber bundle with fiber $X$ that admits a fiberwise positive scalar curvature metric.
Then we have 
\[
\SWbbtot^k(E)=0,\quad
\SWbbhalftot^k(E)=0.
\]
\end{cor}

\begin{proof}
This is a direct consequence of \cref{thm: vanishing by psc} combined with the commutative diagram \eqref{eq: moduli of psc}.
\end{proof}

\begin{cor}
\label{cor: psc non-surj}
Suppose that there is a homology class $\alpha \in H_k(\BDiff^+(X);\Z)$ such that either $\langle\SWbbtot^k(X),\alpha\rangle \neq0$ or $\langle\SWbbhalftot^k(X),\alpha\rangle\neq0$ holds.
Then we have the following:
\begin{enumerate}[label=(\roman*)]
\item The map
\[
\iota_\ast : H_k(\scrM^+(X);\Z) \to H_k(\scrM(X);\Z)
\]
is not surjective.
\item If $\alpha$ lies in the image of the map $H_k(\BDiff_{x_0}(X);\Z) \to H_k(\BDiff^+(X);\Z)$, then the map
\[
(\iota_{x_0})_\ast :H_k(\scrM^+_{x_0}(X);\Z) \to H_k(\scrM_{x_0}(X);\Z)
\]
is not surjective.
\item If $\alpha$ lies in the image of the map $H_k(\BDiff_{\del}(\mathring{X});\Z) \to H_k(\BDiff^+(X);\Z)$, then the map
\[
(\iota_{\del})_\ast:
H_k(\scrM^+_{\del}(\mathring{X});\Z) \to H_k(\scrM_{\del}(\mathring{X});\Z)
\]
is not surjective.
\end{enumerate}
\end{cor}

\begin{proof}
It follow from \cref{thm: vanishing by psc} that $\alpha$ lies in  the cokerel of $\iota_\ast$, which proves the claim (i).
The remaining claims (ii), (iii) follow from the homotopy commutative diagram \eqref{eq: diagram three psc moduli comparison} with homotopy equivalences $\scrM_{x_0}(X) \simeq \BDiff_{x_0}(X)$ and $\scrM_{\del}(\mathring{X}) \simeq \BDiff_{\del}(\mathring{X})$.
\end{proof}

We remark that the non-surjectivity result (\cref{cor: psc non-surj}) for $k=1$ implies the disccoectivity of $\scrR^+(X)$:

\begin{lem}
\label{lem psc disconnected}
Let $X$ be a smooth manifold.
If $\iota_\ast : H_1(\scrM^+(X);\Z) \to H_1(\scrM(X);\Z)$ is not surjective, then $\scrR^+(X)$ is not connected.
\end{lem}

\begin{proof}
Let $\gamma \in H_1(\scrM(X);\Z)$ be a non-zero element in the cokernel of $\iota_\ast$.
Since the Hurewicz map $h : \pi_1(\scrM(X)) \to H_1(\scrM(X);\Z)$ is surjective, we may take an element $\hat{\gamma} \in h^{-1}(\gamma)$.
Since the Hurewicz maps are compatible with $\iota_\ast$, we have that $\hat{\gamma} \notin \mathrm{Im}(\pi_1(\scrM^+(X)) \to \pi_1(\scrM(X)))$.

Now consider the following commutative diagram induced from the fibrations $\Diff^+(X) \to \EDiff^+(X) \times \scrR(X) \to \scrM(X)$ and $\Diff^+(X) \to \EDiff^+(X) \times \scrR^+(X) \to \scrM^+(X)$:
\begin{align*}
\begin{split}
\xymatrix{
0 \ar[r]&
\pi_1(\scrM(X))
\ar[r]^-{\del}_-{\cong}& 
\pi_0(\Diff^+(X))
\ar[r] &
0\\
\pi_1(\scrR^+(X)) \ar[r]\ar[u]&
\pi_1(\scrM^+(X))
\ar[r] \ar[u]& 
\pi_0(\Diff^+(X))
\ar[r]^{i_\ast} \ar@{=}[u]&
\pi_0(\scrR^+(X)). \ar[u]
}
\end{split}
\end{align*}
Here the two rows are exact.
By $\hat{\gamma} \notin \mathrm{Im}(\pi_1(\scrM^+(X)) \to \pi_1(\scrM(X)))$, it is easy to check that $i_\ast \circ \del(\hat{\gamma}) \in \pi_0(\scrR^+(X))$ is non-trivial.
\end{proof}

We shall prove the following \lcnamecref{thm: psc main} in \cref{subsectionProof of the main instability theorem} that generalizes \cref{thm: psc main intro} explained in \cref{section Introduction}: 

\begin{thm}
\label{thm: psc main}
Let $X$ be a simply-connected closed oriented smooth 4-manifold and let $k>0$.
Then there exists a sequence of positive integers $N_1<N_2<\cdots \to +\infty$ such that for all $i$, none of the maps
\begin{align*}
\iota_\ast &: H_k(\scrM^+(X\# N_iS^2\times S^2);\Z) \to H_k(\scrM(X\# N_iS^2\times S^2);\Z),\\
(\iota_{x_0})_\ast &: H_k(\scrM_{x_0}^+(X\# N_iS^2\times S^2);\Z) \to H_k(\scrM_{x_0}(X\# N_iS^2\times S^2);\Z),\\
(\iota_{\del})_\ast &: H_k(\scrM_{\del}^+(\mathring{X}\# N_iS^2\times S^2);\Z) \to H_k(\scrM_{\del}(\mathring{X}\# N_iS^2\times S^2);\Z)
\end{align*}
are surjective.
\end{thm}

\begin{rmk}
For the moduli space $\scrM_{\del}(\mathring{X})$, we may consider stabilization.
Fix a metric on $([0,1] \times S^3)\#S^2\times S^2$ whose restriction to a collar neighborhood of the boundary is a pair of the standard cylindrical metrics on $[0,1] \times S^3$.
Then we may define a stabilization map $\scrR(\mathring{X}) \to \scrR(\mathring{X}\#S^2\times S^2)$ by gluing the fixed metric on $([0,1] \times S^3)\#S^2\times S^2$.
This gives rise to a stabilization map between moduli spaces,
$s_{\scrM} : \scrM_{\del}(\mathring{X}) \to \scrM_{\del}(\mathring{X}\#S^2\times S^2)$.
Seeing the proof of \cref{thm: psc main}, one can find the cokernels of the maps $(\iota_{\del})_\ast$ contain non-zero elements $\alpha_i \in H_k(\scrM_{\del}(\mathring{X}\# N_iS^2\times S^2);\Z)$ that are unstable, i.e. 
\[
\alpha_i \in 
\mathrm{coker}((\iota_{\del})_\ast)
\cap
\ker((s_{\scrM})_\ast).
\]
\end{rmk}

\section{Calculation}
\label{sec Calculation}

\subsection{Multiple mapping torus}
\label{subsection Multiple mapping torus pre}

This \lcnamecref{subsection Multiple mapping torus pre} is devoted to giving preliminary results on multiple mapping tori to carry out calculations of the Seiberg--Witten characteristic class in subsequent discussions. 

We start by setting up some notations. Let $W$ be an oriented 4-manifold with (possibly empty) boundary. We use $\Diff_{\del}(W)$ to denote the group of orientation-preserving diffeomorphisms on $W$ that equals identity near the boundary. (So $\Diff_{\del}(W)$ is just $\Diff^{+}(W)$ when $W$ is closed.) We use $\Diff_{\del}(W)^{\delta}$ to denote the same group equipped with the discrete topology. Consider a smooth oriented bundle $W\rightarrow E\rightarrow B$ with structure group $\Diff_{\partial}(W)$. We also use $E$ to denote the corresponding classifying map $E: B\rightarrow \BDiff_{\partial}(W)$. This map is well-defined up to homotopy. Similarly, when $E$ carries a flat structure, we also use $E$ to denote the corresponding map from $B$ to the moduli space $\BDiff_{\partial}(W)^{\delta}$.

Let $f_1, \ldots, f_k$ be mutually commuting orientation-preserving diffeomorphisms of $W$ supported in $\Int(W)$.
Let 
\[
W \to W_{f_1, \ldots, f_k} \to T^k
\]
denote the multiple mapping torus of $f_1, \ldots, f_k$ with fiber $W$. The bundle $W_{f_1, \ldots, f_k}$ carries a canonical flat structure since it is the quotient of the product bundle $W\times \mathbb{R}^{k}$ by a free $\mathbb{Z}^{k}$-action. Hence we get an induced map 
\[W_{f_1, \ldots, f_k}: T^{k}\rightarrow \BDiff_{\partial}(W)^{\delta}.\]
We give the following technical lemma, which shall be useful later.
\begin{lem}\label{lem: 2-torsion} Let $W'\subset \Int(W)$ be a  compact  embedded 4-dimensional submanifold with boundary and let $U$ be an open neighborhood of $W\setminus W'$ in $W$.  
Suppose that, for any $i\geq 2$, the restriction $f_{i}|_{W'}$ is the identity. Moreover, suppose that there is a smooth isotopy $F^{t}_{1}$ from $f_{1}$ to the identity that fixes $U$ pointwise. Then we have 
\[
W_{f_1,\cdots,f_k,*}[T^{k}]=0\in H_{k}(\BDiff_{\del}(W)^{\delta};\mathbb{Z}).
\]
\end{lem}
\begin{proof}
By our assumption, the restriction $f_{1}|_{W'}$ belong to the unit component of $\Diff_{\del}(W')$, denoted by $\Diff^{0}_{\del}(W')$. By Thurston's result \cite{Thurston74}, the group $\Diff^{0}_{\del}(W')$ is simple. So $f_{1}|_{W'}$ equals a product of commutators of some elements in $\Diff^0_{\del}(W')$. By extended these elements with the identity map on $W\setminus W'$, we obtain diffeomorphisms $g_{1},\cdots,g_{n}, h_{1},\cdots, h_{n}$ that satisfy
\[
f_{1}=\prod_{j=1}^{n}g_{j}h_{j}g^{-1}_{j}h^{-1}_{j}.
\]
Furthermore, for any $1\leq j\leq n$ and $2\leq i\leq k$, both $g_{j}$ and $h_{j}$  commute with $f_{i}$ since their supports do not intersect. Let $\Sigma$ be a compact oriented surface of genus $n$ and with a single boundary component. Consider a homomorphism 
\[
\rho: \pi_{1}(\Sigma\times T^{k-1})\rightarrow \Diff_{\partial}(W)
\]
that sends the generators to $g_{1},\cdots,g_{n},h_{1},\cdots,h_{n},f_{2},\cdots, f_{k}$. Then $\rho$ induces a flat bundle \[W\rightarrow E\rightarrow \Sigma\times T^{k-1}\] whose restriction to $\partial (\Sigma\times T^{k-1})$ is isomorphic to $W_{f_1,\cdots,f_k}$. As a result, the map
\[
W_{f_1,\cdots,f_k}: T^{k}\rightarrow \BDiff_{\del}(W)^\delta
\]
can be extended to $\Sigma\times T^{k-1}$. This directly implies $W_{f_1,\cdots,f_k,*}[T^{k}]=0$. 
\end{proof}

Now we start calculating the families Seiberg--Witten invariants of multiple mapping tori. For the rest of this \lcnamecref{subsection Multiple mapping torus pre},
we fix $k>0$ and set $W=X$ to be a closed oriented smooth 4-manifold with $b^+(X) \geq k+2$. We first note that one may interpret the evaluation
\[
\langle\SWbbhalftot^k(X_{f_1, \ldots, f_k}), [T^k] \rangle \in \Z/2
\]
in the following way. 
Let us fix a section \[\tau : \Spincns(X,k)/\Conj \to \Spincns(X,k).\]
As explained in \cref{ex: perturbation multi map torus}, we can get a families perturbation for $X_{f_1, \dots, f_k}$ from a collection of maps
\begin{align*}
\left(\sigma_{[\fraks]} : [0,1]^k \to \circPi(X)
\right)_{[\fraks] \in \Spincns(X,k)/\Conj}
\end{align*}
with the patching condition \eqref{eq: peruturb general mapping tori}.
For each $\fraks' \in \Spinc(X)$ and each $\sigma_{[\fraks]}$, we may consider the parameterized moduli space
\[
\calM(X, \fraks', \sigma_{[\fraks]})
= \bigcup_{\mathbf{t} \in [0,1]^k} \calM(X, \fraks', \sigma_{[\fraks]}(\mathbf{t})),
\]
where $\calM(X, \fraks', \sigma_{[\fraks]}(\mathbf{t}))$ is just the usual Seiberg--Witten moduli space for the perturbation 
$\sigma_{[\fraks]}(\mathbf{t})$.
If we take generic $\sigma_{[\fraks]}$ so that the parameterized moduli spaces for $\sigma_{[\fraks]}$ over $\del [0,1]^k$ are empty for all $\spinc$ structures, this gives rise to an inductive families perturbation $\sigma : \Spincns(E,k)/\Conj \to \circPi(E,k)'$ in  the sense of \cref{subsection Inductive sections}.
Also we take $\sigma_{[\fraks]}$ so that the parameterized moduli spaces for $\sigma_{[\fraks]}$ is smooth over $\Int([0,1]^k)$ for all $\spinc$ structures.
Then we have
\begin{align}
\label{eq: eval sum}
\langle\SWbbhalftot^k(X_{f_1, \ldots, f_k}), [T^k] \rangle  = \sum_{\fraks \in \tau(\Spincns(X,k)/\Conj)}\#\calM(X, \fraks, \sigma_{[\fraks]})
\end{align}
in $\Z/2$, where $\#$ denotes the mod 2 count of the moduli space.

A main goal of the rest of this \lcnamecref{subsection Multiple mapping torus pre} is to reduce the right-hand side of \eqref{eq: eval sum} to a computable quantity under a certain assumption.
As a first step, it is convenient to introduce the numerical families Seiberg--Witten invariant here.
Given a family $X \to E \to B$ with structure group $\Diff^+(X)$ over a closed manifold $B$ of dimension $k$ and $\fraks \in \Spinc(X,k)$, if the monodromy action of $E$ fixes $\fraks$,
one may define the numerical families Seiberg--Witten invariant \cite{BK20,K21}:
\[
\SW(E, \fraks) \in \Z/2.
\]
This is similar to the families Seiberg--Witten invariant defined by Li--Liu~\cite{LiLiu01} for families of $\spinc$ 4-manifolds, but one may define $\SW(E, \fraks)$ for a family of 4-manifolds (rather than $\spinc$ 4-manifolds) as far as the monodromy action fixes the isomorphism class $\fraks$.
If we have a generic families perturbation as explained in \cref{subsection Inductive sections}, we can calculate $\SW(E, \fraks)$ by counting the parameterized moduli space for the families perturbation.
Using the characteristic class $\SWbb(E,\fraks) \in H^k(B;\Z/2)$ introduced in \cite{K21}, $\SW(E, \fraks)$ can be written as $\SW(E, \fraks)=\left<\SWbb(E,\fraks), [B]\right>$.
The following \lcnamecref{lem: vanish deg family}  immediately follows from the naturality of the characteristic class:

\begin{lem}
\label{lem: vanish deg family}
Let $B'$ be a closed manifold of dimension $<k$, and $\phi : B \to B'$ be a continuous map.
Let $X \to E' \to B'$ be a $\Diff^+(X)$-bundle whose monodromy preserves $\fraks$.
Then we have
$\SW(\phi^\ast E', \fraks)=0$.
\end{lem}

Let $k_0 \in \{1, \dots, k\}$, and let $p_{k_0} : T^k \to T^k$ be a surjective map obtained as the product of $k_0$-th copies of a degree 2 map $S^1 \to S^1$ and $(k-k_0)$-th copies of the identity of $S^1$.
If $(f_i^2)^\ast$ preserves a $\spinc$ structure $\fraks$ for $i \leq k_0$ and $f_i$ preserves $\fraks$ for $i \geq k_0+1$, the monodromy of
the bundle $p_{k_0}^\ast X_{f_1, \dots, f_k}$ preserves $\fraks$, and thus we can define $\SW(p_{k_0}^\ast X_{f_1, \dots, f_k}, \fraks)$.

For $i_1, \ldots, i_l \in \{1, \ldots, k\}$ with $i_1<\cdots<i_l$ and for $a_{i_1}, \ldots, a_{i_l} \in \Z$, let $(f_{i_1, \ldots, i_l}^{a_{i_1}, \ldots, a_{i_l}})^\ast$ denote the composition of pull-backs $(f_{i_1}^{a_{i_1}})^\ast, \ldots,  (f_{i_l}^{a_{i_l}})^\ast$.

\begin{lem}
\label{lem: cancel lemma}
In the above setup, let $\fraks \in \tau(\Spincns(X,k)/\Conj)$ and
let $k_0 \in \{1,\dots,k\}$, and suppose that
\begin{itemize}
\item $f_i^\ast\fraks\neq\fraks,
(f_i^2)^\ast\fraks=\fraks$ if $i \leq k_0$,
\item $f_i^\ast\fraks=\fraks$ for $i\geq k_0+1$, and 
\item $(f_{1, \dots, k}^{a_{1}, \ldots, a_{k}})^\ast\fraks \in \tau(\Spincns(X,k)/\Conj)$ for all $a_i \in \{0,1\}$.
\end{itemize}
Then we have
\begin{align}
\label{eq: lem pair zero}
\sum_{\{(f_{1, \dots, k}^{a_{1}, \ldots, a_{k}})^\ast\fraks|a_1, \ldots, a_k \in \{0,1\}\}}
\#\calM(X, (f_{1, \dots, k}^{a_{1}, \ldots, a_{k}})^\ast\fraks, \sigma_{(f_{1, \dots, k}^{a_{1}, \ldots, a_{k}})^\ast[\fraks]}) 
= \SW(p_{k_0}^\ast X_{f_1, \dots, f_k}, \fraks)
\end{align}
in $\Z/2$.

Moreover, if there are $i \in \{1, \ldots, k_0\}$ and a smooth isotopy $F_i^t$ between $f_i^2$ and the identity such that $f_1^2, \ldots, F_i^t, \ldots, f_{k_0}^2, f_{k_0+1}, \dots, f_k$ mutually commute for all $t$, then $\SW(p_{k_0}^\ast X_{f_1, \dots, f_k}, \fraks)=0$ in $\Z/2$ and hence the left-hand side of \eqref{eq: lem pair zero} is also zero.
\end{lem}

\begin{proof}
We shall use the notation of \cref{ex: perturbation multi map torus}.
Set $o(f_i)=2$ for $i \leq k_0$ and $o(f_i)=1$ for $i \geq k_0+1$.
We prove \eqref{eq: lem pair zero} by showing that the both hand sides of 
\eqref{eq: lem pair zero} are equal to the count of a parameterized moduli space $\#\calM(X,\fraks,\sigma')$ for some families perturbation $\sigma' : [0,1]^{k} \to \circPi(X)$.

First, by the patching condition \eqref{eq: peruturb general mapping tori} and our assumption \[(f_{1, \dots, k}^{a_{1}, \ldots, a_{k}})^\ast\fraks \in \tau(\Spincns(X,k)/\Conj),\] it follows that, for each $i \leq k_0$, $1 \leq i_1<\dots<i_l \leq k$ and $\mathbf{t} = (t_1, \dots, t_k) \in [0,1]^k$ with $t_i=0$,
\[
f_{i}^\ast(\sigma_{f_{i_1}^\ast\dots f_{i_l}^\ast[\fraks]}(\mathbf{t}))
= \sigma_{f_i^\ast f_{i_1}^\ast\dots f_{i_l}^\ast[\fraks]}(\overline{\mathbf{t}}^i).
\]
Thus the families perturbations $\sigma_{f_{i_1}^\ast\dots f_{i_l}^\ast[\fraks]}$ and $\sigma_{f_i^\ast f_{i_1}^\ast\dots f_{i_l}^\ast[\fraks]}$ can be glued via $f_i^\ast$ along $I_{i,0}^{k-1}$ and $I_{i,1}^{k-1}$, where for $a=0,1$,
\begin{align}
\label{eq: k-1 face}
I_{i,a}^{k-1} = \Set{(t_1,\dots,t_k) \in [0,1]^k| t_i=a}.
\end{align}
\cref{fig: pathching perturbations} illustrates this gluing for $k=k_0=2$.
\begin{figure}
\includegraphics[clip,width = 6.3cm]{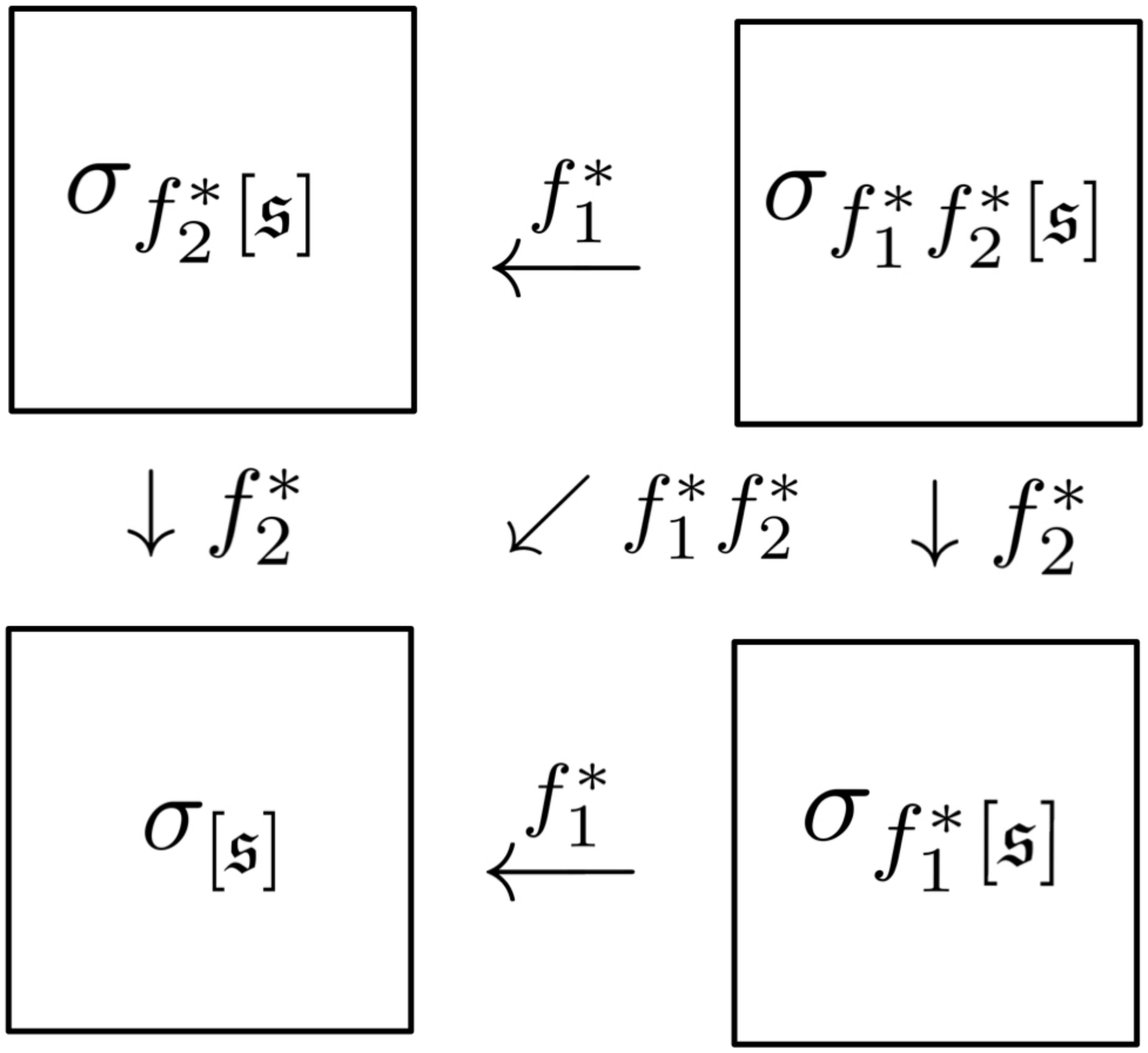}
\caption{Gluing families perturbations for $k=k_0=2$. Families perturbations parameterized by $\bigsqcup_4 [0,1]^2$ are glued and form a families perturbation parameterized by $[0,2]^2$.}
\label{fig: pathching perturbations}
\end{figure}

Therefore the collection
\begin{align}
\label{eq: disjoint per}
\bigsqcup_{\{(f_{1, \dots, k}^{a_{1}, \ldots, a_{k}})^\ast\fraks|a_1, \ldots, a_k \in \{0,1\}\}}
(\sigma_{f_{i_1}^\ast\dots f_{i_l}^\ast[\fraks]} : [0,1]^k \to \circPi(X))
\end{align}
of $2^{k_0}$-tuple of maps is glued via $f_i's$ with $i \leq k_0$ and forms a continuous map $\sigma'' : [0,2]^{k_0} \times [0,1]^{k-k_0} \to \circPi(X)$, and let $\sigma' : [0,1]^{k} \to \circPi(X)$ be the map obtained from $\sigma''$ by rescaling $[0,2]$ by $1/2$.
Thus the left-hand side of \eqref{eq: lem pair zero}, which is the count for the disjoint union \eqref{eq: disjoint per}, coincides with $\#\calM(X,\fraks, \sigma')$.

Next, one may see that, using \eqref{eq: peruturb general mapping tori} and $(f_i^2)^\ast\fraks=\fraks$, the above families perturbation $\sigma'$ satisfies the following condition:
for distinct $i_1, \ldots, i_l \in \{1, \ldots, k\}$ with $i_1<\cdots<i_l$ and $\mathbf{t} = (t_1, \ldots, t_n) \in [0,1]^k$ with $t_{i_1}=\dots=t_{i_l}=0$, we have
\begin{equation}
\label{eq: peruturb general mapping tori deloop}
(f_{i_1, \ldots, i_l}^{o(i_1), \dots, o(i_l)})^\ast\sigma'(\mathbf{t}) =
    \sigma'(\overline{\mathbf{t}}^{i_1, \ldots, i_l}),
\end{equation}
where $\overline{\mathbf{t}}^{i_1, \ldots, i_l} = (t_1',\dots, t_k')$ with $t_j'=t_j$ if $j \notin \{i_1, \ldots, i_l\}$ and $t_j'=1$ if $j \in \{i_1, \ldots, i_l\}$.
Now it is easy to see that $\#\calM(X,\fraks,\sigma')$ coincides with the right-hand side of \eqref{eq: lem pair zero}, which completes the proof of the claim.
This is because, just as in \cref{ex: perturbation multi map torus}, the families perturbation $\sigma'$ with the patching condition \eqref{eq: peruturb general mapping tori deloop} gives rise to a families perturbation for the multiple mapping torus $p_{k_0}^\ast X_{f_1, \dots, f_k} = X_{f_1^2, \ldots, f_{k_0}^2, f_{k_0+1}, \dots, f_k} \to T^k$.

To see the `moreover' part, note that
\[
p_{k_0}^\ast X_{f_1, \dots, f_k} \cong X_{f_1^2, \ldots, f_{k_0}^2, f_{k_0+1}, \dots, f_k}
\cong
X_{f_1^2, \dots, \id, \dots, f_{k_0}^2, f_{k_0+1}, \dots, f_k}
\]
follows from the assumption, where $\id$ is in the $i$-th entry.
Let $\phi : T^k \to T^{k-1}$ be the projection onto the $\{1,\dots,k\} \setminus \{i\}$-factors of $S^1$.
Then the above isomorphism implies that
\[
p_{k_0}^\ast X_{f_1, \dots, f_k} \cong
\phi^\ast X_{f_1^2, \dots, \widecheck{\id}, \dots, f_{k_0}^2, f_{k_0+1}, \dots, f_k},
\]
where $\widecheck{\id}$ indicates that this entry is removed.
It follows from this combined with \cref{lem: vanish deg family} that $\SW(p_{k_0}^\ast X_{f_1, \dots, f_k}, \fraks)=0$.
\end{proof}

\begin{cor}
\label{cor: cancel}
Let $f_1, \ldots, f_k \in \Diff^+(X)$ be mutually commuting diffeomorphisms with the following conditions:
\begin{itemize}
    \item For every $i=1,\dots, k$, $f_i$ preserves $\tau(\Spincns(X,k)/\Conj)$ and the action on $\tau(\Spincns(X,k)/\Conj)$ satisfies the conditions that $f_i^\ast \neq \id$ and  $(f_i^2)^\ast=\id$.
    \item There exist isotopies $(F_{i}^t)_{t\in [0,1]}$ from $f_{i}^2$ to the identity for all $i \in \{1,\dots,k\}$ such that, for any $i\neq j$ and any $t\in [0,1]$, the maps $F_{i}^t$ and $f_{j}$ commute with each other.
    \end{itemize}
Then we have    
\begin{align}
\label{eq: eval sum2}
\langle\SWbbhalftot^k(X_{f_1, \ldots, f_k}), [T^k] \rangle  =
\sum_{\substack{\fraks \in \tau(\Spincns(X,k)/\Conj)\\ f_i^\ast \fraks=\fraks\ (\forall i \in \{1, \dots, k\})}}\SW(X_{f_1, \ldots, f_k}, \fraks).
\end{align}
\end{cor}

\begin{proof}
The diffeomorphisms $f_1, \dots, f_k$ give rise to a $(\Z/2)^k$-action on the set $\tau(\Spincns(X,k)/\Conj)$.
Given $\fraks \in \tau(\Spincns(X,k)/\Conj)$, if there is $i \in \{1, \ldots, k\}$ such that $f_i^\ast\fraks\neq \fraks$, we may apply \cref{lem: cancel lemma} to $\fraks$ (by renaming indices of $f_i$) and obtain that the contribution of the $(\Z/2)^k$-orbit of $\fraks$ in the sum of the right-hand side of \eqref{eq: eval sum} is zero.

As in \cref{ex: trivial monodromy,ex: perturbation single mapping torus 1}, if a $\spinc$ structure $\fraks \in \tau(\Spincns(X,k)/\Conj)$ is preserved by the monodromy of $X_{f_1, \dots, f_k}$, the families perturbation $\sigma_{[\fraks]}$ can be taKan so that $\#\calM(X,\fraks,\sigma_{[\fraks]})=\SW(X_{f_1, \dots, f_k},\fraks)$.
Thus the claim of the \lcnamecref{cor: cancel} follows from \eqref{eq: eval sum}.
\end{proof}

Another important input in calculating the Seiberg--Witten characteristic class is a gluing formula for the numerical families Seiberg--Witten invariant.
To state this, we need to introduce a few notations.
In general, for an oriented closed 4-manifold $X$ and for an oriented fiber bundle $X \to E \to B$ with fiber $X$, a vector bundle $\R^{b^+(X)} \to H^+(E) \to B$ is associated (e.g. \cite{LiLiu01, BK20}).
Roughly, $H^+(E)$ is defined by the action of $\Diff^+(X)$ on the Grassmannian of maximal-dimensional positive-definite subspaces of $H^2(X;\R)$ together with the transition functions of $E$.
While $H^+(E)$ is determined by $E$ only up to isomorphism, once we take a fiberwise metric on $E$, the space of fiberwise self-dual harmonic 2-forms gives a model of $H^+(E)$.

Now let $X$ be an oriented closed smooth 4-manifold with $b^+(X) \geq 2$ and $b_1(X)=0$, and let $B$ be a closed smooth manifold $B$ of dimension $k>0$.
Let $X \to E_1 \to B$ denote the trivial bundle: $E_1= B \times X$.
Let $kS^2 \times S^2 \to E_2 \to B$ be a $\Diff^+(kS^2 \times S^2)$-bundle over $B$ and suppose that $E_2$ admits a section whose normal bundle is trivial.
By choosing such a section and a trivialization of its normal bundle, we may form a fiberwise connected sum $X\#k(S^2\times S^2) \to E_1\#_f E_2 \to B$.
Let $\fraks_S$ denote the spin structure on $kS^2\times S^2$.
For a $\spinc$ structure $\frakt$ with $d(\frakt)=0$ on $X$, the connected sum $\frakt\#\fraks_S$ is invariant under the monodromy of $E_1\#_f E_2$ and we have $d(\frakt\#\fraks_S)=-k$. Thus the numerical families invariant $\SW(E_1\#_f E_2, \frakt\#\fraks_S) \in \Z/2$ makes sense.

\begin{pro}
\label{prop: gluing BK}
In the above setup, suppose that the Stiefel--Whitney class $w_{k}(H^+(E_2)) \in H^k(B;\Z/2)$ is non-zero.
Then, for a spin$^c$ structure $\frakt$ on $X$ with $d(\frakt)=0$, we have
\[
\SW(E_1\#_f E_2, \frakt\#\fraks_S)
= \SW(X,\frakt)
\]
in $\Z/2$.
In particular, we have
\[
\sum_{\frakt \in \tau_0(\Spincns(X,0)/\Conj)}\SW(E_1\#_f E_2, \frakt\#\fraks_S)
= \SWbbhalftot^0(X)
\]
in $\Z/2$.
\end{pro}

\begin{proof}
This is a direct application of
a families gluing formula proven in \cite[Theorem~1.1]{BK20}.
\end{proof}

\subsection{Key calculation}
\label{subsection Multiple mapping torus}

In this \lcnamecref{subsection Multiple mapping torus}, we calculate the Seiberg--Witten characteristic class for a special type of  multiple mapping tori (\cref{thm key computation source}).
This is the key computation to prove main results explained in \cref{section Introduction}.

We again start with some notations. Let $W$ be a smooth 4-manifold with a single boundary component. We use $W_{S}$ to denote the once stabilized manifold
\[
W \cup_{\del W} \left(([0,1] \times \del W)\#S^2\times S^2\right).
\]
Here $([0,1] \times \del W)\#S^2\times S^2$ denotes the inner connected sum.
For $k\geq 1$, we use  $W_{kS}$ to denote the $k$-th stabilized manifold: we use  $W_{kS}$ instead of $W\#kS^2\times S^2$ when we want to specify the way of doing connected sums. Note that we have natural stabilization maps
\begin{equation}\label{eq: stabilization map}
s_{W}: \BDiff_{\del}(W) \rightarrow \BDiff_{\del}(W_{S}) \end{equation}
and
\begin{equation}\label{eq: stabilization map  discrete}
s_{W}^{\delta}: \BDiff_{\del}(W_{kS})^{\delta} \rightarrow \BDiff_{\del}(W_{(k+1)S})^{\delta}
\end{equation}
defined by extending a diffeomorphism with the identity map. We write $s^{\delta}$ and $s$ when $W$ is clear from the context.

We can also define the stabilization of a closed manifold $X$. We pick a smooth embedded $D^{4}$ in a local chart and consider the punctured manifold $\mathring{X}:=X\setminus
\Int(D^4)$. Then we define \[
X_{kS}:=\mathring{X}_{kS}\cup_{S^3}D^4.
\]
Note that the stabilization maps \eqref{eq: stabilization map} and \eqref{eq: stabilization map discrete} are not well-defined for closed $X$. 
Next, we discuss a general construction of multiple mapping tori with fiber $W_{kS}$. Let $S=S^2 \times S^2 \setminus \Int(D^4)$ and let $f \in \Diff_{\del}(S)$ be a relative diffeomorphism that acts as the $(-1)$-multiplication on $H_2(S;\mathbb{Z})$.
Such $f$ can be obtained from the diffeomorphism defined as the product of two copies of the reflection on $S^2$ about the equator, after deforming near $D^4$ by isotopy. Note that $f^2$ is the boundary Dehn twist on $S^2 \times S^2 \setminus \Int(D^4)$. Using the circle action on $S^{2}\times S^2$, one finds a smooth isotopy $\{F^{t}\}_{t\in [0,1]}$ from $f^2$ to the identity.

Note that there are $k$ disjoint copies of embedded $S$ in $W_{kS}$, which we denote by $S_{1},\cdots, S_{k}$. By extending the diffeomorphism $f$ on $S_{i}$ with the identity map on $W\setminus S_{i}$, we obtain diffeomorphisms 
\[
f_{1},\cdots, f_{k}:W_{kS}\rightarrow W_{kS}
\]
that mutually commute with each other.
Furthermore, by extending the isotopy $\{F^{t}\}_{t\in [0,1]}$ on $S_{i}$ with the identity map, we obtain a smooth isotopy 
\begin{equation}\label{eq: isotopy Ft}
\{F^{t}_{i}:W_{kS}\rightarrow W_{kS}\}_{t\in[0,1]}    
\end{equation}
 from $f^2_{i}$ to the identity map, which is supported in $S_{i}$.

Now we take the multiple mapping torus 
\begin{equation}\label{eq: multiple mapping torus}
W_{kS}\rightarrow (W_{kS})_{f_{1},f_2,\cdots,f_{k}}\rightarrow T^{k}
\end{equation}
and denote it by $E^{k}(W)$. The fiber bundle $E^{k}(W)$ carries a natural flat structure and hence induces a map 
\[
E^{k}(W): T^{k}\rightarrow \BDiff_{\partial}(W_{kS})^{\delta}.
\]

\begin{lem}
\label{lem: flat 2-torsion}
For any $W$ with zero or one boundary component, one has 
\[2\cdot E^{k}(W)_{*}[T^{k}]=0\in H_{k}(\BDiff_{\del}(W_{kS})^{\delta};\mathbb{Z}).\]
\end{lem}
\begin{proof}
Consider the multiple mapping torus 
\[
W_{kS}\rightarrow (W_{kS})_{f^2_{1},f_2,\cdots,f_{k}}\rightarrow T^{k}.
\]
It is the pull-back of $E^{k}(W)$ under a degree-$2$ map $T^{k}\rightarrow T^{k}$. Hence we have
\[2\cdot E^{k}(W)_{*}[T^{k}]=((W_{kS})_{f^2_{1},f_2,\cdots,f_{k}})_{*}[T^{k}].\]
We can apply Lemma \ref{lem: 2-torsion} to finish the proof.
\end{proof}

\begin{lem}\label{lem: Ek stabilization}
For any $W$ with a single boundary component, the following diagram commutes up to homotopy \begin{align}\label{diagram: stabilization}
\begin{split}
\xymatrix{
     & \BDiff_{\del}(W_{(k+1)S})^{\delta} \\
    T^k \ar[ru]^-{E^{k}(W_{S})} \ar[r]_-{E^{k}(W)} & \BDiff_{\del}(W_{kS})^{\delta}. \ar[u]_{s_{W_{kS}}^{\delta}}
   }
\end{split}   
\end{align}
Furthermore, let $W'$ be  another smooth manifold with boundary. Then $E^{k}(W)$ and $E^{k}(W')$ are isomorphic as smooth bundles with flat structures if $W$ and $W'$ are diffeomorphic.
Similarly, $E^{k}(W)$ and $E^{k}(W')$ are isomorphic as topological bundles with flat structures if $W$ and $W'$ are homeomorphic.
\end{lem}
\begin{proof}
This is a straightforward checking using definition.
\end{proof}
Now we set $W$ to be a closed manifold $X$ and discuss the families Seiberg--Witten invariant. 
\begin{lem}\label{lem: SW-halftot for mapping tori}
For any oriented closed smooth 4-manifold $X$ with $b^{+}(X)\geq k+2$, we have
\[ 
\langle\SWbbhalftot^{k}(E^{k}(X)),[T^{k}]\rangle=\SWbbhalftot^{0}(X).\]
\end{lem}
\begin{proof} 
First, we fix a section 
\[
\tau_0 : \Spincns(X)/\Conj \to \Spincns(X).
\]
Then, for given $[\fraks] \in \Spincns(X_{kS})/\Conj$, we set $\tau([\fraks])$ to be the $\spinc$ structure $\fraks'$ with $[\fraks']=[\fraks]$ that satisfies $\fraks'|_{X} = \tau_0([\fraks'|_{X}])$.
This defines a section \[\tau_1 : \Spincns(X_{kS})/\Conj \to \Spincns(X_{kS}),\] and
restricting this, we obtain a section \[\tau : \Spincns(X_{kS},k)/\Conj \to \Spincns(X_{kS},k).\]

Note that the isotopy $\{F^{t}_{i}\}$ satisfies the conditions of \cref{cor: cancel}. 
Thus we can reduce the calculation of $\SWbbhalftot^k(E')$ to monodromy invariant $\spinc$ structures;
By \cref{cor: cancel}, we have  
\begin{align}
\label{eq: eval sum4}
\langle\SWbbhalftot^k(E^k(X)), [T^k] \rangle  =
\sum_{\substack{\fraks \in \tau(\Spincns(X_{kS},k)/\Conj)\\ f_i^\ast \fraks=\fraks\ (\forall i \in \{1, \dots, k\})}}\SW(E^{k}(X), \fraks).
\end{align}

On the other hand, the action of $f_{i}$ on the second homology group equals $-1$ on $(S^{2}\times S^{2})_{i}$ and equals $1$ on the other summands. From this, we see that a $\spinc$ structure on $X_{kS}$ is preserved by all $f_{i}$ if and only if it can be written as the connected sum $\mathfrak{t}\#\mathfrak{s}_{S}$ between a $\spinc$ structure $\mathfrak{t}$ on $X$ and the spin structure $\mathfrak{s}_{S}$ on $k(S^2\times S^2)$.
Thus we can apply \cref{prop: gluing BK} and get 
\[
\begin{split}
\sum_{\substack{\fraks \in \tau(\Spincns(X_{kS},k)/\Conj)\\ f_i^\ast \fraks=\fraks\ (\forall i \in \{1, \dots, k\})}}\SW(E^{k}(X), \fraks)&=\sum_{\frakt \in \tau(\Spincns(X,0)/\Conj)}\SW(E^{k}(X), \frakt\#\fraks_{S})\\
&=\sum_{\frakt \in \tau(\Spincns(X,0)/\Conj)}\SW(X, \frakt)
=\SWbbhalftot^0(X).
\end{split}
\]
This combined with \eqref{eq: eval sum4} finishes the proof.
\end{proof}

Next, let us consider a natural map 
\begin{equation}\label{eq: extension map}
\rho : \BDiff_{\del}(\mathring{X}) \to \BDiff^+(X)
\end{equation}
induced from a map extending a relative diffeomorphism by the identity of $D^4$.
We call $\rho$ the extension map, and drop $X$ from our notation.
To state the main result of this \lcnamecref{subsection Multiple mapping torus}, we need one more definition.

\begin{defi}
Let $X$ be a simply-connected, closed, smooth, oriented 4-manifold. We say that $X$ dissolves if it is diffeomorphic to $a\mathbb{CP}^{2}\#b\overline{\mathbb{CP}^2}$ or $c(S^2\times S^2)\#d K3$ for some $a,b,c\in \mathbb{N}$ and $d\in \mathbb{Z}$. We say that $X$ dissolves after $n$ stabilizations if $X\#n(S^{2}\times S^2)$ dissolves. 
\end{defi}

Now we are ready to state the main result of this \lcnamecref{subsection Multiple mapping torus}.

\begin{thm}
\label{thm key computation source}
Let $k > 0$ and let $X$ be a simply-connected closed oriented indefinite smooth 4-manifold with $b^+(X) \geq k+2$.
Assume that $\SWbbhalftot^0(X)=1$ in $\Z/2$, $X$ dissolves after one stabilization, and $X$ is not homotopic to $K3$.
Then there exist an element \[\alpha^{\delta}_{k}(X)\in H_{k}(\BDiff_{\del}(\mathring{X}_{kS})^{\delta};\mathbb{Z}) \]
that satisfies the following properties:
\begin{enumerate}[label=(\roman*)]
\item $\alpha^{\delta}_{k}(X)$ is $2$-torsion,
\item $\alpha^{\delta}_{k}(X)$ belongs to the kernel of the stabilization map 
\[
s^{\delta}_{*}:  H_{k}(\BDiff_{\del}(\mathring{X}_{kS})^{\delta};\mathbb{Z})\rightarrow  H_{k}(\BDiff_{\del}(\mathring{X}_{(k+1)S})^{\delta};\mathbb{Z}).\]
\item Let $\alpha_{k}(X)$ be the image of $\alpha^{\delta}_{k}(X)$ under the forgetful map 
\[
H_{k}(\BDiff_{\del}(\mathring{X}_{kS})^{\delta};\mathbb{Z})\rightarrow H_{k}(\BDiff_{\del}(\mathring{X}_{kS});\mathbb{Z}).
\] Then 
\[0 \neq \langle \SWbbhalftot^{k}(X_{kS}),\rho_{*}(\alpha_{k}(X))\rangle\in \mathbb{Z}/2\]
\item $\alpha_{k}(X) $ and $\alpha_k^{\delta}(X)$ are both  topologically trivial. Namely, we have 
\[
\begin{split}
\alpha^{\delta}_{k}(X) \in \ker(H_k(\BDiff_\del(\mathring{X}_{kS})^{\delta};\Z) &\to H_k(\BHomeo_\del(\mathring{X}_{kS})^{\delta};\Z)),\\
\alpha_{k}(X) \in \ker(H_k(\BDiff_\del(\mathring{X}_{kS});\Z) &\to H_k(\BHomeo_\del(\mathring{X}_{kS});\Z)).
\end{split}
\]
\end{enumerate}
\end{thm}

\begin{rmk}
\label{rmk: MMM class}
\cref{thm key computation source} implies that $\SWbbhalftot$ detects topologically trivial non-zero homology class of $\BDiff^+(X)$ in $\Z$ or $\Z/2$-coefficient.
On the other hand, the (generalized) Mumford--Morita--Miller classes, which are basic characteristic classes for smooth fiber bundles, are generalized also for topological fiber bundles over $\mathbb{F} = \Q$ or $\Z/2$ \cite{ERW14}.
Thus the Mumford--Morita--Miller classes do not detect kernels of $H_\ast(\BDiff^+(X);\mathbb{F}) \to H_\ast(\BHomeo^+(X);\mathbb{F})$.
\end{rmk}

\begin{proof}[Proof of \cref{thm key computation source}]By the assumption on the dissolving and that $X$ is not homotopic to $K3$, there is a 4-manifold $X'$ such that the following conditions hold:
\begin{enumerate}
   \item $X\# S^2\times S^2$ is diffeomorphic to   $X'\# S^2\times S^2$.
    \item $X'$ is homeomorphic to $X$.
    \item The Seiberg--Witten invariant vanishes for any $\spinc$ structure on $X'$. In particular, this implies $\SWbbhalftot^{0}(X')=0$.
\end{enumerate}

Consider the bundles $E^{k}(\mathring{X})$ and $E^{k}(\mathring{X}')$. Since $X'_{kS}$ is diffeomorphic to $X_{kS}$, we have two maps
\[
E^{k}(\mathring{X}), E^{k}(\mathring{X}'): T^k\rightarrow \BDiff_{\del}(\mathring{X}_{kS})^{\delta}.
\]
Then we set 
\[\alpha^{\delta}_{k}(X):=E^{k}(\mathring{X})_{*}[T^{k}]-E^{k}(\mathring{X}')_{*}[T^{k}].\]
Now we verify the four properties in the statement of the \lcnamecref{thm key computation source} as follows:
\begin{itemize}
    \item (i) directly follows from \cref{lem: flat 2-torsion}. 
    \item  By 
\cref{lem: Ek stabilization}, we have \[
s^{\delta}_{*}E^{k}(\mathring{X})_{*}[T^{k}]=E^{k}(\mathring{X}_{S})_{*}[T^k],\quad s^{\delta}_{*}E^{k}(\mathring{X}')_{*}[T^{k}]=E^{k}(\mathring{X}'_{S})_{*}[T^k].
\]
Since $\mathring{X}_{S}$ is diffeomorphic to $\mathring{X}'_{S}$, (ii) is proved.
\item (iii) is equivalent to 
\begin{equation*}
\langle \SWbbhalftot^{k}(E^{k}(X)),[T^{k}]\rangle \neq \langle \SWbbhalftot^{k}(E^{k}(X')),[T^{k}]\rangle.
\end{equation*}
By \cref{lem: SW-halftot for mapping tori}, this is equivalent to 
\[
\SWbbhalftot^{0}(X)\neq \SWbbhalftot^{0}(X'),
\]
which follows from our assumption and $\SWbbhalftot^{0}(X')=0$.
\item 
Since $X$ is homeomorphic to $X'$, $E^{k}(\mathring{X})$ and $E^{k}(\mathring{X}')$ are isomorpic as topological bundles with flat structure. Therefore the maps $E^{k}(\mathring{X})$ and $ E^{k}(\mathring{X}')$ are homotopic after being composed with the forgetful map
\[
\BDiff_\del(\mathring{X}_{kS})^{\delta}\rightarrow \BHomeo_\del(\mathring{X}_{kS})^{\delta}.
\]
This proves that $\alpha^{\delta}_{k}(X)$ is topologically trivial, which implies that $\alpha_{k}(X)$ is topologically trivial as well.
\end{itemize}
\end{proof}

\begin{cor}
\label{thm: Diff Homeo closed 4-manifold}
Let $k > 0$ and let $X$ be a simply-connected closed oriented indefinite smooth 4-manifold with $b^+(X) \geq k+2$.
Assume that $\SWbbhalftot^0(X)=1$ in $\Z/2$, $X$ dissolves after one stabilization, and $X$ is not homotopic to $K3$. Consider the forgetful map $i:\Diff^+(X_{kS})\rightarrow \Homeo^{+}(X_{kS})$.
Then the induced map 
\[
i_{*}:H_k(\BDiff^+(X_{kS});\Z)
\to H_k(\BHomeo^+(X_{kS});\Z)
\]
is not injective and the induced map 
\[
i^{*}:H^k(\BHomeo^+(X_{kS});\Z/2)\rightarrow H^k(\BDiff^+(X_{kS});\Z/2)
\]
is not surjective. Analogous results hold for the maps $i^{\delta}_{*}$ and $i_{\delta}^{*}$ induced by the inclusion $i^\delta:\Diff^+(X_{kS})^{\delta}\hookrightarrow \Homeo^{+}(X_{kS})^{\delta}$.
\end{cor}
\begin{proof}
Set $\alpha=\rho_{*}(\alpha_{k}(X))\in H_k(\BDiff^+(X_{kS});\Z)$.
Then $\alpha\neq 0$ by \cref{thm key computation source} (ii). By \cref{thm key computation source} (iv),  $\alpha$ maps to zero in  $H_k(\BHomeo^+(X_{kS});\Z)$. This implies that $\SWbbhalftot^k(X_{kS})$ is not in the image of $i^{*}$. By using $\alpha^{\delta}_{k}(X)$ instead of $\alpha_{k}(X)$, one can prove analogous results for $i^{\delta}_{*}$ and $i_{\delta}^{*}$.
\end{proof}

\section{Moduli space of formally smooth topological bundles}\label{section BDiffL} The purpose of this section is to define the moduli space $\BHomeo^{+}_{L}(X)$, which classifies topological bundles whose vertical tangent microbundle is lifted to a vector bundle. This moduli space is essential in our discussion of exotic phenomena that is special in dimension 4. We start by recalling some  definitions about microbundles (see Milnor \cite{Milnormicrrobundles}).

\begin{defi} Let $B$ be a topological space. A {\it microbundle} over $B$ of rank $n$ is a triple $\xi=(Y,p,s)$ of a space $Y$ and continuos maps $p:Y\rightarrow B$ and $s: B\rightarrow Y$ such that 
\begin{itemize}
    \item $p\circ s=\id$.
    \item For any point $x\in B$, there exists a neighborhood $U\subset B$ of $x$, a neighborhood $V\subset Y$ of $s(U)$, and a homeomorphism $\phi: \mathbb{R}^{n}\times U\xrightarrow{\cong} V$ that make the following diagram commute:
\[
\xymatrix{ U\ar@{=}[d]\ar[r]^-{s_0} &\mathbb{R}^{n}\times U\ar[r]^-{p_2}\ar[d]^{\phi}& U\ar@{=}[d]\\
U\ar[r]^{s} &V\ar[r]^-{p}& U}.
\]
Here $s_{0}(x)=(0,x)$ and $p_2(t,x)=x$.
\end{itemize}
\end{defi}
\begin{defi}\label{defi: vertical tangent microbundle} Let $X\hookrightarrow E\rightarrow B$ be a topological bundle (i.e. a fiber bundle whose structure group is in $\Homeo(X)$). The \emph{vertical tangent microbundle} $\scrT(E/B)$ is the microbundle over $E$ given by the triple $(E\times_{B}E,p,s)$. Here $p: E\times_{B}E\rightarrow E$ is defined as
$p(x_{1},x_{2}):=x_{2}$ and $s: E\rightarrow E\times_{B}E$ is defined as $s(x):=(x,x)$.
\end{defi}

\begin{defi} Given two microbundles $\xi=(Y,p,s)$ and $\xi'=(Y',p',s')$ over  $B$, an {\it isomorphism} $\rho: \xi\rightarrow \xi'$ is a homeomorphism $V\xrightarrow{\cong} V'$ between a neighborhood $V\subset Y$ of $s(B)$ and a neighborhood $V'\subset Y'$ of $s'(B)$ that is compatible with the maps $p,p',s,s'$. We treat two isomorphisms as the same if they coincide on a neighborhood of $s(B)$.
\end{defi}
\begin{defi} A {\it linear structure} on a microbundle $\xi$ is an isomorphism $\rho:W\rightarrow \xi$ from some vector bundle $W$. We treat two linear structures $\rho: W\rightarrow \xi$ and $\rho': W'\rightarrow \xi$ as the same if there exists an isomorphism $\psi: W\rightarrow W'$ of vector bundles such that $\rho\circ \psi=\rho'$. 
\end{defi}

To avoid delicate point-set topology issues in our construction of moduli spaces, we use simplical sets, which has been a standard tool in studying moduli spaces of manifolds. For readers' convenience, we summarize some key properties here and refer to \cite{Goerss09} for a detailed introduction. 
Let $\bf{\Delta}$ be the simplex category whose objects are the ordered sets 
$[n]:=\{0,\cdots, n\}$ and whose morphisms are order preserving maps between these ordered sets. A simplical set $S$ is a contravariant functor
\[
S: \bf{\Delta}\rightarrow \bf{Set}
\] 
to the category of sets. That means we associate to $[k]$ a set $S[k]$ (called the set of $k$-simplices) and associate a structure map 
\[\iota^{*}: S[l]\rightarrow S[k].\]
to an order preserving map $\iota:[k]\rightarrow [l]$. Any simplical set $S$ has a decomposition into its path-components.

A map $S\rightarrow T$ between two simplical sets is by definition a natural transformation between the functors. We say the map is a Kan fibration if it satisfies the so-called ``horn-lifting property''. Let $*$ be the simplical set with a single $k$-simplex for each $k$. We say $S$ is a Kan complex if the constant map $K\rightarrow *$ is a Kan fibration. Kan complexes are those simplical sets which resemble the simplical set of singular simplices of a topological space. In particular, given a Kan complex $S$ and a $0$-simplex $x$, one can define its homotopy group $\pi_{k}(S,x)$ for any $k$. A map $S\rightarrow T$ between Kan complexes is called a weak equivalence if it induces a bijection between path-components and induces an isomorphism on all homotopy groups. A Kan complex $S$ is weakly contractible if it is path-connected and has vanishing homotopy groups.


A simplical group is a contravariant functor 
\[
G: \bf{\Delta}\rightarrow \bf{Group}
\]
to the category of groups. A left action (resp. right action) of a simplical group $G$ on a simplical set $S$ is a left action (resp. right action) of $G[k]$ on $S[k]$ for each $k$. These actions are required to be compatible with structure maps. We say the action is levelwise free if the action of $G[k]$ on $S[k]$ is free for each $k$. And we define the quotient simplical set $S/G$ by setting $(S/G)[k]:=S[k]/G[k]$.

Given a simplical set $S$ with a right $G$-action and simplical set $T$ with a left $G$-action, we define  $S\times_{G}T$ by setting
\[
(S\times_{G}T)[k]:=S[k]\times_{G[k]} T[k].
\]

There is a way to associate to a simplical set $S$ a CW complex $|S|$, called the \emph{geometric realization} of $S$.  We briefly recall its construction. For each $k\in \mathbb{N}$, we consider 
the standard $k$-simplex 
\[
\Delta^{k}:=\{(x_0,\cdots,x_{k})\mid x_{i}\geq 0, \sum^{n}_{i=0}x_{i}=1\}\subset \mathbb{R}^{k+1}.
\]
Vertices of $\Delta^{k}$ are canonically identified with elements in $[k]$. Given any order preserving map $\tau: [l]\rightarrow [k]$, we have an induced simplicial map $\Delta^\tau: \Delta^{l}\rightarrow \Delta^{k}.$ The geometric realization is defined as the quotient space 
\[
|S|:=(\bigsqcup_{k\in \mathbb{N}}(\Delta^{k}\times S[k]))/\sim.
\]
Here we use discrete topology on $S[k]$ and the equivalence relation is generated by 
\[
(\Delta^{\tau}(x),\sigma)\sim (x, \tau^*(\sigma))
\]
for any morphism $\iota: [k]\rightarrow [l]$, any $x\in \Delta^{k}$ and any $\sigma\in S[l]$. This is a functorial construction, meaning that  any map $f:S\rightarrow T$ between simplical sets induces a continuous map $|f|: |S|\rightarrow |T|$. We note that the homotopy group of a Kan complex coincides with its geometric realization.

Now we start defining various simplical sets associated to a connected, closed, oriented smooth manifold $X$. For now, we do not assume that $X$ has dimension $4$. 

\begin{defi}
The simplical homeomorphism group $\Homeo^{+}(X)$ is defined by setting $\Homeo^{+}(X)[k]$ to be the group of orientation-preserving self-homeomorphisms on $\Delta^{k}\times X$ that covers $\id_{\Delta^{k}}$. Similarly, we define the simplical diffeomorphism group $\Diff^{+}(X)$ by setting $\Diff^{+}(X)[k]$ to be the group of orientation-preserving self-diffeomorphims on $\Delta^{k}\times X$ that covers $\id_{\Delta^{k}}$. Here we treat $\Delta^{k}$ as a smooth manifold with corner. Given a morphism $\iota: [k]\rightarrow [l]$ in $\bf{\Delta}$, the associated structure maps of $\Homeo^{+}(X)$ and  $\Diff^{+}(X)$ are both defined using the bundle map 
\[
\Delta^{\tau}\times \id_{X}: \Delta^{k}\times X\rightarrow \Delta^{l}\times X.
\]
\end{defi}

\begin{defi}
The simplical set $\Emb(X,\mathbb{R}^{\infty})$ is defined as follows. For each $k\in \mathbb{N}$ and $n\in \mathbb{N}$, we let $\Emb(X,\mathbb{R}^{n})[k]$ be the set of locally flat topological embeddings $\Delta^{k}\times X\hookrightarrow \Delta^{k}\times \mathbb{R}^{n}$ that covers $\id_{\Delta^{k}}$. Then we define 
\[\Emb(X,\mathbb{R}^{\infty})[k]:=\underset{n\rightarrow \infty}{\operatorname{colimit}{}}\Emb(X,\mathbb{R}^{n})[k].
\]
Here the colimit is defined using standard inclusion $\mathbb{R}^{n}\hookrightarrow \mathbb{R}^{n+1}$. Given a morphism $\tau:[k]\rightarrow [l]$ in $\bf{\Delta}$, the structure map
\[
\tau^*: \Emb(X,\mathbb{R}^{\infty})[l]\rightarrow \Emb(X,\mathbb{R}^{\infty})[k]
\]
is defined as follows. Given an element in $\Emb(X,\mathbb{R}^{\infty})[l]$ represented by an embedding $f: \Delta^{k}\times X\hookrightarrow \Delta^{k}\times \mathbb{R}^{n}$, we let $\tau^{*}(f)\in \Emb(X,\mathbb{R}^{\infty})[k]$ be the unique embedding that makes the following diagram commutes. 
\begin{equation}
\label{diagram: pullback embedding}
\begin{split}
\xymatrix{
\Delta^{k}\times X\ar[r]^{\tau^{*}(f)}\ar[d]_{\Delta^{\tau}\times \id_{X}}& \Delta^{k}\times \mathbb{R}^{n}\ar[d]^{\Delta^{\tau}\times \id_{\mathbb{R}^{n}}}\\
\Delta^{l}\times X \ar[r]^{f} & \Delta^{l}\times \mathbb{R}^{n}.}
\end{split}
\end{equation}
\end{defi}

There is a levelwise free right action of $\Homeo^{+}(X)$ on $\Emb(X,\mathbb{R}^{\infty})$ defined by precomposing an embedding $f\in \Emb(X,\mathbb{R}^{\infty})[k]$ with a self-diffeomorphism $h\in \Homeo^{+}(X)[k]$. Via the forgetful map $\Diff^{+}(X)\rightarrow \Homeo^{+}(X)$, we also get a levelwise free right action of $\Diff^{+}(X)$ on $\Emb(X,\mathbb{R}^{\infty})$.

\begin{defi}
The simplical set $L(\scrT X)$, which corresponds to the space of linear structures on $\scrT X$, is defined as follows. For each $k\in \mathbb{N}$, we consider the vertical tangent microbundle $\scrT(\Delta^{k}\times X/\Delta^{k})$ of the trivial bundle $\Delta^{k}\times X\rightarrow \Delta^{k}$. We let $L(\scrT X)[k]$ be the set of linear structures on $\scrT(\Delta^{k}\times X/\Delta^{k})$. Given a morphism $\tau:[k]\rightarrow [l]$ in $\bf{\Delta}$, we define the 
structure map
\[\tau^{*}: L(\scrT X)[l]\rightarrow L(\scrT X)[k]
\]
by pulling back a linear structure via the bundle map 
\[
\Delta^{\tau}\times \id_{X}: \Delta^{k}\times X\rightarrow \Delta^{l}\times X.
\]
\end{defi}

\begin{defi}
The simplical set $R(TX)$, which corresponds to the space of Riemannian metrics on $TX$, is defined as follows. For each $k\in \mathbb{N}$, let $R(TX)[k]$ be the set of smooth metric on the vertical tangent bundle $T(\Delta^{k}\times X/\Delta^{k})$. Given a morphism $\tau:[k]\rightarrow [l]$ in $\bf{\Delta}$, we define the 
structure map
\[\tau^{*}: R(T X)[l]\rightarrow R(T X)[k] 
\]
by pulling back a metric via the bundle map 
\[\Delta^{\tau}\times \id_{X}: \Delta^{k}\times X\rightarrow \Delta^{l}\times X.
\]
\end{defi}

\begin{lem}\label{lem: Ken} The following conclusion holds.
\begin{enumerate}
    \item $\Homeo^{+}(X), \Diff^{+}(X)$ are Kan complexes.
    \item $L(\scrT X)$ is a Kan complex.
    \item $\Emb(X,\mathbb{R}^{\infty})$ is a weakly contractible Kan complex.
    \item  $R(TX)$ is a weakly contractible Kan complex.
    \item $\Emb(X,\mathbb{R}^{\infty})\times_{\Homeo^{+}(X)} L(\scrT X)$ is a Kan complex.
    \item $\Homeo^{+}(X)\times_{\Diff^{+}(X)}R(TX)$ is a Kan complex.
    \item $\Emb(X,\mathbb{R}^{\infty})/\Homeo^{+}(X)$ is a Kan complex and its geometric realization is a model for   $\BHomeo^{+}(X)$.
    \item $\Emb(X,\mathbb{R}^{\infty})\times_{\Diff^{+}(X)}R(TX)$ is a Kan complex and its geometric realization is a model for $\BDiff^{+}(X)$.
\end{enumerate}
\end{lem}
\begin{proof} The underlying simplical set of any simplical group is a Kan complex \cite[Lemma I.3.4]{Goerss09}, which proves (1).

We let $\Lambda^{n}_{i}$ be the ``horn'' obtained by removing interior of the $i$-th face from $\partial\Delta^n$. By the definition of Kan complex, (2) is equivalent to showing that any linear structure $\rho\in \scrT(\Lambda^{n}_{i}\times X/\Lambda^{n}_{i})$ can be extended to a linear structure $\rho'$ on $\scrT(\Delta^{n}\times X/\Delta^{n})$. This is straightforward since we can construct $\rho'$ as the pullback of $\rho$ via the continuous retraction $\Delta^n\rightarrow   \Lambda^{n}_{i}$. 

(3) is proved as in \cite[Lemma 2.2]{kupers2015}.

The simplicial set $R(TX)$ is a Kan complex since any metric on $T(\Lambda^{n}_{i}\times X/\Lambda^{n}_{i})$ can be extended to a metric on $T(\Delta^{n}\times X/\Delta^{n})$. And $R(TX)$ is weakly contractible since any metric on $T(\partial\Delta^n\times X/\partial\Delta^n)$ can be extended to a metric on $T(\Delta^{n}\times X/\Lambda^{n}_{i})$. This proves (4).

By (1) and (4), the product $\Emb(X,\mathbb{R}^{\infty})\times L(\scrT X)$ is a Kan complex. By \cite[Corollary V.2.6]{Goerss09}, the map \[\Emb(X,\mathbb{R}^{\infty})\times L(\scrT X)\rightarrow \Emb(X,\mathbb{R}^{\infty})\times_{\Homeo^{+}(X)} L(\scrT X)\] is a surjective Kan fibration. These two facts together imply that $\Emb(X,\mathbb{R}^{\infty})\times_{\Homeo^{+}(X)} L(\scrT X)$ is also a Kan complex, namely (5) (See \cite[Proof of Lemma V.3.7]{Goerss09}.)

(6) is proved by an argument similar to that used to show (5).

(7) is proved as in \cite[Lemma 2.7]{kupers2015}. And (8) can be proved similarly.
\end{proof}

\begin{defi} A {\it formally smooth $X$-bundle} is a topological fiber bundle $X\rightarrow E\rightarrow B$ equipped with a linear structure on its vertical tangent microbundle. Given two formally smooth $X$-bundles $\theta_0, \theta_1$ over $B$, we say that they are {\it isomorphic} if there is an isomorphism between the bundles that respects the linear structures. And we say that they are {\it concordant} if there exists a formally smooth $X$-bundle $\theta$ over $I\times B$ such that $\theta|_{\{i\}\times X}$ is isomorphic to $\theta_i$ for $i=0,1$.
\end{defi}

\begin{defi}\label{defi: moduli space for formally smooth bundle}
We define the {\it moduli space of formally smooth manifold $X$} as follows \[\BHomeo^{+}_{L}(X):=|\Emb(X,\mathbb{R}^{\infty})\times_{\Homeo^{+}(X)} L(\scrT X)|.\]
\end{defi}

Definition \ref{defi: moduli space for formally smooth bundle} is justified by the following proposition.  
\begin{pro} There exists a formally smooth $X$-bundle $\tilde{\theta}$ over $\BHomeo^{+}_{L}(X)$ that satisfies the following universal property: Given any simplical complex $K$ and any formally smooth $X$-bundle over $|K|$, it is isomorphic to the pullback $g^{*}(\widetilde{\theta})$ under some map $g:|K|\rightarrow \BHomeo^{+}_{L}(X)$. Furthermore, $g_0^*(\widetilde{\theta})$ and $g_1^*(\widetilde{\theta})$ are concordant if and only if $g_0$ and $g_1$ are homotopic.
\end{pro}
\begin{proof} We first construct $\widetilde{\theta}$. We denote $\Emb(X,\mathbb{R}^{\infty})\times_{\Homeo^{+}(X)} L(\scrT X)$ by $S$. For any $k$-simplex $\sigma=[(f,\rho)]$ of $S$, we consider the bundle $
X\rightarrow E(\sigma)\xrightarrow{\pi(\sigma)} \Delta^{k}$.
Here $E(\sigma)\subset \Delta^{k}\times \mathbb{R}^{\infty}$ is the image of $f$ and $\pi(\sigma)$ is the projection map. The pushforward $f_{*}(\rho)$ gives a linear structure $\rho_{\sigma}$ on $\scrT(E(\sigma)/\Delta^{k})$. Given any morphism $\tau:[l]\rightarrow [k]$, by (\ref{diagram: pullback embedding}), there is a bundle map 
\[\tau_{*}: E(\tau^{*}(\sigma))\rightarrow E(\sigma)
\]
that covers the map $\Delta^{\tau}:\Delta^{l}\rightarrow \Delta^{k}$. We consider the quotient space 
\[
\widetilde{E}:=(\bigsqcup_{k\in \mathbb{N}}\bigsqcup_{\sigma\in S[k]}E(\sigma))/\sim,
\]
where $\sim$ is generated by $\tau_{*}(x)\sim x$ for any morphism $\tau: [l]\rightarrow [k]$, any $\sigma\in S[k]$ and any point $x\in E(\tau^{*}\sigma)$. By gluing $\{\pi(\sigma)\}_{\sigma}$ together, one obtains a topological bundle 
\[
\widetilde{\pi}: E\rightarrow |S|
\]
whose vertical tangent microbundle has a linear structure $\widetilde{\rho}$ obtained by gluing  $\{\rho_{\sigma}\}_{\sigma}$ together. The pair $(\widetilde{\pi},\widetilde{\rho})$ defines $\widetilde{\theta}$. 

Now we prove the universal property of $\widetilde{\theta}$. Consider a formally smooth $X$-bundle $\theta$ over $|K|$, represented by a topological bundle $X\rightarrow E\xrightarrow{\pi} |K|$ equipped with a linear structure $\rho'$ on $\scrT(E/|K|)$. We fix an order $o$ on $0$-simplices of $K$. For each $k$-simplex $\eta$ of $K$, we use $o$ to fix a homeomorphism $\Delta^{k}\cong |\eta|\subset |K|$. And we pick a trivialization $\pi^{-1}(|\eta|)\cong \Delta^{k}\times X$. Then the restriction of $\rho$ gives a linear structure $\rho_{\eta}$ on $\scrT(\Delta^{k}\times X/X)$. Furthermore, by Lemma \ref{lem: Ken} (3), we can inductively construct an embedding $f_{\eta}\in \Emb(X,\mathbb{R}^{\infty})[k]$ for each $\eta$ so that the following compatibility condition holds: for any $l$-simplex $\eta'\subset \eta$, the following diagram commutes: 
\begin{equation}
\label{diagram: compatibility}
\begin{split}
\xymatrix{\pi^{-1}(|\eta'|)\ar@{^{(}->}[d]\ar[r]^{\cong}&\Delta^{l}\times X\ar@{^{(}->}[r]^{f_{\eta'}}& \Delta^{l}\times \mathbb{R}^{\infty}\ar[r]^{\cong}& |\eta'|\times \mathbb{R}^{\infty}\ar@{^{(}->}[d]\\
\pi^{-1}(|\eta|)\ar[r]^{\cong}&\Delta^{k}\times X\ar@{^{(}->}[r]^{f_{\eta}}& \Delta^{k}\times \mathbb{R}^{\infty}\ar[r]^{\cong}& |\eta|\times \mathbb{R}^{\infty}.
}
\end{split}
\end{equation}
Set $\sigma(\eta)=[(f_{\eta},\rho_{\eta})]\in S[k]$. We use $o$ to determine a homeomorphism \[g_{\eta}:|\eta|\xrightarrow{\cong} \Delta^{k}\times \{\sigma_{\eta}\}\subset \Delta^{k}\times S[k].\] By (\ref{diagram: compatibility}), these homeomorphisms can be glued together to form a map $g:|K|\rightarrow |S|$. It is straightforward to check that $g^{*}(\widetilde{\theta})=\theta$. 

Given $g_0,g_1: |K|\rightarrow |S|$, it is obvious that $g_0^{*}(\widetilde{\theta})$ and $g_1^{*}(\widetilde{\theta})$ are concordant if $g_0$ and $g_1$ are homotopic. Furthermore, suppose that $g_0^{*}(\widetilde{\theta})$ and $g_1^{*}(\widetilde{\theta})$ are concordant. By the simplical approximation theorem \cite[Section 5.3, Section 6.9]{Rourke}, we may homotope $g_0$ and $g_1$ so that they are both simplical. By a similar argument as above, we construct a simplical map $I\times |K|\rightarrow |S|$ that restricts to $g_i$ on $\{i\}\times |K|$. This completes the proof.
\end{proof}

We note that, for any $g\in R(TX)[k]$, the exponential map gives a linear structure \[\rho_{g}: T(\Delta^{k}\times X/ \Delta^{k})\rightarrow \scrT(\Delta^{k}\times X/\Delta^k).\] This give a map $R(TX)\rightarrow  L(\scrT X)$ and further induces a map 
\[
\operatorname{Exp}: \Homeo^{+}(X)\times_{\Diff^{+}(X)} R(TX)\rightarrow L(\scrT X)
\]
defined by $\operatorname{Exp}[(f,g)]:=f_{*}(\rho_{g})$ for any $f\in \Homeo^{+}(X)[k]$ and $g\in R(TX)[k]$. 

The following theorem is a reformulation of Kirby--Siebenmann's classification theorem of smooth structures. See \cite[Theorem V.2.2]{Kirby77}.
\begin{thm}[Kirby--Siebenmann]\label{thm: KS}
Suppose that $X$ is 
\textbf{not} $4$-dimensional. Then $\operatorname{Exp}$ is a weak equivalence to the union of path components that intersect the image of $\operatorname{Exp}$. 
\end{thm}
We consider the following forgetful map
\[
\begin{split}
\id\times_{\Diff^{+}(X)} \operatorname{Exp}:\Emb(X,\mathbb{R}^{\infty})\times_{\Homeo^{+}(X)}(\Homeo^{+}(X)\times_{\Diff^{+}(X)} R(TX))&\rightarrow \\
\Emb(X,\mathbb{R}^{\infty})\times_{\Homeo^{+}(X)}&L(\scrT X).
\end{split}
\]
Passing to the geometric realization, we obtain the map
\[
\iota:\BDiff^{+}(X)\rightarrow \BHomeo^{+}_{L}(X). 
\]


\begin{cor}\label{cor: iotawekequivalence}
Suppose that $X$ is not $4$-dimensional. Then the map $\iota$ is a weak equivalence to the path component of $\BHomeo^+_L(X)$ containing the image of $\iota$.
\end{cor}
\begin{proof} This directly follows from Theorem \ref{thm: KS} and \cite[Lemma 2.6 (iii)]{kupers2015}.
\end{proof}

Now we turn to dimension $4$. What we need is a following refinement of \cref{thm: Diff Homeo closed 4-manifold}.

\begin{pro}\label{pro: Diff HomeoL closed 4-manifod} Let $k>0$ and let $X$ be a closed oriented 4-manifold that satisfies the conditions of \cref{thm: Diff Homeo closed 4-manifold}. Consider the forgetful map $\iota: \BDiff^{+}(X_{kS})\rightarrow \BHomeo^{+}_{L}(X_{kS})$.
Then the induced map 
\[
\iota_{*}:H_k(\BDiff^+(X_{kS});\Z)
\to H_k(\BHomeo^+_{L}(X_{kS});\Z)
\]
is not injective and the induced map 
\[
\iota^{*}:H^k(\BHomeo_{L}^+(X_{kS});\Z/2)\rightarrow H^k(\BDiff^+(X_{kS});\Z/2)
\]
is not surjective. 
\end{pro}

Before proving \cref{pro: Diff HomeoL closed 4-manifod}, we shall see the following. Set $I=[0,1]$.

\begin{lem}\label{lem: homeo respect linear structure}
Let $X$ and $X'$ be two simply-connected smooth $4$-manifolds which are homeomorphic to each other. Then there exist  Riemannian metrics $g$ on $X$ and $g'$ on $X'$, smoothly embdded disks $D^{4}\subset X$ and $D^{4'}\subset X'$, and a homeomorphism $f: X\rightarrow X'$ such that the following conditions hold:
\begin{itemize}
    \item $f$ restricts to a smooth isometry from a neighborhood of $D^{4}$ to a neighborhood of $D^{4'}$. 
    \item There is a linear structure $\rho$ on $\scrT(I\times X/I)$ such that $\rho|_{\{0\}\times X}$ equals $\rho_{g}$, $\rho|_{\{1\}\times X}$ equals $f^{*}(\rho_{g'})$, and $\rho|_{I\times D^{4}}$ equals the pullback of $\rho_{g}|_{D^{4}}$.
\end{itemize}
\end{lem}
\begin{proof}  Pick smoothly embedded disks $D^{4}\hookrightarrow X$ and $D^{4'}\hookrightarrow X'$. Let $f: X\rightarrow X$ be a homeomorphism. After an isotopy, we may assume $f(D^{4})\subset D^{4'}$. By the annulus theorem \cite{Quinn82}, $D^{4'}\setminus f(D^{4})$ is homeomorphic to $I\times S^3$. This allows us to further modify $f$ so that it restricts to a diffeomorphism from a neighborhood of $D^{4}$ to  a neighborhood of $D^{4'}$. Then we pick a metric $g'$ on $X'$ and extend the metric $f_{*}(g)$ on $D^{4'
}$ to a metric $g'$ on $X'$. The first condition is now satisfied and it remains to find a linear structure $\rho$ that extends the given linear structure $\rho'$ over the subspace
\[
B=(\{0,1\}\times X)\cup (I\times D^{4})\subset I\times X.
\]
By the Kister--Mazur theorem \cite{Kister64}, the microbundle   $\scrT(I\times X/I)$ is equivalent to a topological bundle $\mathbb{R}^{4}\rightarrow E\rightarrow I\times X$. The given linear structure $\rho'$ on $E|_{B}$ gives a reduction of its structure group from $\operatorname{Top}(4)$ to $O(4)$. (Here $\operatorname{Top}(4)$ denotes the group of homeomorphisms of $\mathbb{R}^{4}$ preserving the origin.) And it suffices to extend this reduction to the whole $I\times X$. There is a series of obstruction classes for the existence of such extension, which sit inside the groups 
\[
H^{i+1}(I\times X, B;\pi_{i}(\operatorname{Top}(4)/O(4)))\]
for $0\leq i\leq 4$. Since $\pi_{i}(\operatorname{Top}(4)/O(4))))\cong 0,0,0,\mathbb{Z}/2,0$ for $i=0,\dots, 4$ respectively (see \cite[Theorem 8.7A]{FreedmanQuinn}), the only possible obstruction is in the group 
\[H^{4}(I\times X, B;\mathbb{Z}/2)\cong H_{1}(X;\mathbb{Z}/2).\]
Since $X$ is simply-connected, this group is trivial as well. So the extension exists and the proof is completed.
\end{proof}
\begin{proof}[Proof of \cref{pro: Diff HomeoL closed 4-manifod}] Let $\alpha\in H_k(\BDiff^+(X_{kS});\Z)$ be the same as in the proof of \cref{thm: Diff Homeo closed 4-manifold}. Then it suffices to show that $\iota_{*}(\alpha)=0$. Recall that $\alpha$ is presented as the difference of two smooth $X_{ks}$-bundles $E^{k}(X)$ and $E^{k}(X')$ over $T^{k}$. Here $X'$ is a smooth 4-manifold homeomorphic to $X$. We just need to show that $E^{k}(X)$ and $E^{k}(X')$ are concordant as formally smooth $X_{ks}$-bundles. 

Note that $E^{k}(X)$ is obtained by gluing together the bundles $E^{k}(D^4)$ and $(X\setminus \mathring{D}^{4})$ along their common boundary. The same operation applied to $X'$ gives the bundle  $E^{k}(X)$. Applying Lemma \ref{lem: homeo respect linear structure}, we obtain the homeomorphism $f:X\rightarrow X'$, a metric $g$ on $X$ and a metric $g'$ on $X'$.  

We have an isomorphism $\tilde{f}:E^{k}(X)\rightarrow  E^{k}(X')$ of topological bundles, obtained by gluing 
\[\operatorname{id}_{T^{k}}\times f|_{X\setminus \mathring{D}^{4}}: T^{k}\times (X\setminus \mathring{D}^{4})\rightarrow  T^{k}\times (X'\setminus \mathring{D}^{4})\]
 with the identity map on $E^{k}(D^4)$. 
 
Next, we fix a families metric $\tilde{g}_0$ on $E^{k}(D^4)$. By  gluing $g$ and $g'$ with $\tilde{g}_0$ fiberwise, we obtain a families metric $\tilde{g}$ on $E^{k}(X)$ and a families metric $\tilde{g}'$ on $E^{k}(X')$. We use $\rho_{\tilde{g}}$ and $\rho_{\tilde{g}'}$ to denote the induced linear structures. By gluing together the linear structure $\rho$ provided by Lemma \ref{lem: homeo respect linear structure} and the linear structure $\scrT(E^{k}(D^4)\times I/T^{k}\times I)$ pulled back from $\tilde{g}_0$, we get a concordance from $\rho_{\tilde{g}}$ and $\tilde{f}^{*}(\rho_{\tilde{g}'})$. This shows that $E^{k}(X)$ and $E^{k}(X')$ are concordant as formally smooth bundles.
\end{proof}

\section{Proof of the main results}\label{subsectionProof of the main instability theorem}

Now we have prepared all necessary results on families Seiberg--Witten theory to prove the main results described in \cref{section Introduction}.
The only remaining piece of the proofs of those is the following result on geography involving the usual (i.e. unparameterized) Seiberg--Witten invariant:
\begin{thm}\label{thm: 4-mfds that dissolves} For any integer $l\in \mathbb{Z}$, there exists a sequence of simply-connected $4$-manifolds $\{X^{l}_{i}\}^{\infty}_{i=1}$ that satisfies the following properties:
\begin{enumerate}
\item $\sign(X^l_i)=l$,
\item $b_2(X^l_i) \to +\infty$ as $i \to +\infty$,
\item $X^{l}_{i}$ is of Seiberg--Witten simple type,
\item $
\SWbbhalftot^{0}(X^{l}_{i})=1\in \mathbb{Z}/2,
$ and
\item $X^l_i$ dissolves after one stabilization.
\item Suppose $l$ is divisible by $16$. Then we can pick $X^l_i$ to be spin or nonspin. 
\end{enumerate}
\end{thm}

To avoid a digression, we postpone the proof of \cref{thm: 4-mfds that dissolves} to \cref{construction of 4-manifolds}.
Let us give the proofs of the main results in \cref{section Introduction} assuming  \cref{thm: 4-mfds that dissolves}.

\begin{proof}[Proof of \cref{thm: main cal}]
We set $l=\sign(X)$ and denote by $X_i$ the 4-manifold $X^l_i$ in \cref{thm: 4-mfds that dissolves}. Here, we choose $X_{i}$ to be of the same type (i.e. spin/non-spin) as $X$.
By Wall's theorem \cite{Wall64}, there exists a positive integer $n(X)$ such that $X_{n(X)S}$ dissolves. There exists an integer $i_0$ such that for any $i>i_0$, the manifold $X_{i}$ is indefinite, not homotopy equivalent to $K3$ and satisfies 
\[
b_{2}(X_{i})+2k\geq b_{2}(X)+2n(X).
\]
For any $i>0$, we set
\[
N_{i}:=\frac{b_2(X_{i_0+i})-b_{2}(X)}{2}+k.
\]
Then we have $N_i \to +\infty$. Furthermore, $(X_{i_0+i})_{kS}$ is diffeomorphic to $X_{N_{i}S}$ since they have the same signature, Euler number, type and both dissolve. We apply \cref{thm key computation source} to the manifold $X_{i_0+i}$ and obtain homology classes 
\[
\alpha^{\delta}_{k}(X_{i_0+i})\in H_{k}(\BDiff_{\del}(\mathring{X}_{N_{i}S})^{\delta};\mathbb{Z})\text{ and }\alpha_{k}(X_{i_0+i})\in H_{k}(\BDiff_{\del}(\mathring{X}_{N_{i}S});\mathbb{Z}).
\]
We set $\alpha^{\delta}_{i}=\alpha^{\delta}_{k}(X_{i_0+i})$ and $\alpha_{i}=\alpha_{k}(X_{i_0+i})$. It remains to verify that they satisfy the desired property. By  \cref{thm key computation source}~(iv), $\alpha^{\delta}_{i}$ and $\alpha_{i}$ are both topologically trivial. 
By \cref{thm key computation source}~(ii), $\alpha^{\delta}_{i}\in \ker s^{\delta}_{N_{i},*}$. This implies that $\alpha_{i}\in \ker s_{N_{i},*}$. 
By \cref{thm key computation source}~(iii), we have 
\begin{align}
\label{eq: non-vanishing pairing proof of main}
 \langle \SWbbhalftot^{k}(X_{N_iS}),\rho_{*}(\alpha_{i})\rangle\neq 0.
\end{align}
This combined with \cref{cor: vanishing} show that \[\alpha_{i}\notin \operatorname{Image}(s_{N_{i}-1,*}\circ\cdots\circ s_{N_{i}-k-1,*}).\]
 This further implies that 
 \[\alpha^{\delta}_{i}\notin \operatorname{Image}(s^{\delta}_{N_{i}-1,*}\circ\cdots\circ s^{\delta}_{N_{i}-k-1,*}),
 \]
which completes the proof.
\end{proof}


\begin{proof}[Proof of \cref{thm: Diff Homeo sequence general}] We let $N_{i}$ be chosen as in the proof of \cref{thm: main cal}.
Then we apply \cref{pro: Diff HomeoL closed 4-manifod} to the manifold $X_{i_0+i}$ to conclude that $\iota_{*}$ is non-injective and $\iota^{*}$ is non-surjective. The proof is finished by the fact that $(X_{i_0+i})_{kS}$ is diffeomorphic to $X_{N_{i}S}$.
\end{proof}

\begin{rmk}
\label{rmk: all ex dissolve}
In \cref{thm: Diff Homeo sequence general}, all $X\#N_i S^2 \times S^2$ are are connected sums of $\CP^2, S^2\times S^2$, $K3$ and their orientation reversals.
This is because we choose $n_i$ in the proof of \cref{thm: Diff Homeo sequence general} so that $X\#n_iS^2\times S^2$ dissolve, and so do $X\#N_iS^2\times S^2$.
\end{rmk}

\begin{proof}[Proof of \cref{thm: abelianisation noninjective}] We set $N_i$ to be the same as in \cref{thm: main cal} and let $\beta_{i}=\rho_{*}(\alpha_i)$. It suffices to show that $\beta_i$ can be represented by an exotic diffeomorphism, since all other assertions directly follow from the corresponding assertions in  \cref{thm: main cal,thm: Diff Homeo sequence general}. 

By the construction of $\alpha_i$ in the proof of \cref{thm: main cal,thm key computation source}, we see that $\beta_{i}$, when treated as a homology class in $H_{1}(\BDiff^{+}(X_{N_i}S);\mathbb{Z})$, can be expressed as the difference 
\[
(X_{N_{i}S})_{f_1,*}[T^1]-(X_{N_{i}S})_{f'_1,*}[T^1]\in H_{1}(\BDiff^{+}(X_{N_i}S);\mathbb{Z}).
\]
Here $(X_{N_{i}S})_{f_1}$ and $(X_{N_{i}S})_{f'_1}$ are mapping tori for diffeomorphisms  $f_1,f_1'\in \Diff^{+}(X_{N_{i}S})$. Furthermore, by \eqref{FreeQuinnPerron}, there is a homeomorphism $g:X_{N_{i}S}\rightarrow  X_{N_{i}S}$ such that $f_1$ is topologically isotopic to $g\circ f'_1\circ g^{-1}$. 
We may assume $g$ is actually a diffeomorphism since the map (\ref{eq: MCG surjective}) is surjective for $N=N_i$. Then \[f\circ g\circ f'^{-1}\circ g^{-1}: X_{N_iS}\rightarrow X_{N_iS}\]
is an exotic diffeomorphism that represents $\beta_i$. 
\end{proof}

\cref{thm: psc main intro} is just a part of  \cref{thm: psc main}, which we prove here:

\begin{proof}[Proof of \cref{thm: psc main}]
Let $\alpha_i \in H_k(\BDiff_{\del}(\mathring{X}\# N_i S^2\times S^2);\Z)$ be the homology class that was constructed in the proof of \cref{thm: main cal}.
Recall that
the pairing of $\alpha_i$ with $\SWbbhalftot^k$ is non-zero, \eqref{eq: non-vanishing pairing proof of main}.
Thus the claim of the \lcnamecref{thm: psc main} follows from \cref{cor: psc non-surj}.
\end{proof}

\begin{rmk}\label{rmk: twisted stablization details} As we mentioned in \cref{rmk: twisted stabilization}, the whole argument can be adapted to twisted stabilizations  (i.e., taking connected sum with $\mathbb{CP}^{2}\#\overline{\mathbb{CP}^2}$ instead of $S^2\times S^2$). We summarize the adaptations needed here. For $k\geq 0$, we use $X_{k\tilde{S}}$ to denote $X\#k(\mathbb{CP}^{2}\#\overline{\mathbb{CP}^2})$. 
For simply-connected nonspin $X$, there is a diffeomorphism between $X_{kS}$ and $X_{k\tilde{S}}$ for any $k$
\cite[Corollary~1]{Wall64D}.
For such $X$, the multiple mapping torus $E^{k}(X)$ is a smooth bundle with fiber $X_{k\tilde{S}}$.  (See \eqref{eq: multiple mapping torus}.)  Note that we do not adopt the construction of $E^{k}(X)$ so the monodromies are still diffeomorphisms supported on copies of on embedded $(S^2\times S^2)\setminus \Int(D^4)$, not embedded $(\mathbb{CP}^{2}\#\overline{\mathbb{CP}^2})\setminus \Int(D^4)$. With this in mind, one can repeat the proof of \cref{thm key computation source} for non-spin $X$ and show that \cref{thm key computation source} remains true if one replaces $X_{kS}$ with $X_{k\widetilde{S}}$ and instead consider the twisted stabilization map
\[
\tilde{s}^{\delta}_{*}:  H_{k}(\BDiff_{\del}(\mathring{X}_{kS})^{\delta};\mathbb{Z})\rightarrow  H_{k}(\BDiff_{\del}(\mathring{X}_{(k+1)S})^{\delta};\mathbb{Z}).\]
Furthermore, by applying \cref{thm: vanishing2} to $N=(k+1)(\mathbb{CP}^{2}\#\overline{\mathbb{CP}^2})$, one obtains a variation of \cref{cor: vanishing} for twisted stabilizations. With these results proved, we can repeat the proofs of \cref{thm: main cal,thm: Diff Homeo sequence general,thm: abelianisation noninjective} word by word for twisted stabilizations on a non-spin $X$. For a spin $X$, we simply replace it with $X_{\tilde{S}}$ and reduce it to the non-spin case.
\end{rmk}






\subsection{Comparison with Torelli groups}
\label{subsection Comparison with Torelli groups}

Given an oriented smooth manifold $X$, let $\TDiff(X)$ denote a subgroup of $\Diff(X)^+$ defined to be the group of diffeomorphisms that act trivially on homology $H_\ast(X;\Z)$.
This subgroup $\TDiff(X)$ is often called the {\it Torelli group}.
Tools from families gauge theory known before have been effectively used to study $\TDiff(X)$ rather than $\Diff^+(X)$ (see, such as, \cite{Rub99}).
In this \lcnamecref{subsection Comparison with Torelli groups}, we remark that $H_\ast(\BTDiff(X))$ and $H_\ast(\Diff^+(X))$ are significantly different, as well as $\pi_1(\Diff^+(X))$ and $H_1(\Diff^+(X))$ are.

For a closed oriented smooth 4-manifold $X$ with a fixed $\spinc$ structure $\fraks$ and a homology orientation $\mathcal{O}$,
let $\Diff(X,\fraks,\mathcal{O})$ be the group of diffeomorphisms that preserve orientation of $X$, $\fraks$ and $\mathcal{O}$.
Note that $\TDiff(X) \subset \Diff(X,\fraks,\mathcal{O})$.
We are going to use a characteristic class 
\begin{align}
\label{eq: SW fixing O}
\SWbb(X,\fraks,\mathcal{O}) \in H^{-d(\fraks)}(\BDiff(X,\fraks,\mathcal{O});\Z)
\end{align}
which was defined in \cite{K21}.
Let $\mathcal{B}(X)$ denote the set of Seiberg--Witten basic classes of $X$ of formal dimension zero.

\begin{pro}
\label{pro: Torelli vs the whole strong}
Let $X$ be a simply-connected closed oriented indefinite smooth 4-manifold with $b^+(X) \geq 2$.
Assume that there are infinitely many smooth oriented 4-manifolds $\{X_i\}_{i=1}^\infty$ such that $X_i \#S^2\times S^2$ is diffeomorphic to $X\#S^2\times S^2$ for every $i$, and $\#\mathcal{B}(X_i) \to +\infty$.
Then the kernels of both of the following natural maps
\begin{align}
H_1(\BTDiff(X\#S^2\times S^2);\Z) &\to H_1(\BDiff^+(X\#S^2\times S^2);\Z),\label{eq: ker H1s'}\\
\pi_1(\BDiff^+(X\#S^2\times S^2)) &\to H_1(\BDiff^+(X\#S^2\times S^2);\Z).\label{eq: ker pi1H1'}
\end{align}
contain subgroups isomorphic to $\Z^\infty=\oplus_{\Z}\Z$.
\end{pro}

\begin{proof}
The proof is an adaptation of Ruberman's argument \cite{Rub99} in the Seiberg--Witten setup with a help of the characteristic class \eqref{eq: SW fixing O}.
Set $Z=X\#S^2\times S^2$.
By assumption, there is a sequence $\{\frakt_i\}_i \subset \Spinc(X,0)$ such that $\frakt_i \in \mathcal{B}(X_i) \setminus \mathcal{B}(X_{i-1})$ for all $i$.
Let $\fraks_i \in \Spinc(Z,1)$ be the connected sum of $\frakt_i$ with the spin structure on $S^2 \times S^2$.
By repeating the construction of the bundle $E \to S^1$ in \cref{thm key computation source} using $(X,X_i)$ in place of $(X,X')$, we obtain a diffeomorphism $f_i \in \TDiff(Z)$ such that:
\begin{itemize}
\item $f_i$ is decomposed into $f_i=g_i \circ g_i'$, where the squares of $g_i, g_i'$ are smoothly isotopic to the identity.
\item $\SW(E_i,\fraks_i) \neq 0$, where $E_i \to S^1$ is the mapping torus of $f_i$ with fiber $Z$.
\item $\SW(E_i,\fraks_j) = 0$ for all $j>i$.
\end{itemize}
As we can see from the proof of  \cref{thm: main cal}, we have that $(E_i)_\ast([S^1]) \in H_1(\BDiff^+(Z);\Z)$ are 2-torsions.
On the other hand, let $E_i^T : S^1 \to \TDiff(X)$ be the classifying map of $E_i$ with structure group $\TDiff(Z)$.
Then, by evaluating $\{(E_i^T)_\ast([S^1])\}_i$ on 
the homomorphism
\[
\bigoplus_{\fraks \in \Spinc(Z,1)} \left<\SWbb(Z,\fraks,\mathcal{O}),-\right>
: H_1(\BTDiff(Z);\Z)
\to \bigoplus_{\fraks \in \Spinc(Z,1)}\Z \cong \Z^\infty,
\]
we have that $\{(E_i^T)_\ast([S^1])\}_i$ are linearly independent in $H_1(\BTDiff(Z);\Z)$.
Thus $\{2(E_i^T)_\ast([S^1])\}_i$ generates a subgroup in the kernel of \eqref{eq: ker H1s'} isomorphic to $\Z^\infty$.
The claim on the map \eqref{eq: ker pi1H1'} immediately follows from the above argument combined with the following commutative diagram, where $h$ denote the Hurewicz maps:
\begin{align*}
\xymatrix{
    \pi_1(\BTDiff(Z))\ar[r]^-{h} \ar[d]^-{\cong} & H_1(\BTDiff(Z);\Z) \ar[r] & H_1(\BDiff^+(Z);\Z) \\
    \pi_0(\TDiff(Z))\ar[r]_-{\subset} & \pi_0(\Diff^+(Z)) \ar[r]^{\cong} &
    \pi_1(\BDiff^+(Z))\ar[u]_-{h}.
   }
\end{align*}
\end{proof}

It is easy to find an example of $X$ and $X_i$ that satisfy the assumptions of \cref{pro: Torelli vs the whole strong} by the knot surgery formula due to Fintushel--Stern \cite{FS98}.

\section{Manifolds that dissolve after one stabilization}\label{construction of 4-manifolds}

In this section, we prove \cref{thm: 4-mfds that dissolves}. This will finish our argument. 
Our proof is inspired by the paper \cite{Gompf91,Hanke03,Akhmedov15}, where the authors constructed many examples of symplectic 4-manifolds that dissolve after one stabilization. The idea is to start with some symplectic 4-manifolds with a unique Seiberg--Witten basic class up to sign (e.g., minimal algebraic surfaces of general type) and do fiber sum along embedded surfaces with self-intersection $0$ and genus $\geq 2$. (We avoid doing fiber sum along tori or spheres because they usually annihilate the half-total Seiberg--Witten invariants.) 

We start with some notations. Let $M_{1}, \cdots ,M_{n}$ be oriented closed smooth 4-manifolds. Let $F$ be an oriented genus $g$ closed surface and let  $F_{i}\subset M_{i}$ be a smoothly embedded surface $F_{i}$ of the same genus and self-intersection $0$.  We let  $D^{2}_{1},\cdots ,D^{2}_{n}$ be disjoint disks in $S^{2}$.
For each $i$, we pick an open tubular neighborhood $\nu(F_{i})$ and an orientation-reversing diffeomorphism 
$$
\varphi_{i}: D^{2}_{i}\times F\rightarrow \overline{\nu(F_{i})}
$$
that covers a diffeomorphism $F\rightarrow F_{i}$. Then the fiber sum is defined as 
$$
M_{1}\#_{F}\cdots \#_{F}M_{n}:=((S^{2} \setminus \bigsqcup_{i} D^{2}_{i})\times F)\bigcup_{\varphi_{i}|_{\partial D^{2}_{i}\times F}} (\bigsqcup_{i} (M_{i}\setminus \nu(F_{i}))). 
$$
Of course, the result depends on the choice of gluing function $\varphi_{i}$. Using $\varphi_{i}$, we also get a boundary parametrization $\partial(M_{i}\setminus \nu(F_{i})) \cong S^{1}\times F$. One has 
\[
\begin{split}
\sign(M_{1}\#_{F}\cdots \#_{F}M_{n})&=\sum^{n}_{i=1}\sign(M_{i}),\\    
\chi(M_{1}\#_{F}\cdots \#_{F}M_{n})&=\sum^{n}_{i=1}\chi(M_{i})+(n-1)(4g-4).
\end{split}
\]
Furthermore, by suitably choosing $\varphi_{i}$, one can always make $M$ not spin. And suppose $M_{i}\setminus \nu(F_{i})$ is spin for all $i$. Then one can also choose  $\varphi_{i}$ to make $M_{1}\#_{F}\cdots \#_{F}M_{n}$ spin. In addition, suppose that, for each $i$, $M_{i}$ is symplectic and $F_{i}$ is a symplectic surface. By a theorem of Gompf \cite{Gompf95}, for any choice of gluing function, the manifold $M_{1}\#_{F}\cdots \#_{F}M_{n}$ always carries a symplectic structure.

\begin{pro}\label{pro: Seiberg--Witten product}
Let $M=M_{1}\#_{F}M_{2}$ be the fiber sum of closed oriented smooth 4-manifolds $M_i$ along a surface $F$ of genus $g\geq 2$. We assume that $b^{+}(M_{i})>1, b_{1}(M_{i})=0$ and $M_{i}$ is of Seiberg--Witten simple-type. Let $\mathfrak{s}_{i}$ be a $\spinc$ structure on $M_{i}\setminus \nu(F_i)$. Let $\mathfrak{s}^{k}$ be the unique $\spinc$ structure on $F\times S^1$ that is pulled back from $F$ and satisfies $\langle c_{1}(\mathfrak{s}^{k}), [F]\rangle =2k$. We consider the Seiberg--Witten invariants
$$
\SW(M,\mathfrak{s}_{1},\mathfrak{s}_{2}):= \sum_{\substack{\{\mathfrak{s}\text{ s.t. } \mathfrak{s}|_{M_{i}\setminus \nu(F_{i})}=\mathfrak{s}_{i},\\ d(\mathfrak{s})=0\}}} \SW(M,\mathfrak{s})
$$
and 
$$
\SW(M_{i},\mathfrak{s}_{i}):= \sum_{\substack{\{\mathfrak{s}\text{ s.t. } \mathfrak{s}|_{M_{i}\setminus \nu(F_{i})}=\mathfrak{s}_{i},\\ d(\mathfrak{s})=0\}}} \SW(M_{i},\mathfrak{s}).
$$ Then we have 
$$
\SW(M,\mathfrak{s}_{1},\mathfrak{s}_{2})= \begin{cases}\pm \SW(M_1,\mathfrak{s}_{1})\cdot \SW(M_2,\mathfrak{s}_{2})&\text{ if }\mathfrak{s}_{1}|_{F\times S^{1}}=\mathfrak{s}_{2}|_{F\times S^{1}}=\mathfrak{s}^{\pm (g-1)},\\
0 &\text{ otherwise.}
\end{cases}
$$
\end{pro}
\begin{proof}
Set $Y=F\times S^1$. We focus on the case $\mathfrak{s}_{1}|_{Y}=\mathfrak{s}_{2}|_{Y}=\mathfrak{s}^{0}$ because the other cases are proved by Mu\~{n}oz--Wang \cite[Theorem 1.2]{MunozWang} \footnote{In the first version of \cite{MunozWang}, the authors also included a proof for the remaining case, by assuming an inequality on the rank of the monopole Floer homology of $F\times S^1$ for the torsion $\spinc$ structure. By the computation of Jabuka--Mark \cite[Theorem 9.4]{Jabuka08} on Heegaard-Floer homology, this inequality is now known to hold.} and Morgan--Szab\'o \cite[Corollary 9.10]{Morgan96}. The argument here is adapted from \cite{MunozWang} and \cite{Jabuka08} (which proved an analogous result for Ozsv\'ath--Szab\'o's mix invariants). 
We start by considering the graded abelian group
\[
X(g,g-1):=\oplus_{i=0}^{g-1}\Lambda^{i}H^{1}(F;\mathbb{Z})\otimes\frac{U^{i-g+1}\cdot \mathbb{Z}[U]}{U\cdot \mathbb{Z}[U]}.
\]
Here the degree of $U$ equals $-2$. We use $X(g,g-1)_{k}$ to denote the degree-$k$ component of $X(g,g-1)$.
Let $\eta$ be the 1-cycle in $Y$ given by $*\times S^{1}$. Then one can use $\eta$ to define a local system $\Gamma_{\eta}$ over the configuration space for $Y$. The fiber of $\Gamma_{\eta}$ equals the field of formal Laurent series with rational coefficients,
$$
\mathcal{L}:=\mathbb{Q}[[t^{-1},t].
$$
Under this local system, one has a canonical isomorphism between the two versions of monopole Floer homology groups
$$
\widecheck{HM}_{*}(Y,\mathfrak{s}^0;\Gamma_{\eta})\cong \widehat{HM}_{*}(Y,\mathfrak{s}^0;\Gamma_{\eta})
$$
and we just denote it by $\mathcal{H}$.
The Heegaard Floer homology of $Y$ is computed by Jabuka--Mark \cite[Theorem 9.4]{Jabuka08}.\footnote{In \cite{Jabuka08}, the authors did the completion in the other direction and used the field $\mathbb{Q}[t^{-1},t]]$ instead. We follow the convention of \cite{KM07} and do negative completion.} We can translate their result to monopole Floer homology using the work of Kutluhan--Lee--Taubes \cite{Kutluhan20}. (See \cite[Theorem 3.1]{Lee19} for a version for twisted coefficients) or alternatively, the work of Colin--Honda--Gighini \cite{Colin2011} and Taubes \cite{Taubes10}. This gives us an isomorphism of vector spaces  over $\mathcal{L}$ with relative $\mathbb{Z}$-grading
\begin{equation}\label{eq: HM for surface times a circle}
\mathcal{H}\cong X(g,g-1)\otimes \mathcal{L}.
\end{equation}
We use $\mathcal{H}_{k}$ to denote the component corresponding to  $X(g,g-1)_{k}\otimes \mathcal{L}$. Then $\mathcal{H}_{0}\cong \mathcal{H}_{2g-2}\cong \mathcal{L}$.
The orientation-reversing diffeomorphism on $Y$ given by the reflection in the $S^{1}$-component induces a nondegenerate  pairing 
\[
\langle\cdot,\cdot\rangle: \mathcal{H}\otimes_{\mathcal{L}}\mathcal{H}\rightarrow \mathcal{L}.
\]
Elements in $\mathcal{H}_{k}$ and $\mathcal{H}_{l}$ can have nonzero pairing only if $k+l=2g-2$. 
The vector space $\mathcal{H}$ is also a graded module over the ring 
$$
\mathbb{A}(F)=\Lambda^{*}H_{1}(F;\mathbb{Z})\otimes \mathbb{Z}[U].
$$
The action of $U$ is just the multiplication. For any $\gamma\in H_{1}(F;\mathbb{Z})$ and any $\alpha\in \Lambda^{i}H^{1}(F;\mathbb{Z})$, the action is given by 
\[
\gamma\cdot (\alpha\otimes U^{l})=\iota_{\gamma}\alpha\otimes U^{l}+ (\operatorname{PD}(\gamma)\wedge \alpha)\otimes U^{l+1}+\text{terms with negative $t$-degree.}
\]
Here $\iota$ denotes the contraction. In particular, we see that for any element in $\alpha\in \mathcal{H}\setminus \mathcal{H}_{0}$, there exists an homogeneous $\tau\in \mathbb{A}(F)$ with positive degree such that $\tau\cdot \alpha$ is a nonzero element in $\mathcal{H}_0$.
To prove the claim, we consider the manifold $N=F\times D^{2}$ with boundary $Y$ and the relative 2-cycle $\nu_{N}=*\times D^2$ bounded $\eta$. Let $\mathfrak{s}_{N}$ be the unique $\spinc$ that extends $\mathfrak{s}^0$. Then the relative Seiberg--Witten invariant of $N$ gives a map that preserves the structure of $\mathbb{Z}$-graded $\mathbb{A}(F)$-modules
\[
\phi^{SW}_{N,\nu_{N}}(\mathfrak{s}_{N},-):\mathbb{A}(F)\rightarrow \mathcal{H}.
\]
We claim that
\begin{equation}\label{eq: relative invariant for D2 cross F}
\phi^{SW}_{N,\nu_{N}}(\mathfrak{s}_{N},1)\in \mathcal{H}_{2g-2}\setminus \{0\}.    
\end{equation}
To see this, we consider the manifold $\widetilde{N}=N\cup_{Y}N=F\times S^2$. We let $\mathfrak{s}_{\widetilde{N}}$ be unique torsion $\spinc$ structure on $\widetilde{N}$ and let $\nu_{\widetilde{N}}$ be the 2-cycle $*\times S^2$. We use $\operatorname{PD}[F]\in H^{2}(\widetilde{N};\mathbb{Z})$ to fix a component of the positive cone 
$$
\{\alpha\in H^{2}(\widetilde{N};\mathbb{R})\mid \alpha\cdot \alpha\geq 0, \alpha\neq 0\}.
$$
This allows us to define the Seiberg--Witten invariant for the negative chamber 
\[
\SW_{-}(\widetilde{N},\mathfrak{s},-):\mathbb{A}(F)\rightarrow \mathbb{Z}
\]
for any $\spinc$ structure $\mathfrak{s}$ on $\widetilde{N}$. By applying the gluing formula for Seiberg--Witten invariants \cite[Propposition 27.5.1]{KM07}, one obtains the equality
\[
\sum_{k\in \mathbb{Z}}\SW_{-}(\widetilde{N},\mathfrak{s}_{\widetilde{N}}+n\operatorname{PD}[F],U^{g-1})t^{2n}=\langle \phi^{SW}_{N,\nu_{N}}(\mathfrak{s}_{N},1),U^{g-1}\cdot \phi^{SW}_{N,\nu_{N}}(\mathfrak{s}_{N},1)\rangle.
\]
Using the positive scalar curvature metric on $\widetilde{N}$, we obtain that 
\[
\SW_{-}(\widetilde{N},\mathfrak{s}_{\widetilde{N}}+n\operatorname{PD}[F],U^{g-1})=\SW_{+}(\widetilde{N},\mathfrak{s}_{\widetilde{N}}-n\operatorname{PD}[F],U^{g-1})=0
\]
for any $n\geq 0$. Using the wall-crossing formula (see \cite[Theorem 1.2]{Li95} and \cite[Theorem 16]{Okonek96}), we further obtain that
\[
\SW_{-}(\widetilde{N},\mathfrak{s}_{\widetilde{N}}+n\operatorname{PD}[F],U^{g-1})=\begin{cases}
0 \quad&\text{if }n\geq 0,\\
(-n)^{g} \quad&\text{if }n\leq -1.
\end{cases}
\]
In particular, this implies $U^{g-1}\cdot \phi^{SW}_{N,\nu_{N}}(\mathfrak{s}_{N},1)\neq 0$. Hence the homogeneous element $\phi^{SW}_{N,\nu_{N}}(\mathfrak{s}_{N},1)$ must be a nonzero element in $\mathcal{H}_{2g-2}$. This proves the claim (\ref{eq: relative invariant for D2 cross F}).
For $i=1,2$, we let $\nu_{i}$ be a relative 2-cycle in $M_{i}\setminus \nu(F)$ bounded by $\eta$ and consider the relative Seiberg--Witten invariant 
\[
\phi^{SW}_{M_{i}\setminus \nu(F)}(\mathfrak{s}_{i},-): \mathbb{A}(F)\rightarrow \mathcal{H}.
\]
Then the gluing theorem applied to $M_{i}$ gives the relation
\[
\sum_{\{\mathfrak{s}\mid \mathfrak{s}|_{M_{i}\setminus\nu(F)}\}}\SW(M_{i},\mathfrak{s},\alpha)t^{\langle c_{1}(\mathfrak{s}),\nu_{i}\cup \nu_{N}\rangle}=\langle \alpha\cdot \phi^{SW}_{M_{i}\setminus \nu(F)}(\mathfrak{s}_{i},1), \phi^{SW}_{N,\nu_{N}}(\mathfrak{s}_{N},1))\rangle
\]
for any $\alpha\in \mathbb{A}(F)$. We actually have $\phi^{SW}_{M_{i}\setminus \nu(F)}(\mathfrak{s}_{i},1)\in \mathcal{H}_0$. Otherwise we can find some $\alpha$ with positive degree such that \[\alpha\cdot \phi^{SW}_{M_{i}\setminus \nu(F)}(\mathfrak{s}_{i},1)\in \mathcal{H}_0.\]
Then by (\ref{eq: relative invariant for D2 cross F}), we have 
\[
\langle \alpha\cdot \phi^{SW}_{M_{i}\setminus \nu(F)}(\mathfrak{s}_{i},1), \phi^{SW}_{N,\nu_{N}}(\mathfrak{s}_{N},1))\rangle\neq 0.
\]
That implies $\SW(M_{i},\mathfrak{s},\alpha)\neq 0$ for some $\mathfrak{s}$, which contradicts our assumption that $M_{i}$ is of Seiberg--Witten simple type. 
Since any two elements in $\mathcal{H}_0$ pairs trivially, we get
\[
\langle \phi^{SW}_{M_{1}\setminus \nu(F)}(\mathfrak{s}_{1},1),\phi^{SW}_{M_{2}\setminus \nu(F)}(\mathfrak{s}_{2},1)\rangle=0.
\]
As the last step, we apply the gluing formula to the manifold $M$ and obtain that 
\[
\begin{split}
\sum_{\substack{\{\mathfrak{s}\text{ s.t. } \mathfrak{s}|_{M_{i}\setminus \nu(F_{i})}=\mathfrak{s}_{i},\\ d(\mathfrak{s})=0\}}}\SW(M,\mathfrak{s})\cdot t^{\langle c_{1}(\mathfrak{s}),\nu_1\cup \nu_2\rangle}&=\langle \phi^{SW}_{M_{1}\setminus \nu(F)}(\mathfrak{s}_{1},1),\phi^{SW}_{M_{2}\setminus \nu(F)}(\mathfrak{s}_{2},1)\rangle\\
&=0\in \mathcal{L}.
\end{split}
\]
Comparing the coefficients of both sides, we finish the proof.
\end{proof}


\begin{pro}\label{pro: fiber sum has nontrivial SW-tot} For $2\leq n$ and $1\leq i\leq n$, let $M_{i}$ be a closed symplectic 4-manifold with $b^{+}(M_{i})>1$. Let $F_{i}\subset M_{i}$ be an embedded symplectic surface of genus $g\geq 2$ and self-intersection $0$. Suppose that $M_{i}$ has a unique Seiberg--Witten basic class up to sign and suppose $H_{1}(M_{i}\setminus \nu(F_{i});\mathbb{R})=0$. Take the fiber sum $M=M_{1}\#_{F}\cdots  \#_{F}M_{n}$.
 Then we have  
 \begin{equation}\label{eq: SW-halftot of fiber sum}
 \SWbbhalftot^{0}(M)=1.    
 \end{equation}
\end{pro}
\begin{proof}
 By the work of Taubes, the canonical $\spinc$ structure $\mathfrak{s}_{J,i}$ on $M_{i}$ and its conjugate both have Seiberg--Witten invariant $\pm 1$.
 By our assumptions, all the other $\spinc$ structures on $M_{i}$ have trivial Seiberg--Witten invariants.
 Note that  $\mathfrak{s}_{J,i}$ is not torsion since the  adjunction formula $\langle c_{1}(\mathfrak{s}_{J,i}),[F_{i}]\rangle =2g-2$
holds.

 Let $\mathfrak{s}_{i}$ be a $\spinc$ structure on $M_{i}\setminus \nu(F_{i})$. Consider the  quantity
$$
\SW(M,\mathfrak{s}_{1},\cdots,\mathfrak{s}_{n}):= \sum_{\substack{\{\mathfrak{s}\text{ s.t. } \mathfrak{s}|_{M_{i}\setminus \nu(F_{i})}=\mathfrak{s}_{i},\\ d(\mathfrak{s})=0\}}} \SW(M,\mathfrak{s}).
$$
Note that the fiber sum of $M_{i}$ carries a symplectic structure so is of Seiberg--Witten simple-type. Via a repeated application of Proposition \ref{pro: Seiberg--Witten product}, we see that 
\begin{equation}\label{eq: SW fiber sum}
\SW(M,\mathfrak{s}_{1},\cdots,\mathfrak{s}_{n})=\begin{cases}
\pm 1 &\text{ if }\mathfrak{s}_{i}=\mathfrak{s}_{J,i}|_{M_i\setminus\nu(F_{i})} \text{ for all }i,\\
\pm 1 &\text{ if }\mathfrak{s}_{i}|_{M_{i}}=\overline{\mathfrak{s}_{J,i}}|_{M_i\setminus\nu(F_{i})} \text{ for all }i,\\
0 &\text{ otherwise.}
\end{cases}
\end{equation}
To deduce (\ref{eq: SW-halftot of fiber sum}) from (\ref{eq: SW fiber sum}), it suffices to show that any $\spinc$ structure $\mathfrak{s}$ that satisfies 
$$
\mathfrak{s}|_{M\setminus \nu(i)}=\overline{ \mathfrak{s}|_{M\setminus \nu(i)}} \text{ for any }i
$$
does not contribute to $\SWbbhalftot^{0}(M)$. Indeed, for such $\mathfrak{s}$, we have 
$$
c_{1}(\mathfrak{s})\in \ker (H^{2}(M;\mathbb{R})\rightarrow \oplus^{n}_{i=1} H^{2}(M_{i}\setminus \nu(F_{i});\mathbb{R})).
$$
This implies that $\operatorname{PD}(c_{1}(\mathfrak{s}))$ can be represented by a cycle in $(S^{2} \setminus \bigsqcup_{i} D^{2}_{i})\times F$. Therefore, $c^{2}_{1}(\mathfrak{s})=0$ and we have 
$$
d(\mathfrak{s})=-\frac{2\chi(M)+3\sign(M)}{4}=-\sum^{n}_{i=1}\frac{2\chi(M_{i})+3\sign(M_{i})}{4}-(n-1)(2g-2).
$$
By Taubes's result \cite{Taubes00}, we have $2\chi(M_{i})+3\sign(M_{i})\geq 0$. So $d(\mathfrak{s})< 0$ and it does not contribute to $\SWbbhalftot^{0}(M)$.
\end{proof}

As the next step, we will construct three symplectic manifolds $M_{1}, M_{2}, M_{3}$, which serve the building blocks in our construction of the manifold $X^{l}_{i}$.

The manifold $M_{1}$ is constructed by Persson--Peters--Xiao \cite{Persson}. We sketch its construction as follows: Let $Y_{0}=\mathbb{CP}^{1}\times \mathbb{CP}^{1}$ and let $\pi: Y_0\rightarrow \mathbb{CP}^{1}$ be the projection to the second factor.  Consider a triple sequence of double covers 
$$
Y_3\xrightarrow{p_{3}}Y_2 \xrightarrow{p_2}Y_1\xrightarrow{p_1} Y_0
$$
with branched loci being $B'_3\subset Y_2, p^{-1}_{1}(B_2)\subset Y_1$ and $B_1\subset Y_0$ respectively. Here $B_1, B_2$ are singular curves in $Y_0$ and $B'_3$ is a singular curve in $Y_2$ that is linearly equivalent to $p^{-1}_{2}p^{-1}_{1}(B_3)$ for some $B_3\in Y_0$. The algebraic surface $Y_{3}$ has certain singularities coming from the singular points of the branched loci. And $M_{1}$ is defined to be the resolution of $Y_3$. By suitable choosing $B_1, B_2, B_3$ and $B'_3$, the manifold $M_1$ can be made into a  simply-connected algebraic surface of general type with $\sign(M_1)>0$. 
Suppose the bi-degree of $B_{i}$ equals $(a_{i},b_{i})$. Then one can ensure that $M_{1}$ is spin by choosing $B_{i}$ so that both $a_{1}+a_{2}+a_{3}$ and $b_{1}+b_{2}+b_{3}$ are divisible by $4$. Furthermore, the map $p: M_1\rightarrow \mathbb{CP}^{1}$ defined by the composition
$$
M_1\rightarrow Y_3\xrightarrow{p_{3}}Y_2 \xrightarrow{p_2}Y_1\xrightarrow{p_1} Y_0\xrightarrow{\pi} \mathbb{CP}^{1}
$$
makes $M_1$ into a fibration. We let $F_{1}$ be a generic fiber. Then $F_{1}$ is a $\mathbb{Z}/2\oplus \mathbb{Z}/2\oplus \mathbb{Z}/2$ branched cover of $\mathbb{CP}^{2}$. 
Using the Riemann--Hurewicz formula, we see that 
\[
1<g(F_1)=2(b_{1}+b_{2}+b_{3})-7\equiv 1\mod 8.
\]
In \cite{Persson}, the authors proved that $M_1$ is simply-connected. The same argument implies the simply-connectedness of $M_{1}\setminus F_{1}$. 
Indeed, using \cite[Corollary F]{Persson}, the authors found a disk $D^2\subset \mathbb{CP}^{1}$ such that $p^{-1}(D^2)$ is simply-connected. Since the fibration $p: M_1\rightarrow D^2$ has no multiple fibers, by the same argument as \cite[Lemma 3.20]{Persson81}, this implies 
$
\pi_{1}(M_1\setminus F_{1})=1.
$

The manifold $M_2$ is constructed in \cite[Proof of Corollary D]{Persson} in a similar flavor. For any positive integers $a,b$, we let $Y_{a,b}$ be the double cover of $Y_0$ branched over  
$$(\mathbb{CP}^{1}\times \{2a \text{ generic points}\})\cup (\{2b \text{ generic points}\}\times \mathbb{CP}^{1}).$$
We let $M'_2$ be the fiber product 
$$
\xymatrix{
M'_2 \ar[d] \ar[r] & Y_{a_{4},b_{4}}\ar[d]\\
Y_{a_{5},b_{5}} \ar[r] & Y_0.}
$$
Then $M'_{2}$ has discrete singular points coming from double points of the branched loci. We let $M_2$ be the minimal resolution of these singularities. Topologically, that means we remove tubular neighborhoods of these singular points (which are cones over $\mathbb{RP}^{3}$) and glue copies of the disk bundle for the $TS^{2}$.  It is proved in \cite{Persson} that $M_2$ is a simply-connected surface of general type, with signature 
$$
\sign(M_2)=-8(a_4b_4+a_5b_5)<0.
$$
Furthermore, $M_2$ is spin if $a_4+a_5,b_4+b_5$ are even. 

As before, we consider the composition 
$$
M_2\rightarrow Y_0\xrightarrow{\pi} \mathbb{CP}^{1}
$$
and let $F_2$ be a generic fiber. Then $F_2$ is a $\mathbb{Z}/2\oplus \mathbb{Z}/2$ cover of $\mathbb{CP}^{1}$ branched at  $2(b_4+b_5)$ points, each with two preimages. From this, we compute that $g(F_2)=2(b_4+b_5)-3$. We set 
$$
b_{4}=\frac{g(F_1)+7}{4},\quad  b_{5}=\frac{g(F_1)-1}{4}.
$$
Then $g(F_2)=g(F_1)$. Note that 
\[
\frac{\sign(M_{2})}{16}=-a_{4}-\frac{b_{5}(a_{4}+a_{5})}{2}.
\]
By picking suitable $a_4,a_5$, we may assume that  $\frac{\sign(M_1)}{16}, \frac{\sign(M_2)}{16}$ are coprime to each other. Since $M_2$ is a minimal algebraic surface of general type, it has a unique Seiberg--Witten basic class up to sign.
One can use a similar argument as $M_1$ to show that $\pi_{1}(M_{2}\setminus F_{2})=1$.

The manifold $M_3$ is constructed by Fintushel--Park--Stern \cite{Fintushel2002}. Unlike $M_1, M_2$, the manifold $M_{3}$ is not complex. In \cite{Fintushel2002}, the authors constructed many symplectic manifolds with one Seiberg--Witten basic class up to sign. In particular, for any positive integer $p$, they constructed a simply-connected symplectic $4$-manifold $X_{p}$ by taking the fiber sum of two rational surfaces $R(2p-3)$ and $S(p)$ along a symplectic surface of genus $2p-5$. And one has 
$$\sign(X_{p})=25-14p.$$ We let $M_3=X_{p}$ with $p=\frac{g(F_{2})+5}{2}$ and let $F_3$ be the surface on which the fiber sum is conducted. It is proved in \cite{Fintushel2002} that $\pi_{1}(M_3\setminus F_3)=1$. Furthermore, one can find an embedded surface in $S_{p}\setminus F_{3}\subset X_{p}\setminus F_{3}$ that has an odd self-intersection number. 

So far, we have constructed the building blocks $(M_{i},F_{i})$ and established some of their important properties. We summarize these properties in the following lemma.

\begin{lem}\label{lem: building block}
There exist closed symplectic 4-manifolds $M_1, M_{2}, M_{3}$ and embedded symplectic surfaces $F_{i}\subset M_{i}$ of self-intersection $0$ such that the following conditions are satisfied:
\begin{enumerate}
\item $g(F_{1})=g(F_2)=g(F_3)\geq 2$.
\item $M_{i}$ has a unique Seiberg--Witten basic class up to conjugacy.
    \item $M_{1}, M_{2}$ are spin, and $M_3$ is not spin.
    \item $\pi_{1}(M_i\setminus F_i)=1$ for $i=1,2,3$.
    \item $\sign(M_1)>0$ while $\sign(M_{2}), \sign(M_3)<0$.
    \item $\frac{\sign(M_{1})}{16}$ and $\frac{\sign(M_{2})}{16}$ are coprime to each other.
    \item $\sign(M_3)$ is odd.
    \item $M_{3}\setminus F_{3}$ contains an embedded surface with odd self-intersection number.
\end{enumerate}
\end{lem}

The last ingredient we need is the following proposition

\begin{pro}[{Mandelbaum \cite{Mandelbaum79}, see also \cite[Corollary 5]{Gompf91}}
]\label{pro: fiber sum stabilization} Let $M, M'$ be smooth oriented closed 4-manifolds. Let $F\subset M$, $F'\subset M'$ be closed, connected, orientable surfaces with self-intersection $0$ and of the same genus $g$. 
Suppose that $M$ is spin and suppose $\pi_{1}(M')=\pi_{1}(M
\setminus F)=1$. Then we have 
$$
(M\#_{F} M')\#(S^{2}\times S^{2})
\cong M\# M'\# 2g(S^2\times S^2)
$$
if $M\#_{F} M'$ is spin, and 
$$
(M\#_{F} M')\#(S^{2}\times S^{2})
\cong M\# M'\# 2g(\mathbb{CP}^{2}\# \overline{\mathbb{CP}^2})
$$
if $M\#_{F} M'$ is not spin.
\end{pro}
\begin{proof}[Proof of Theorem \ref{thm: 4-mfds that dissolves}] For $a,b,c\geq 0$, we use $M(a,b,c)$ to denote the fiber sum of $a$ copies of $M_1$, $b$ copies of $M_2$ and $c$ copies of $M_{2}$ along the embedded surfaces $F_{1}, F_{2}, F_{3}$. Then $M(a,b,c)$ is always simply-connected. By Gompf's result \cite{Gompf95}, $M(a,b,c)$ is symplectic for any $a,b,c$ and can be chosen to be spin when $c=0$. When $c\neq 0$, $M(a,b,c)$ is not spin by by Lemma \ref{lem: building block} (8).

By Wall's result \cite{Wall64}, there exists $n\gg 0$ such that $M_{1}, M_{2}, M_{3}$ all dissolves after $n$ stabilizations.

Now we construct nonspin $X_{l}^{i}$ for any $l\in \mathbb{Z}$. By Lemma \ref{lem: building block} (5), (6), (7), there exists integers $a,b$ and a positive integer $c$ such that $$
a\cdot\sign(M_{1})+b\cdot\sign(M_2)+c\cdot\sign(M_{3})=l.
$$
Then we set  
$$
X^{l}_{i}= M(a-(n+i)\sign(M_2), b+(n+i)\sign(M_{1}), c).
$$
By Proposition \ref{pro: fiber sum has nontrivial SW-tot}, $\SWbbhalftot^0(X^{l}_{i})=1$. By Proposition \ref{pro: fiber sum stabilization}, $X^{l}_{i}$ dissolves after one stabilization. (Again note that for nonspin, simply-connected manifolds, connected summing with $S^2\times S^2$ is the same as connected summing with $\mathbb{CP}^{2}\# \overline{\mathbb{CP}^2}$.)

Next, we construct spin $X^{l}_{i}$ assuming  
 $16|l$. In this case, we can find $a,b$ such that 
 $$
a\cdot\sign(M_{1})+b\cdot\sign(M_2)=l.
$$
We set 
$$
X^{l}_{i}= M(a-(n+i)\sign(M_2), b+(n+i)\sign(M_{1}), 0).
$$
Again by Proposition \ref{pro: fiber sum stabilization}, $X^{l}_{i}$ dissolves after one stabilization.
\end{proof}


\bibliographystyle{plain}
\bibliography{mainref}

\end{document}